\documentclass{amsart}

\usepackage{amsfonts,amssymb,verbatim,amsmath,amsthm,latexsym,textcomp,amscd}
\usepackage{latexsym,amsfonts,amssymb,epsfig,verbatim}
\usepackage{amsmath,amsthm,amssymb,latexsym,graphics,textcomp, hyperref}
\usepackage{graphicx}
\usepackage{color}
\usepackage{url}

\begin{document}

\newtheorem{theorem}{Theorem}[section]
\newtheorem{prop}[theorem]{Proposition}
\newtheorem{lemma}[theorem]{Lemma}
\newtheorem{cor}[theorem]{Corollary}
\newtheorem{definition}[theorem]{Definition}
\newtheorem{defn}[theorem]{Definition}
\newtheorem{conj}[theorem]{Conjecture}
\newtheorem{claim}[theorem]{Claim}
\newtheorem{defth}[theorem]{Definition-Theorem}
\newtheorem{obs}[theorem]{Observation}
\newtheorem{rmk}[theorem]{Remark}
\newtheorem{qn}[theorem]{Question}

\newcommand{\hhat}{\widehat}
\newcommand{\boundary}{\partial}
\newcommand{\C}{{\mathbb C}}
\newcommand{\integers}{{\mathbb Z}}
\newcommand{\natls}{{\mathbb N}}
\newcommand{\ratls}{{\mathbb Q}}
\newcommand{\reals}{{\mathbb R}}
\newcommand{\proj}{{\mathbb P}}
\newcommand{\lhp}{{\mathbb L}}
\newcommand{\tube}{{\mathbb T}}
\newcommand{\cusp}{{\mathbb P}}
\newcommand\AAA{{\mathcal A}}
\newcommand\BB{{\mathcal B}}
\newcommand\CC{{\mathcal C}}
\newcommand\DD{{\mathcal D}}
\newcommand\EE{{\mathcal E}}
\newcommand\FF{{\mathcal F}}
\newcommand\GG{{\mathcal G}}
\newcommand\HH{{\mathcal H}}
\newcommand\II{{\mathcal I}}
\newcommand\JJ{{\mathcal J}}
\newcommand\KK{{\mathcal K}}
\newcommand\LL{{\mathcal L}}
\newcommand\MM{{\mathcal M}}
\newcommand\NN{{\mathcal N}}
\newcommand\OO{{\mathcal O}}
\newcommand\PP{{\mathcal P}}
\newcommand\QQ{{\mathcal Q}}
\newcommand\RR{{\mathcal R}}
\newcommand\SSS{{\mathcal S}}
\newcommand\TT{{\mathcal T}}
\newcommand\UU{{\mathcal U}}
\newcommand\VV{{\mathcal V}}
\newcommand\WW{{\mathcal W}}
\newcommand\XX{{\mathcal X}}
\newcommand\YY{{\mathcal Y}}
\newcommand\ZZ{{\mathcal Z}}
\newcommand\CH{{\CC\Hyp}}
\newcommand\MF{{\MM\FF}}
\newcommand\PMF{{\PP\kern-2pt\MM\FF}}
\newcommand\ML{{\MM\LL}}
\newcommand\PML{{\PP\kern-2pt\MM\LL}}
\newcommand\GL{{\GG\LL}}
\newcommand\Pol{{\mathcal P}}
\newcommand\half{{\textstyle{\frac12}}}
\newcommand\Half{{\frac12}}
\newcommand\Mod{\operatorname{Mod}}
\newcommand\Area{\operatorname{Area}}
\newcommand\ep{\epsilon}
\newcommand\Hypat{\widehat}
\newcommand\Proj{{\mathbf P}}
\newcommand\U{{\mathbf U}}
 \newcommand\Hyp{{\mathbf H}}
\newcommand\D{{\mathbf D}}
\newcommand\Z{{\mathbb Z}}
\newcommand\R{{\mathbb R}}
\newcommand\Q{{\mathbb Q}}
\newcommand\E{{\mathbb E}}
\newcommand\til{\widetilde}
\newcommand\length{\operatorname{length}}
\newcommand\tr{\operatorname{tr}}
\newcommand\gesim{\succ}
\newcommand\lesim{\prec}
\newcommand\simle{\lesim}
\newcommand\simge{\gesim}
\newcommand{\simmult}{\asymp}
\newcommand{\simadd}{\mathrel{\overset{\text{\tiny $+$}}{\sim}}}
\newcommand{\ssm}{\setminus}
\newcommand{\diam}{\operatorname{diam}}
\newcommand{\pair}[1]{\langle #1\rangle}
\newcommand{\T}{{\mathbf T}}
\newcommand{\inj}{\operatorname{inj}}
\newcommand{\collar}{\operatorname{\mathbf{collar}}}
\newcommand{\bcollar}{\operatorname{\overline{\mathbf{collar}}}}
\newcommand{\I}{{\mathbf I}}

\newcommand{\bbar}{\overline}
\newcommand{\UML}{\operatorname{\UU\MM\LL}}
\newcommand{\EL}{\mathcal{EL}}
\newcommand\MT{{\mathbb T}}
\newcommand\Teich{{\mathcal T}}

\title{Cannon-Thurston Maps for Surface Groups}

\author[Mahan Mj]{Mahan Mj}

\address{RKM Vivekananda University, Belur Math, WB-711 202, India}

\email{mahan.mj@gmail.com; mahan@rkmvu.ac.in}


\thanks{Research partially supported by  DST Research grant
DyNo. 100/IFD/8347/2008-2009. }   

\date{}

 \begin{abstract}

We prove the existence of Cannon-Thurston maps for simply and doubly degenerate surface Kleinian groups. As a consequence we prove that
  connected
limit sets of finitely generated Kleinian groups are locally connected.
 
\end{abstract}

\maketitle

\tableofcontents

\section{Introduction}  
 Let $\Gamma$ be a finitely generated Kleinian group, i.e. a 
 finitely generated discrete subgroup of $Isom (\Hyp^3) (= PSl_2(C))$, the isometry group of hyperbolic 3-space. Then $\Gamma$ acts on the boundary Riemann sphere $S^2$
(of $\Hyp^3$)
by conformal automorphisms. The limit set of $\Gamma$, denoted by
$\Lambda_\Gamma$, is the collection of accumulation points of any $\Gamma$-orbit in $S^2$. The limit set is independent of the $\Gamma$-orbit
chosen. In particular,
for any $z \in  \Lambda_\Gamma$, the orbit $\Gamma .z$ is dense in $\Lambda_\Gamma$.  The complement $S^2 \setminus \Lambda_\Gamma$ is called the domain of discontinuity
of $\Gamma$ and is denoted $D_\Gamma$. The action of $\Gamma$ on $D_\Gamma$ is properly discontinuous. 
Thus, the limit set  $\Lambda_\Gamma$ may be thought of as the locus of chaotic dynamics
for the action of $\Gamma$ on $S^2$ and it would  be desirable to `tame' it. 

\smallskip

\noindent {\bf Motivation and Statement of Results:}\\
Towards this,
  Thurston raises the following  question (see \cite[Problem 14]{thurstonbams}):

\begin{qn} \label{thurston-bams-qn}        Suppose $\Gamma$ has the property that $(\Hyp^3 \cup D_\Gamma)/{\Gamma}$ is compact. Then is it
true that the limit set of any other Kleinian group $\Gamma^{\prime}$ isomorphic to $\Gamma$ is the
continuous image of the limit set of $\Gamma$, by a continuous map taking the
fixed points of an element $\gamma$ to the fixed points of the corresponding element $\gamma^\prime$? \end{qn}

Essentially the same question is raised by Cannon and Thurston
in Section 6 of \cite{CT, CTpub} in the specific context of surface Kleinian groups:

\begin{qn}
Suppose that a surface group $\pi_1 (S)$ acts freely and properly
discontinuously on ${\mathbb{H}}^3$ by isometries such that the quotient manifold has no accidental parabolics. Does the inclusion
$\tilde{i} : \widetilde{S} \rightarrow {\mathbb{H}}^3$ extend
continuously to the boundary?
\label{ctqn}
\end{qn}

The authors of \cite{CT} point out that for a simply degenerate surface Kleinian group,
this is equivalent, via the Caratheodory extension Theorem, to asking if the limit set is locally connected. The most general question in this context
is the following:

\begin{qn} Let $\Gamma$ be a finitely generated Kleinian group such that the limit set $\Lambda_\Gamma$ is connected. Is 
$\Lambda_\Gamma$ locally connected?
\label{lcqn}
\end{qn}
 
It is a classical fact of general topology that a compact, connected, locally connected metric space $X$ is homeomorphic to a Peano continuum, i.e.
$X$ is a continuous image of the closed interval $[0,1]$. Hence, asking if the limit set is locally connected is equivalent to asking if 
there is some parametrization by $[0,1]$.  Question \ref{thurston-bams-qn} makes this precise by asking for an explicit parametrization.
For surface Kleinian groups, Question \ref{ctqn} asks for a parametrization of $\Lambda_\Gamma$ by a circle.
 In this paper, we give a  positive answer to  Question \ref{ctqn}.  

\smallskip

\noindent {\bf Theorems \ref{crucial} and \ref{crucial-punct}:} {\it
Let $\rho$ be a representation of a surface group $H (= \pi_1(S))$  into
$PSl_2(C)$ without accidental parabolics. Let $M$ denote the (convex
core of) ${\mathbb{H}}^3 / \rho
(H)$.  Further suppose that $i: S \rightarrow M$, taking
 parabolic to parabolics, induces a homotopy
equivalence.
  Then the inclusion
  $\tilde{i} : \widetilde{S} \rightarrow \widetilde{M}$ of universal covers extends continuously
  to a map $\hat{i} : \widehat{S} \rightarrow \widehat{M}$ between the compactifications
  of universal covers. Hence the limit set
  of $\rho (H)$ is locally connected.}

\smallskip

In \cite{mahan-kl} we  extend the techniques of this paper to answer Question \ref{thurston-bams-qn} affirmatively.
The continuous boundary extensions above are called Cannon-Thurston maps. 

\smallskip

Combining Theorems \ref{crucial} and \ref{crucial-punct} with a theorem of Anderson and Maskit \cite{and-mask}, we have the following
affirmative answer to Question \ref{lcqn}:

\smallskip

\noindent {\bf Theorem \ref{lcfinal}:}
{\it Let $\Gamma$ be a finitely generated Kleinian group with connected limit set $\Lambda$. Then $\Lambda$ is locally connected.}

\smallskip

Note that the limit set of a finitely generated Kleinian group $\Gamma$ is connected if and only if the boundary 
of the {\it convex core} of ${\mathbb{H}}^3/\Gamma$ is incompressible away from cusps.

\smallskip

\noindent {\bf Relationship with The Ending Lamination Theorem:}\\
Seminal work of Minsky \cite{minsky-elc1} and Brock-Canary-Minsky \cite{minsky-elc2}, building on work of Masur-Minsky \cite{masur-minsky, masur-minsky2},
has resolved Thurston's Ending Lamination Conjecture. The Ending Lamination Theorem roughly says that for a simply or doubly degenerate surface Kleinian group
$\Gamma$
without accidental parabolics, the isometry type of the manifold $M = {\Hyp}^3/\Gamma$ is determined by its end-invariants. For a doubly degenerate group,
the end-invariants are two ending laminations, one each for the two geometrically infinite ends of $M$.
For a simply  degenerate group,
the end-invariants  are an ending lamination  corresponding to the  geometrically infinite end of $M$ and a conformal structure corresponding to the 
geometrically finite end of $M$. The ending lamination corresponding to a  geometrically infinite end may be regarded as a purely topological piece of data
associated to the end. Thus, in the context of geometrically infinite Kleinian groups,
 the Ending Lamination Theorem roughly says that {\it `Topology implies Geometry'}: an analog of Mostow Rigidity
 for infinite covolume Kleinian groups.

Theorems \ref{crucial} and \ref{crucial-punct} prove the existence of Cannon-Thurston maps for surface Kleinian groups, but leave unanswered the question
about the point preimages of these maps. In \cite{mahan-elct}, we relate the point preimages of Cannon-Thurston maps for simply and doubly
degenerate surface Kleinian groups to ending laminations. In particular, the ending lamination corresponding to a degenerate end can be recovered
from the Cannon-Thurston map. 
More generally, since topological conjugacies are compatible with Cannon-Thurston maps,
a topological conjugacy of $\Gamma -$ actions on limit sets comes from a
biLipschitz homeomorphism of quotient manifolds.
Hence  the Ending Lamination Theorem \cite{minsky-elc1, minsky-elc2} in conjunction with 
Theorems \ref{crucial} and \ref{crucial-punct} and the main result 
of \cite{mahan-elct} show that the geometry of $M$ can be recovered from the action of $\Gamma$ on the limit set $\Lambda_\Gamma$.  This justifies the slogan: {\it `Dynamics on the Limit Set determines Geometry in the Interior.'}\\

\smallskip

\noindent {\bf History:} \\
 Several authors have  contributed to the theme of this paper. We shall give
below a brief account of the history
of the problem along with some further developments that use the results of this paper.

Cannon and Thurston \cite{CTpub},
Minsky \cite{minsky-jams}, Alperin, Dicks and Porti \cite{ADP}, Cannon
and Dicks \cite{cd1, cd2},  Klarreich \cite{klarreich}, McMullen \cite{ctm-locconn}, Bowditch \cite{bowditch-stacks, bowditch-ct} and the author \cite{mitra-trees, mitra-ct, mahan-pared, mahan-ibdd, mahan-bddgeo, mahan-amalgeo} 
have obtained partial positive answers to Questions \ref{thurston-bams-qn} and \ref{ctqn}. We describe some of this history in brief.

In \cite{abikoff-ct}, Abikoff gave an approach to a negative answer to Question \ref{lcqn}. However, around 1980, Thurston realized that
this approach would not work.
 Then, in a foundational paper, Cannon and Thurston \cite{CT} gave the first examples furnishing a positive answer to Question \ref{thurston-bams-qn}
for geometrically infinite surface Kleinian groups; hence the term `Cannon-Thurston map'.
In approximate chronological order, the existence of Cannon-Thurston maps in the context of Kleinian groups was proven \\
\begin{enumerate} 
\item by Floyd \cite{Floyd} for geometrically finite Kleinian groups.
\item  by Cannon and Thurston \cite{CT} \cite{CTpub} for fibers of closed hyperbolic 3 manifolds fibering over the circle
and for simply degenerate groups with asymptotically periodic ends. 
\item by Minsky \cite{minsky-jams} for closed surface groups of bounded geometry (see also \cite{mitra-trees, mahan-bddgeo}).
\item by the author \cite{mitra-trees}, and independently by Klarreich \cite{klarreich} using different methods, for hyperbolic 3-manifolds
of bounded geometry  with
 an incompressible core and  without parabolics.
\item  by Alperin-Dicks-Porti \cite{ADP} for fibers of the figure eight knot complement regarded as a fiber bundle over the circle. 
 \item 
 by McMullen \cite{ctm-locconn} for punctured torus groups (see also \cite{mahan-ibdd}).
 \item  by Bowditch \cite{bowditch-stacks, bowditch-ct} for punctured surface groups of bounded geometry (see also \cite{mahan-pared}).
\item by Miyachi \cite{miyachi-ct} for handlebody groups of bounded geometry (see also \cite{souto-ct}).
 \item by the author \cite{mahan-pared} for  hyperbolic 3-manifolds
of bounded geometry  with core
 incompressible away from cusps.
 \item by the author \cite{mahan-ibdd, mahan-amalgeo} for special unbounded geometries.\\
\end{enumerate}

\noindent {\bf Further Developments:} \\
In 
\cite{mahan-elct}, we give an explicit parametrization of the limit set of a surface Kleinian group by describing the point pre-images of the Cannon-Thurston map
and relating them to ending laminations.
In 
a further followup paper \cite{mahan-kl}, we answer Question \ref{thurston-bams-qn} affirmatively and completely for all finitely generated Kleinian groups,
using some preliminary work in \cite{mahan-red}.
The techniques of this paper can thus be strengthened  to show that Cannon-Thurston maps exist in general for finitely generated Kleinian groups, 
thus answering a conjecture of McMullen \cite{ctm-locconn}. 

\smallskip

\noindent {\bf Acknowledgments:} I would like  to thank Jeff Brock,
Dick Canary and Yair
Minsky
for their help
during the course of this work. In particular,  Minsky and
 Canary brought a couple of critical gaps in  previous versions of this paper to my notice. I would also like to thank Benson Farb for innumerable exciting
conversations on relative hyperbolicity when we were graduate students.  I would like to thank 
Caroline Series  and Al Marden for a number of comments and corrections and Kingshook Biswas for comments
that have improved the exposition. Special thanks are due to the referee
for  not unwelcome pressure to re-organize a couple of unwieldy manuscripts (\cite{mahan-amalgeo} and \cite{mahan-split})
into a relatively more streamlined version and also for painstakingly going through all the details and making several
important and constructive suggestions. 
A major revision of this paper was done when the author was visiting Universite Paris-Sud (Orsay) under the ARCUS Indo-French program.

\smallskip

\noindent {\bf Dedication:} This paper is fondly dedicated to Gadai
and  Sarada for their support and indulgence.

\subsection{Broad Scheme of Proof} Let $M$ be a hyperbolic 3-manifold homotopy equivalent to a closed hyperbolic
surface $S$. We think of $S$ as an embedded incompressible surface in $M$.
Let $\til{S}$ and $ \til{M} (= \Hyp^3)$ denote the  universal covers of $S, M$ 
respectively and 
  $\tilde{i} : \widetilde{S} \rightarrow \widetilde{M}$ the  inclusion of universal covers.

Given a hyperbolic geodesic segment $\lambda$ in $\til S$ lying outside 
a large ball about a fixed reference point $o \in \til{S}$, our aim is to show that the geodesic in $\Hyp^3$
joining the endpoints of $\til{i} (\lambda)$  lies outside a large ball 
about $\til{i} (o)$
in $\Hyp^3$. This is sufficient to prove the existence of Cannon-Thurston maps
(Lemma \ref{contlemma}). Instead of proving this directly, our objective will be to construct a  set  $\LL_\lambda$ (called a `ladder') 
containing $\til{i} (\lambda)$  such that

\begin{enumerate} 
\item If $\lambda$ lies outside a large ball in $\til S$, then the ladder $\LL_\lambda$ lies outside a a large ball in $ \til{M}$.
It is much easier to show (and follows from an essentially elementary argument) that $\LL_\lambda$ lies outside a a large ball than to find the
exact (or even approximate) location of the geodesic in $\Hyp^3$
joining the endpoints of $\til{i} (\lambda)$. Hence this approach.
\item $\LL_\lambda$ is quasiconvex  with respect to a modified (pseudo) metric $d_G$ on  $\til{M}$, thus forcing the
$d_G$ geodesic joining the endpoints of $\lambda$ to lie $d_G-$close to $\LL_\lambda$.
\item The pseudometric $d_G$ is constructed in such a way that $\LL_\lambda$ still controls 
(cf. Lemma \ref{ea-strong}) the location of the hyperbolic geodesic $\beta^h$ in $\Hyp^3$ 
joining the endpoints of $\til{i} (\lambda)$, thus forcing $\beta^h$ to lie outside a large ball in $\Hyp^3$.
\end{enumerate}

For ease of notation we shall often identify  any point or subset of $\til S$ with its image under $\til i$.

\subsubsection{The Ladder:}  One of the main steps in proving the  sufficient condition of Lemma \ref{contlemma}
(and hence concluding the existence of Cannon-Thurston maps) is 
to  construct a quasiconvex `hyperbolic ladder' as in \cite{mitra-trees} and \cite{mitra-ct} containing $\lambda$.
Suppose that   a sequence $\{ S_i \}$  of  {\it disjoint, equispaced, embedded, bounded geometry} surfaces exiting an end $E$ of $M$
has been `judiciously' constructed. We shall describe a little later what `judicious' means.
We think of  $\{ S_i \}$ as a sequence of surfaces exiting a vertical end $E$.
Identify $S$ with the base 
surface $S_0$. 

Choose a basepoint  in $S$ and fix  a lift  $p$ of the base-point in $\til S$ as the origin. Let ${r}$ be a quasigeodesic
ray in $M$, starting at  $p$,
 exiting  $E$ and making linear progress as it exits $E$. Suppose $\lambda = [a,b] \subset \til{S}$ is a geodesic
in the {\it intrinsic} metric on $\til S$ joining two lifts $a (=a_0)$ and $b(=b_0)$ of $p$. Let $r_a, r_b$ be the lifts
of $r$ starting at $a, b$ respectively. Let $a_i $ (resp. $b_i$) be the point at which $r_a$ (resp. $ r_b$) intersects
$\til{ S_i}$. Let $\lambda_i$ be the geodesic in the  intrinsic metric on $\til{S_i}$ joining $a_i, b_i$. The ladder associated to
the sequence  $\{ S_i \}$  and the geodesic $\lambda$ is
$\LL_\lambda = \bigcup_i \lambda_i$.  To prove quasiconvexity
of $\LL_\lambda$ we construct a  retraction $\Pi_\lambda$ of $\bigcup_i \til{ S_i}$ onto $\LL_\lambda = \bigcup_i \lambda_i$ by defining 
$\Pi_\lambda$ on $ \til{ S_i}$ as the nearest point retraction onto $\lambda_i$ in the {\it intrinsic} metric on $ \til{ S_i}$.
We would like to ensure that $\Pi_\lambda$ is coarsely Lipschitz. 
The construction of $\LL_\lambda$ and $\Pi_\lambda$ is detailed in Section \ref{ladder}. 

The ladder $\LL_\lambda$ has the following property that
we want: If $\lambda$  lies outside 
a large ball about the origin in $\til S$, then  $\LL_\lambda$  lies outside 
a large ball about the origin in $\til M$.
 
This construction works exactly for 3-manifolds of bounded geometry, where the $S_i$'s may be chosen such that
\begin{enumerate}
\item Equispaced Condition:  the regions between $S_i$ and $S_{i+1}$ are uniformly biLipschitz to $S_i \times [0,1]$ (for all $i$). 
\item Quasi-isometry Condition:   The map from $\til{S_i}$ to $\til{S_{i+1}}$ that takes $(x,i)$ to $(x, i+1)$ is a uniform quasi-isometry.
\end{enumerate}

Both of these break down in general. In fact,
 quasiconvexity of $\LL_\lambda$  is not in general true in the
hyperbolic metric on $\til M$ for the choice of the sequence $\{ S_i \}$ we describe below. 

The technical tool we shall use  to address this issue in this  paper is {\it electric geometry}
and relative hyperbolicity (Section \ref{relh}). Let  $\HH = \{ H_i\}$ be a collection of quasiconvex
subsets of $\Hyp^3$. The electric (pseudo) metric obtained by electrocuting elements of $\HH$
 essentially allows one to travel for free within any $H_i$. However, this metric has the crucial feature that
electric geodesics control hyperbolic geodesics (Lemma \ref{ea-strong}) and hence allows 
recovery
of  hyperbolic geodesics from electric geodesics. We emphasize that it is quasiconvexity of  $H_i$'s that allows this recovery.

\subsubsection{A motivational special case of Split Geometry:} We describe first a special  case of the model geometry  of a geometrically infinite
unbounded geometry end $E$. This will be a particular case of what is referred to as `split geometry' later on in the paper, and is representative in
a sense to be explicated. (The model geometry described here was called `graph amalgamation geometry' in \cite{mahan-amalgeo}).
Suppose we have  the following situation: 

\begin{enumerate}
\item there exists a sequence
$\{ S_i \}$  of   disjoint,  embedded, bounded geometry surfaces exiting  $E$. These are ordered in a natural way along $E$, i.e. $i<j$ implies that
$S_j$ is contained in the unbounded component of $E \setminus S_i$. The topological product region  between $S_i$ and $S_{i+1}$ is denoted $B_i$.
\item corresponding to each such product region $B_i$, there exists a Margulis tube $T_i$ such that $T_i \subset B_i$. Further, 
 $T_i \cap S_i$ and  $T_i \cap S_{i+1}$ are annuli on $S_i$ and $S_{i+1}$ respectively, with core curves homotopic to the  core curve of $T_i$. 
\end{enumerate}

We think of the Margulis tube $T_i$ as `splitting' the block $B_i$ and hence the surfaces $S_i$ and $S_{i+1}$.  
(See figure below.) The complementary components $K_{ij}$ 
of $B_i \setminus T_i$  and their lifts $\til{K_{ij}}$ to $\til E$ will play a special role later.

\begin{center}

\includegraphics[height=6cm]{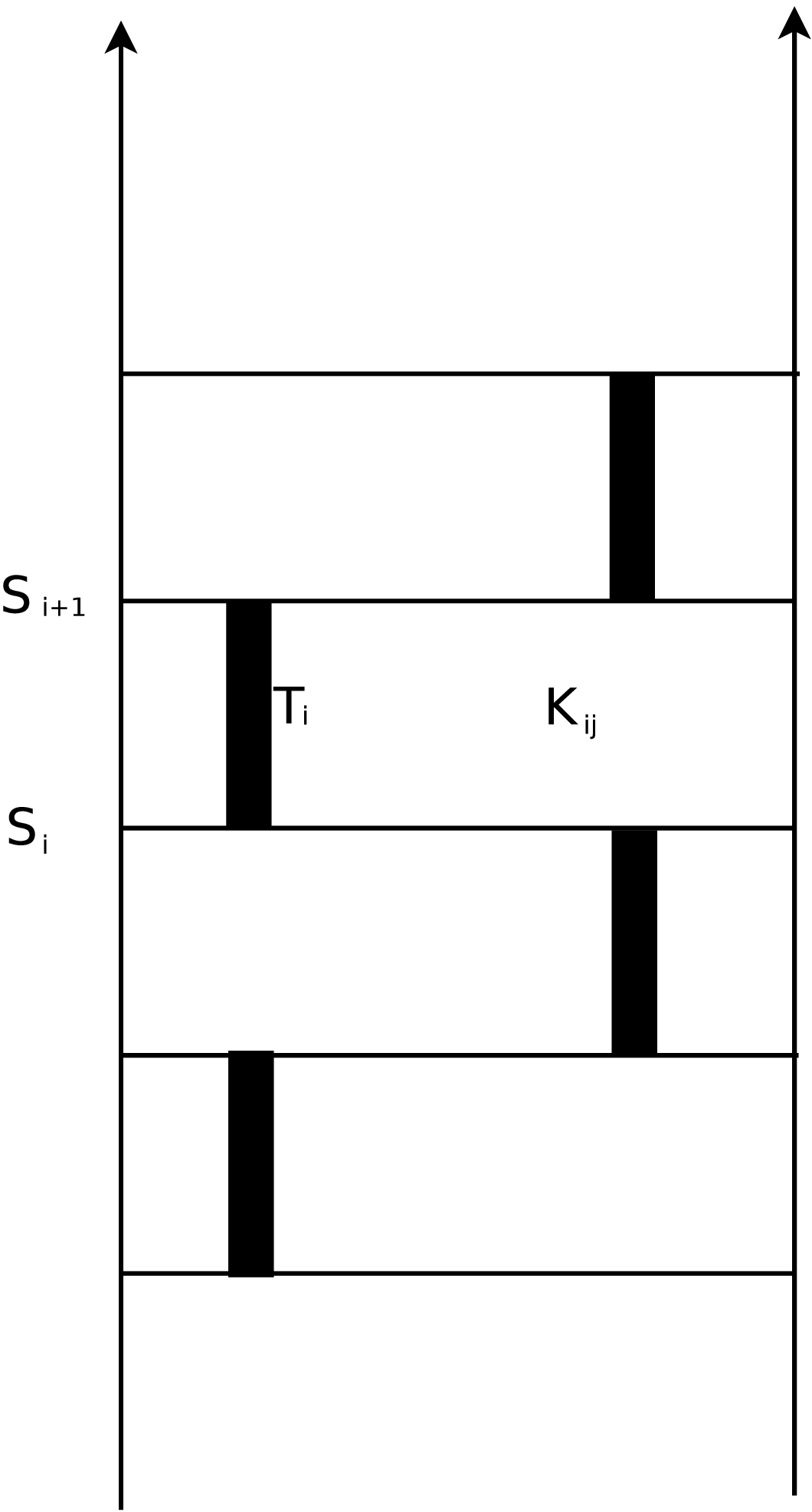}

\smallskip

\underline{Figure:  {\it A special case of split geometry} }

\end{center}

Note that we have no control on the geometry of the 
complementary components $K_{ij}$. So the only thing we can do with them is to electrocute them and lose
the geometry contained {\it within} any such component. Electrocuting $K_{ij}$'s forces $S_i$ and $S_{i+1}$ to be equispaced (about distance one apart from each other).
It is in this modified electric metric that the sequence $\{ S_i \}$ satisfies both the Equispaced Condition
and the Quasi-isometry Condition above and the ladder construction can go through.

 It will turn out that the universal covers $\til{K_{ij}} \subset \til{E}$ are quasiconvex in a certain weak sense.
Thus, we can electrocute such components and still hope to recover hyperbolic geodesics from electric geodesics using Lemma \ref{ea-strong}.

\subsubsection{Choice of the sequence of  surfaces:}  We shall first describe a couple of restrictive assumptions
on  a degenerate end that reduce it to the above model geometry.
We shall then state (very briefly) how one needs to modify the above model to obtain a model geometry for a general degenerate end.\\
{\bf The Special Case:} We  give a brief sketch of the  simplifying assumptions on a general degenerate end
that leads us to a model geometry and a choice of a sequence $\{ S_i \}$ as above.
First one needs a linear order on incompressible (but not necessarily embedded) surfaces in $E$. It is at this stage that we need Minsky's
model manifold from \cite{minsky-elc1}, and more generally the hierarchy machinery from \cite{masur-minsky2}.
 The model manifold of \cite{minsky-elc1} does not quite furnish a sequence of complete surfaces exiting $E$, but rather
a sequence of pants decompositions of $S$ exiting $E$.  
A sequence $\{ P_m \}$ of pants decompositions of $S$ exiting $E$ means the following. Fix an isometry type $\mathbb{P}$ of a pair of pants.
Let $\tau_m$ denote the  simple multicurve on $S$ forming the boundary curves of
the pants decomposition $P_m$.  Then  the complement of a thin (open) annular neighborhood of $\tau_m$ in $S$ can be identified with $P_m$.
We demand that this complementary region (identified with $P_m$) can be embedded in $E$ such that $P_m$ with
the inherited metric is of uniformly bounded geometry, i.e.
each pair of pants (component)
in $P_m$ is uniformly biLipschitz to $\mathbb{P}$.
We demand further that each such embedding of $P_m$ can be extended to a topological embedding $S_m$ of $S$ and these topological embeddings 
$\{ S_m \}$ exit the end $E$.

 The sequence $\{ P_m \}$ of (boundary curves of) pants decompositions exiting $E$ is often referred to
as a  {\it resolution}. The sequence $\{ P_m \}$ in \cite{minsky-elc1} is chosen in such a way that
the boundary curves of the pants decompositions $\{ P_m \}$ occurring in the resolution have short geodesic realization in $E$.

Each pants decomposition gives a simplex in
the curve complex $CC(S)$ of the surface $S$. Hence
the resolution furnishes a special kind of a path of simplices in   $CC(S)$.
Associated to such a path is a 
{\it geodesic} of simplices in   $CC(S)$ called  a {\it tight geodesic} \cite{minsky-elc1}.  A tight geodesic furnishes a `tight sequence' 
$\cdots , \tau_i, \tau_{i+1}, \cdots$
of multicurves on the surface $S$. (Note the difference between the suffixes $i$ and $m$ at this stage.
This indicates that we are actually passing to a subsequence.) This material is detailed in Section \ref{min}.

\smallskip

\noindent {\it Simplifying Assumption 1:}
Assume, for simplicity, that for all $i$, the length 
of exactly one curve in $\tau_i$ is sufficiently small, less than the
Margulis constant in particular.

\smallskip

Call the short curve $\tau_i$ for convenience.
The surface $S_i$  corresponds (roughly) to the first occurrence of the vertex $\tau_i$ in the resolution. Since  $\tau_i$
is short, the Margulis tube $T_i$ corresponding to it {\it splits} both $S_i$ and $S_{i+1}$. 

\smallskip

\noindent {\it Simplifying Assumption 2:}
Assume further that 
that the surfaces $S_i$  have injectivity radius uniformly bounded below, i.e. the tube $T_i$ is
trapped entirely  between $S_i$ and $S_{i+1}$.

\smallskip

The product region $B_i$ between 
$S_i$ and $S_{i+1}$ will be called a {\it split block} as it is split by $T_i$. 
This situation (an end $E$ satisfying  Simplifying Assumptions 1 and 2) gives us the model geometry (special case of split geometry) described above.\\

\noindent {\bf The General Case:} 
The construction of the sequence  $\{ S_i \}$  in general (without the simplifying assumptions of the special case) is described in detail in Section 
\ref{splitgeometry}. Here we content ourselves by providing a couple of caveats.

Note first that $B_i \setminus T_i$ might be very far
from a metric product. Thus electrocution is a necessity to make the $S_i$'s equispaced.

We point out further that
in general (when Simplifying Assumption 2 is no longer valid) the Margulis tube $T_i$ may not be entirely contained in $B_i$, but may extend into $B_{i+1}$ or $B_{i-1}$.
As a result the  surface $S_i$ may have a thin part contained entirely in $T_i$, destroying the product structure of $B_i$. 

To address this issue, we shall excise the interiors of Margulis tubes and  `weld' the `vertical sides' (see previous figure)
of $T_i$ together. The resulting manifold is called  the {\it welded model manifold} $M_{wel}$. $M_{wel}$ is thus a quotient space of $M$ homeomorphic  to $M$ itself.
In the previous schematic figure, the thick dark vertical rectangle denotes a section of the Margulis tube $T_i$.
The quotient map identifies the vertical sides of   this vertical rectangle
and collapses the horizontal $I-$direction  to a  point ($T_i$ should be thought of as a product of the dark vertical rectangle with a circle).
We shall
also construct a new (pseudo) metric $d_{tel}$  on $B_i$ after welding the vertical sides (the `welded blocks').
This process is called {\it tube electrocution} and is carried out
on the  welded model manifold $M_{wel}$
rather than the model manifold itself in Section \ref{electcns}. The pseudometric
$d_{tel}$ on the welded manifold $M_{wel}$
roughly gives zero length to all horizontal circles of $T_i$ and a uniformly bounded length  to the vertical direction.

\subsubsection{ Split Geometry and Graph Quasiconvexity:} 
Lifts $\til K$ of components of $B_i \setminus T_i$ to the universal cover $\til{ M_{wel}}$ are called split components. We construct an auxiliary
metric $d_G$ called the graph (pseudo) metric on $\til{ M_{wel}}$ by {\it electrocuting} the family  of split components in $\til{ M_{wel}}$.
What this means is that for each split component $\til{K} \subset \til{ M_{wel}}$ we attach a copy of $\til{K} \times [0,\frac{1}{2}]$, identifying
$\til{K} \times \{ 0 \}$ with $\til K  \subset \til{ M_{wel}}$ and equipping $\til{K} \times \{ \frac{1}{2} \}$ with the zero metric.
(This is  slightly different from Farb's coning construction \cite{farb-relhyp}.) A crucial fact we prove in
Sections \ref{sec-gqc} and \ref{sec-gqc1} is that the hyperbolic convex hull $CH(\til{K})$ has uniformly bounded diameter
in the graph metric $d_G$. We describe this by
saying that $\til K$ is uniformly {\it graph-quasiconvex} as any hyperbolic geodesic joining points in 
 $\til K$ lies in a uniformly bounded neighborhood of  $\til K$ in the $d_G$-metric.
It follows that $(\til{ M_{wel}}, d_G)$ is a (Gromov)-hyperbolic metric space. Equivalently,
$\til{ M_{wel}}$  is {\it weakly hyperbolic} relative to the collection of split components. Note that we cannot in general use strong relative
hyperbolicity as two adjacent split components $\til{K_1}, \til{K_2} \subset\til{ M_{wel}}$ intersect along a lift of a welded Margulis tube.
This issue is responsible for much of the strife in the Recovery Step  below (Sections \ref{recovery1} and \ref{recovery2}). 

Gromov-hyperbolicity of $(\til{ M_{wel}}, d_G)$ ensures quasiconvexity of the ladder  $\LL_\lambda$ in $(\til{ M_{wel}}, d_G)$ whose construction is described
above.
This is proven in  Section \ref{ladder}.

\subsubsection{ Recovery of Hyperbolic Geodesics:} 
There is a fair bit of technical difficulty at this stage. The graph metric is constructed on the
welded model manifold. So we have to have a way of getting back to the model manifold from the 
welded model manifold. To do this, we note that the complement of the Margulis tubes in the model
manifold and the complement of the welded tubes in the welded model
manifold are the same. This allows us to construct a pseudometric quasi-isometric to $d_G$ on the
model manifold $M$ itself. Abusing notation slightly, we call this pseudometric $d_G$ also.

The split components of $\til M$ are obtained from those of $\til{ M_{wel}}$ by adjoining certain Margulis tubes.
Weak relative hyperbolicity of $\til M$  relative to the collection of split components
gives us control over hyperbolic geodesics in terms of geodesics in $(\til{ M_{wel}}, d_G)$. The process of recovering a hyperbolic
geodesic from a geodesic in $(\til{ M_{wel}}, d_G)$ is detailed in Sections \ref{recovery1} and \ref{recovery2}. A more detailed sketch of the
scheme of recovery is given in Section \ref{scheme}.

\subsubsection{A Flowchart of Main Ideas:} Here is a mnemonic flow-chart of the above scheme that may be useful:
\begin{itemize}
\item $\til{M} \longrightarrow \til{M_{wel}}$  (welding)
\item  $\til{M_{wel}} \longrightarrow (\til{M_{wel}},d_{tel})$ (tube-electrocution)
\item   $(\til{M_{wel}},d_{tel}) \longrightarrow (\til{M_{wel}},d_{G})$ (split-component-electrocution)
 \item   $(\til{M_{wel}},d_{G}) \longrightarrow (\til{M},d_{G}) \longrightarrow \til{M}$ (recovery)
\end{itemize}

The principal purpose behind carrying out each of these steps is given below in brief:

\begin{itemize}
\item Welding allows us to construct a sequence of {\em bounded geometry} surfaces  exiting the end(s) of $M_{wel}$, though such a sequence might not exist in $M$.
The  sequence of  bounded geometry surfaces permits us to construct the ladder $\LL_\lambda$ in $\til{M_{wel}}$.
\item Tube electrocution and split-component-electrocution ensure both the Equispaced Condition and the Quasi-isometry condition. In a certain sense therefore,
the two electrocution steps allow us to reduce the problem to a model satisfying Conditions (1) and (2) of the bounded geometry case. We can (as in the 
bounded geometry case) show that
$\LL_\lambda$ is quasiconvex in $(\til{M_{wel}},d_{G})$. 
\item Quasiconvexity of $\LL_\lambda$ furnishes a $d_G-$quasigeodesic in $\til{M_{wel}}$ contained in $\LL_\lambda$ joining the end-points of $\lambda$.
\item Finally, the recovery step allows us to come back from $(\til{M_{wel}},d_{G})$ to $\til M$ via $(\til{M},d_{G})$.
\end{itemize}

\subsubsection{ Outline of the paper:}
We recall the notions of relative hyperbolicity and electric geometry (cf. \cite{farb-relhyp}) in Section 2 and derive some consequences
that will be useful in this paper. In Section 3, we collect together features of the model manifold constructed by Minsky in \cite{minsky-elc1} and proven to be a biLipschitz model
for simply and doubly degenerate manifolds by Brock-Canary-Minsky in \cite{minsky-elc2}. In Section 4, we select out a sequence of split surfaces
from the split surfaces occurring in the model manifold and proceed to `fill' the intermediate spaces between successive split surfaces by special blocks homeomorphic to $S \times I$.
This gives us a `split geometry' model for simply and doubly degenerate manifolds. We make crucial use of electric geometry and relative hyperbolicity
at this stage. In Section 5, we construct a quasiconvex (Gromov) `hyperbolic ladder' in the (Gromov)
hyperbolic electric space constructed in Section 4
and use it to construct a  quasigeodesic in the electric metric joining the endpoints of $\lambda$. In Section 6, we recover
information about the hyperbolic geodesic joining the endpoints of $\lambda$ from the electric geodesic constructed in Section 5. In Section 7 we put all
the ingredients together to prove the existence of Cannon-Thurston maps for closed surface Kleinian groups (Theorem \ref{crucial}). In Section 8
we describe the modifications necessary for punctured surfaces.

\subsubsection{ Notation:} We shall in general use $N$ (resp. $N^h$) to denote (the convex core of) a simply or doubly degenerate hyperbolic 3-manifold {\it without 
(resp. with) cusps}. For a manifold $N^h$ with cusps, $N$ will also denote $N^h$ minus an open neighborhood of the cusps. $M$ will denote the model manifold
(Section \ref{min}).

Similarly  $S$ (resp. $S^h$) shall denote a closed (resp. finite volume with cusps) hyperbolic surface. 
 For a surface $S^h$ with cusps, $S$ will also denote $S^h$ minus an open neighborhood of the cusps. We shall sometimes use $S$ to denote a biLipschitz
homeomorphic image of a hyperbolic $S$.
Thus $M, N$ will both be homeomorphic to $S \times J$ where $J = [0, \infty )$ or $\mathbb R$ according as $N$ is simply or doubly degenerate.

Since we shall not have specific need for manifolds with cusps till the last Section of this paper, $N$ will denote (the convex core of) a simply or doubly degenerate hyperbolic 3-manifold
 without 
 cusps unless otherwise mentioned.

$d$ will denote the hyperbolic (or biLipschitz to hyperbolic) metric on $S$. $d_M$ will denote the metric on the model manifold.

\subsection{Gromov Hyperbolic Metric Spaces and Cannon-Thurston Maps}

We start off with some preliminaries about  hyperbolic metric
spaces  in the sense
of Gromov \cite{gromov-hypgps}. For details, see \cite{CDP}, \cite{GhH}. Let $(X,d_X)$
be a (Gromov) hyperbolic metric space. The 
{\bf Gromov boundary} of 
 $X$, denoted by $\partial{X}$,
is the collection of equivalence classes of geodesic rays $r:[0,\infty)
\rightarrow X$ with $r(0)=x_0$ for some fixed ${x_0}\in{X}$,
where rays $r_1$
and $r_2$ are equivalent if $sup\{ d_X(r_1(t),r_2(t))\}<\infty$.
Let $\widehat{X} = X\cup\partial{X}$ denote the natural
 compactification of $X$ topologized the usual way (cf.\cite{GhH} pg. 124).

We denote the $k$-neighborhood of a subset $Z$ of $(X,d_X)$ by $N_k(Z,d_X)$ or simply $N_k(Z)$ when $d_X$ is
understood.

\begin{defn} {\rm  A subset $Z$ of $(X,d_X)$ is said to be 
{\bf $k$-quasiconvex}
 if any geodesic joining points of  $ Z$ lies in a $k$-neighborhood $N_k(Z,d_X)$ of $Z$.
A subset $Z$ is  quasiconvex if it is $k$-quasiconvex for some
$k$. } \end{defn}

For  simply connected real hyperbolic
manifolds this is equivalent to saying that the convex hull $CH(Z)$ of the set
$Z$ lies in a  bounded neighborhood of $Z$. We shall have occasion to
use this alternate characterization.

\begin{defn} {\rm
A map $f$ from one metric space $(Y,{d_Y})$ into another metric space 
$(Z,{d_Z})$ is said to be
 a {\bf $(K,\epsilon)$-quasi-isometric embedding} if
 
\begin{center}
${\frac{1}{K}}({d_Y}({y_1},{y_2}))-\epsilon\leq{d_Z}(f({y_1}),f({y_2}))\leq{K}{d_Y}({y_1},{y_2})+\epsilon$
\end{center}

If  $f$ is a quasi-isometric embedding, 
 and every point of $Z$ lies at a uniformly bounded distance
from some $f(y)$ then $f$ is said to be a {\bf quasi-isometry}.

A $(K,{\epsilon})$-quasi-isometric embedding that is a quasi-isometry
will be called a $(K,{\epsilon})$-quasi-isometry.

A {\bf $(K,\epsilon)$-quasigeodesic}
 is a $(K,\epsilon)$-quasi-isometric embedding
of
a closed interval in $\mathbb{R}$. A $(K,K)$-quasigeodesic will also be called
a $K$-quasigeodesic. A $(K,K)$-quasigeodesic will simply be called a $K$-quasigeodesic} \end{defn}

We shall say that two paths $\alpha, \beta$ in $X$ {\it `$C$-track' each other} in $A \subset X$ if $\alpha \cap A$ and
 $\beta \cap A$ lie in a $C$ neighborhood of each other. The following Lemma says that quasigeodesics starting
and ending close by track each other.

\begin{lemma} \cite{GhH} Let $(X,d)$ be a $\delta$-hyperbolic metric space. Then for any $K,\epsilon, D $ there exists $C
= C(\delta, K,\epsilon, D)$ such that if $\alpha  , \beta$ are two  $(K,\epsilon )$-quasi-geodesics whose starting points 
(as also ending points) are at most $D$ apart, then $\alpha \subset N_C(\beta , d)$. \label{fellow} \end{lemma}

The conclusion of Lemma \ref{fellow} above is also
summarized by saying that  $\alpha  , \beta$  {\it $C-$fellow travel each other} and this property of quasi-geodesics
is called the $C-$fellow traveler property.

\smallskip

Let $(X,{d_X})$ be a (Gromov) hyperbolic metric space and $Y$ be a subspace that is (Gromov)
hyperbolic with the inherited path metric $d_Y$.
By 
adjoining the Gromov boundaries $\partial{X}$ and $\partial{Y}$
 to $X$ and $Y$, one obtains their compactifications
$\widehat{X}$ and $\widehat{Y}$ respectively.

Let $ i :Y \rightarrow X$ denote inclusion.

\begin{defn} {\rm  Let $X$ and $Y$ be (Gromov) hyperbolic metric spaces and
$i : Y \rightarrow X$ be an embedding. 
 A {\bf Cannon-Thurston map} $\hat{i}$  from $\widehat{Y}$ to
 $\widehat{X}$ is a continuous extension of $i$. } \end{defn}

The following  lemma (Lemma 2.1 of \cite{mitra-ct})
 says that a Cannon-Thurston map exists
if for all $M > 0$ and $y \in Y$, there exists $N > 0$ such that if $\lambda$
lies outside an $N$ ball around $y$ in $Y$ then
any geodesic in $X$ joining the end-points of $\lambda$ lies
outside the $M$ ball around $i(y)$ in $X$.
For convenience of use later on, we state this somewhat
differently and include the  proof from \cite{mahan-bddgeo}  for completeness.

\begin{lemma}
A Cannon-Thurston map from $\widehat{Y}$ to $\widehat{X}$
 exists if  the following condition is satisfied:\\
Given ${y_0}\in{Y}$, there exists a non-negative function  $M(N)$, such that 
 $M(N)\rightarrow\infty$ as $N\rightarrow\infty$ and for all geodesic segments
 $\lambda$  lying outside an $N$-ball
around ${y_0}\in{Y}$  any geodesic segment in $X$ joining
the end-points of $i(\lambda)$ lies outside the $M(N)$-ball around 
$i({y_0})\in{X}$.
\label{contlemma}
\end{lemma}

\begin{proof}
 Suppose $i:Y\rightarrow{X}$
does not extend continuously . Since $i$ is proper, there exist 
sequences $x_m$, $y_m$ $\in{Y}$ and $p\in\partial{Y}$,
such that $x_m\rightarrow p$
and $y_m\rightarrow p$ in $\widehat{Y}$, but $i(x_m)\rightarrow u$
and $i(y_m)\rightarrow v$ in $\widehat{X}$, where 
$u,v\in\partial{X}$ and $u\neq v$.

Since $x_m\rightarrow p$ and $y_m\rightarrow p$, any geodesic in $Y$
joining $x_m$ and $y_m$ lies outside an $N_m$-ball ${y_0}\in{Y}$,
 where $N_m\rightarrow\infty$ as $m\rightarrow\infty$. Any
bi-infinite geodesic in $X$ joining  $u,v\in\partial{X}$
has to pass through some $M$-ball around $i({y_0})$ in $X$ as
$u\neq v$. There exist constants $c$ and $L$ such that for all $m > L$
any geodesic joining $i(x_m)$ and $i(y_m)$ in $X$ 
passes through an $(M+c)$-neighborhood
 of $i({y_0})$. 
Since $(M+c)$ is a constant not depending on the index $m$
this proves the lemma. \end{proof}

The above result can be interpreted as saying that a Cannon-Thurston map 
exists if the space of geodesic segments in $Y$ embeds properly in the
space of geodesic segments in $X$.

\section{Relative Hyperbolicity} \label{relh}

In this section, we shall recall first certain notions of relative
hyperbolicity due to Farb \cite{farb-relhyp}, Klarreich
\cite{klarreich}, Bowditch \cite{bowditch-relhyp} and the author \cite{mahan-ibdd}. 

\subsection{Electric Geometry}

We collect together certain facts about the electric metric that Farb
proves in \cite{farb-relhyp}. 

\begin{defn} {\rm
Given a metric space $(X,d_X)$ and a collection $\mathcal{H}$ of subsets, let  
$\EE(X,\HH ) = X \bigsqcup_{H \in \HH} (H \times [0,\frac{1}{2}])$ be the 
identification space obtained by  identifying $(h,0) \in H \times [0,\frac{1}{2}]$ 
  with $h \in X$. Each $\{ h \} \times [0,\frac{1}{2}]$ is declared to be
isometric to the  interval $[0, \frac{1}{2} ]$
 and  $H \times \{ \frac{1}{2} \}$ is equipped  with the zero metric. 

A path $\sigma : I \rightarrow \EE(X,\HH )$ is said to be {\bf distinguished} if $\sigma (I) \cap \{ h \} \times (0,\frac{1}{2})$
is either empty or all of $\{ h \} \times (0,\frac{1}{2})$. 
The distance between two points in $\EE(X,\HH )$ is defined to be the infimum of the lengths of  distinguished
paths between them.

The resulting pseudo-metric space 
$\EE(X,\HH )$ is the {\bf electric space} associated to $X$ and the collection $\HH$. \\
We shall say that $\EE(X,\HH )$  is constructed from $X$ by {\bf electrocuting} the collection $\HH$ and the induced pseudo-metric $d_e$ will 
be called the {\bf
electric metric}.\\
If  $\EE(X,\HH )$ is (Gromov) hyperbolic, we say that $X$ is {\bf weakly hyperbolic} relative to $\HH$. }
\end{defn}

The notion of electrocution above is slightly different from the coning construction introduced by Farb in  \cite{farb-relhyp}, inasmuch as
Farb \cite{farb-relhyp}) collapses $H \times \{ \frac{1}{2} \}$ to a point. Thus ours is a geometric
(as opposed to topological) version of Farb's construction. {\it All paths in $\EE(X,\HH )$ will henceforth
be assumed to be distinguished.}

If $X$ is a geodesic metric space and each  $H \in \mathcal{H}$ is closed, then $(\EE(X,\HH ),d_e)$  is a geodesic (pseudo) metric space.
 Geodesics and quasigeodesics in $(\EE(X,\HH ),d_e)$ will be referred to as electric geodesics and electric quasigeodesics 
respectively.

Note that since  $\EE(X,\HH ) = X \bigsqcup_{H \in \HH} (H \times [0,\frac{1}{2}])$, $X$ can be naturally identified with a subspace of $\EE(X,\HH )$.
Paths in $(X,d_X)$ (in particular geodesics and quasigeodesics) can therefore be regarded as paths in $\EE(X,\HH )$.

A collection $\mathcal{H}$ of subsets of $(X,d_X)$ is said to be $D$-separated if $d_X(H_1, H_2) \geq D$
for all $H_1, H_2 \in \HH; H_1 \neq H_2$.
$D$-separatedness is only a technical restriction as the collection $\{ H \times \{ \frac{1}{2} \}: H \in \HH \}$ is $1$-separated in $\EE(X,\HH )$.

\begin{defn} \label{bt} {\rm  
Given a collection $\mathcal{H}$
of $C$-quasiconvex, $D$-separated sets in a (Gromov) hyperbolic metric space $(X,d_X)$ and a number $\epsilon$ we
shall say that a geodesic (resp. quasigeodesic) $\gamma$ is a geodesic
(resp. quasigeodesic) {\bf without backtracking} with respect to
$\epsilon -$ neighborhoods if $\gamma$ does not return to $N_\epsilon
(H, d_X)$ after leaving it, for any $H \in \mathcal{H}$. 
A geodesic (resp. quasigeodesic) $\gamma$ is a geodesic
(resp. quasigeodesic)  without backtracking if it is a geodesic
(resp. quasigeodesic) without backtracking with respect to
$\epsilon$ neighborhoods for $\epsilon = 0$.} \end{defn}

\noindent  {\bf Notation:} For any pseudo metric space $(Z, \rho)$ and $A\subset Z$, 
we shall use the notation $N_R(A, \rho) = \{ x \in Z: \rho(x, A) \leq R \}$ as for metric spaces.

\begin{lemma} (Lemma 4.5 and Proposition 4.6 of \cite{farb-relhyp}; Theorem 5.3 of \cite{klarreich}; \cite{bowditch-relhyp})\\
Given $\delta , C, D$ there exists $\Delta$ such that
if $(X,d_X)$ is a $\delta$-hyperbolic metric space with a collection
$\mathcal{H}$ of $C$-quasiconvex $D$-separated sets.
then,\\
1) {\rm Electric quasi-geodesics electrically track (Gromov) hyperbolic
  geodesics:} For all $P > 0$, there exists $K > 0$ such that if $\beta$ is any electric $P$-quasigeodesic from $x$ to
  $y$, and  $\gamma$ is a geodesic in $(X,d_X)$ from $x$ to $y$, 
then $\beta \subset N_K ( \gamma, d_e )$. \\
2) $\gamma \subset N_K ((N_0 ( \beta, d_e)), d_X)$. \\
3) {\rm Relative Hyperbolicity:} 
  $X$  is weakly hyperbolic relative to $\HH$. $\EE(X, \HH )$ is $\Delta$-hyperbolic.\\
\label{farb1A}
\end{lemma}

Let $(X,d_X)$ be a $\delta$-hyperbolic metric
space, and $\mathcal{H}$ a family of $C$-quasiconvex, $D$-separated,
 collection of subsets. Then $X$ is weakly hyperbolic relative to $\mathcal{H}$ \cite{bowditch-relhyp}.
Let $\alpha = [a,b]$ be a geodesic in $(X,d_X)$ and $\beta $ 
an electric 
quasigeodesic without backtracking (in $\EE(X,\HH )$)
 joining $a, b$. Order from the left the collection of maximal subsegments of $\beta$
contained entirely in some $H \times \frac{1}{2}$ for  some $H \in \mathcal{H}$. 
Since $\beta$ is a distinguished path (by our blanket assumption about paths in $\EE(X,\HH )$),
 any such maximal subsegment can be extended by adjoining vertical subsegments at its end-points
to obtain a path of the form $\{ p \} \times [0,\frac{1}{2}]\cup [p \times \frac{1}{2},
 q\times \frac{1}{2}] \cup \{ q \} \times [0,\frac{1}{2}]$. We shall refer to
 these subpaths of $\beta$ as  extended  maximal subsegments.
Replace, as per the above ordering,   extended  maximal subsegment with end-points $p,q$ (say)
by a   geodesic $[p,q]$ in $(X,d_X)$. (Note here that as per the definition of $\EE(X,\HH )$, $(p,0) \in \EE(X,\HH )$
is identified with $p\in X$; similarly for $(q,0)$ and $q$.)
The resulting
{\bf connected}
path $\beta_q$ is called an {\em electro-ambient representative} of $\beta$ in
$X$. Also, if $\beta $ is
an electric  $P$-quasigeodesic (resp. $(K,\epsilon)$-quasigeodesic) without backtracking (in $\EE(X,\HH )$), then $\beta_q$ is called an {\bf electro-ambient $P$-quasigeodesic  (resp. electro-ambient $(K,\epsilon)$-quasigeodesic).}
If  $\beta $ is
an electric  geodesic (i.e. a $(1,0)$-quasigeodesic) without backtracking (in $\EE(X,\HH )$), then $\beta_q$ is simply called an {\bf electro-ambient
quasigeodesic.}

\begin{rmk} \label{ea-clarfn} We emphasize a point about the terminology we use here. An electro-ambient
quasigeodesic in our sense is the same as an  electro-ambient $(1,0)$-quasigeodesic, {\bf not} an
 electro-ambient $(K,\epsilon)$-quasigeodesic for some $K,\epsilon$.
\end{rmk}

Note that $\beta_q$ {\it need
  not be a (Gromov) hyperbolic quasigeodesic}. However, the proof of Proposition
  4.3 of Klarreich \cite{klarreich} gives the following.

\begin{lemma} (See Proposition 4.3 of \cite{klarreich}, also see Lemma
  3.10 of \cite{mahan-ibdd}) 
Given $\delta$, $C, P$ there exists $C_3$ such that the following
holds: \\
Let $(X,d_X)$ be a $\delta$-hyperbolic metric space and $\mathcal{H}$ a
family of $C$-quasiconvex,  collection of quasiconvex
subsets. Let $(\EE(X,\HH ),d_e)$ denote the electric space obtained by
electrocuting elements of $\mathcal{H}$.  Then, if $\alpha , \beta_q$
denote respectively a (Gromov) hyperbolic geodesic and an electro-ambient
$P$-quasigeodesic with the same end-points in $X$, then $\alpha$ lies in a
(Gromov hyperbolic $d_X-$) 
$C_3$ neighborhood of $\beta_q$.
\label{ea-strong}
\end{lemma}

For the convenience of the reader, we illustrate the content of Lemma \ref{ea-strong} by the figure below:

\begin{center}

\includegraphics[height=4cm]{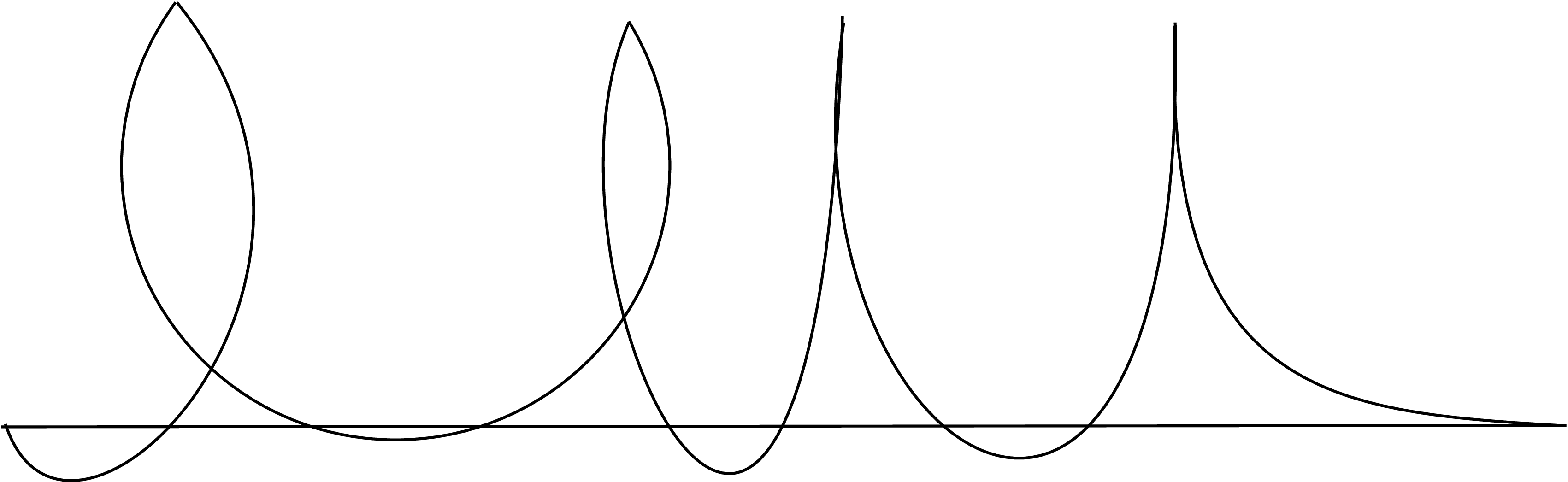}

\smallskip

\underline{Figure:  {\it Hyperbolic geodesic lies in a neighborhood of an electro-ambient quasigeodesic} }

\end{center}

In the above figure, the straight line below indicates a hyperbolic geodesic, and the broken line built up of curves depicts
 an electro-ambient quasigeodesic.

\smallskip

\noindent {\bf  Proof of Lemma \ref{ea-strong}:} The proof follows Proposition
  4.3 of Klarreich \cite{klarreich} closely. Let $\alpha = [a,b] (\subset X)$  be a geodesic and let $\beta = \overline{ab} (\subset \EE(X,\HH ))$ 
be an electric $P$-quasigeodesic with the same end-points. Further, suppose that, for each $H \in \HH$,  $\beta \cap (H\times  \{\frac{1}{2}\})$ is \\
a) a maximal subsegment of $\beta$
contained in $H\times  \{\frac{1}{2}\}$, \\
b) $\beta \cap (H\times  \{\frac{1}{2}\})$ is a geodesic in $H\times  \{\frac{1}{2}\}$ with respect to the intrinsic metric on $H (= H\times \{\frac{1}{2}\})$.\\

For the purposes of this proof, we shall need to deal with two metrics (more precisely a metric and a pseudometric) on the {\it topological space}
$\EE(X,\HH )$: \\
a)  The first is the electric (pseudo) metric $d_e$ described above. \\
b) The other is the (genuine) metric on $X \bigcup_{H \in \HH} H \times [0,  \frac{1}{2}]$ obtained as a quotient space of $X$ along with
copies of $H \times [0,  \frac{1}{2}]$. We call this metric $d_q$.\\

Thus $d_e$ is obtained from $d_q$ by redefining distance between points on $H\times  \{\frac{1}{2}\}$ to be zero.

Recall that (by construction) the electro-ambient quasigeodesic  $\beta_q$ is obtained from $\beta$ by "projecting" maximal subsegments of $\beta$ to $X$.
It therefore suffices to show that $\alpha$ lies in a (uniformly) bounded neighborhood of $\beta$ in $(\EE(X,\HH ), d_q)$. 

Let $a=a_0, a_1, \cdots, a_n, a_{n+1}=b$ be a sequence of points on $\beta$ such that for all $i$, $\overline{a_{2i} a_{2i+1}}(\subset \overline{a b})$
 are maximal subsegments in $H_i \times \{\frac{1}{2}\}$
for some $H_i \in \HH$. Also, assume that $n$ is maximal, i.e. for all $i$, $\overline{a_{2i-1} a_{2i}}$ is a union of three segments: \\
a) a vertical segment of the form $a_{2i-1} \times [0,  \frac{1}{2}]$ traced from $a_{2i-1} \times \{ \frac{1}{2} \}$ to $a_{2i-1} \times \{ 0 \}$, \\
b) a geodesic in $(X, d_X)$ from $a_{2i-1}$ (identified with $a_{2i-1} \times \{ 0 \}$) to $a_{2i}$   (identified with $a_{2i} \times \{ 0 \}$),\\
c) a vertical segment of the form $a_{2i}\times [0,  \frac{1}{2}]$ traced from $a_{2i} \times \{ 0 \}$ to $a_{2i} \times \{ \frac{1}{2} \}$.\\

Note first that the collection $\{ H\times \{\frac{1}{2}\}\}$, $H \in \HH$ is automatically $1-$ separated. Hence $d_e(a_{2i-1}, a_{2i}) \geq 1$.

With this setup, the proof is a small reworking of Proposition
  4.3 of  \cite{klarreich}. Choose an $R>0$. Let $z \in [a,b]$ be a point for which no point of  $\beta = \overline{ab}$ lies within $R$ of $z$.
Let $(p,q)$ be a maximal subsegment of $[a,b]$ containing $z$ such that no point of  $\beta = \overline{ab}$ lies within $R$ of $(p,q)$.

Let $p_1 \in \overline{ab}$ and  $q_1 \in \overline{ab}$ be points in  $\overline{ab}$ closest to $p, q$ respectively (with respect to the metric $d_q$).
Let $\overline{p_1q_1}$ be the subpath of $\overline{ab}$ between $p_1,q_1$. Also, let $a_j, \cdots , a_{j+l}$ be the collection of vertices
in $\overline{ab}$ between $p_1,q_1$. Then the proof of Proposition
  4.3 of  \cite{klarreich} shows that there exists $R_0$ depending on $\delta$, $C, P$ such that for all $R \geq R_0$, $l (=(j+l)-j)$ is bounded in terms of 
$R, \delta$, $C, P$.  (This is essentially because $l$ grows like $d_X(p,q)e^R$, cf. \cite{farb-relhyp}.)
Let $l(R)$ be this bound for $l$.

Choosing $R=R_0$, we find that $(p,q)$ is contained in a $(2R_0 + (l(R_0) + 4) \delta)-$ neighborhood of $\beta = \overline{ab}$. This completes
the proof. $\Box$
 
\begin{definition}  \cite{farb-relhyp} {\rm
Two paths   $\beta , \gamma$ in $(X,d_X)$ with  the same endpoints  are said to have \emph{similar intersection patterns} 
with $\HH$ if  there exists $ D>0$, depending only on $(X,\HH)$, such that:
\begin{itemize}
\item {\it Similar Intersection Patterns 1:}  If
  precisely one of $\{ \beta , \gamma \}$ meets some
    $H \in \mathcal{H}$, then the $d_X$-distance  from the first entry point
  to the 
 last exit point is at most $D$. 
\item {\it Similar Intersection Patterns 2:}  If
 both $\{ \beta , \gamma \}$ meet some    $H \in \mathcal{H}$,
 then the distance  from the first entry point of
 $\beta$ to that of $\gamma$ is at most $D$, and  similarly for the last exit points. 
\end{itemize}
}
\end{definition}

\begin{definition}  \cite{farb-relhyp} {\rm
Suppose that $X$ is 
 weakly hyperbolic relative to $\mathcal{H}$.
Suppose that any two electric quasigeodesics without backtracking and with the same endpoints  have similar intersection patterns with $\mathcal{H}$.
Then $(X,\HH)$ is said to satisfy  {\bf bounded penetration}  and 
 $X$ is said to be
{\bf strongly hyperbolic} relative to  $\mathcal{H}$.}
\end{definition}

The next condition ensures that $(X,\HH)$ is  {\bf strongly hyperbolic} relative to  $\mathcal{H}$.

\begin{definition} {\rm A collection $\mathcal{H}$ of  sets in a $\delta$-hyperbolic metric space $X$ is said to be 
{\bf uniformly D-separated} if $d(H_i, H_j) \geq D$ for all $H_i, H_j \in \HH; H_i\neq H_j$.\\
A collection $\mathcal{H}$ of uniformly
$C$-quasiconvex sets in a $\delta$-hyperbolic metric space $X$
is said to be {\bf mutually D-cobounded} if 
 for all $H_i, H_j \in \mathcal{H}$, $\pi_i
(H_j)$ has diameter less than $D$, where $\pi_i$ denotes a nearest
point projection of $X$ onto $H_i$. A collection is {\bf mutually
  cobounded} if it is mutually D-cobounded for some $D$. 
}
  \end{definition}

Coboundedness ensures strong relative hyperbolicity.

\begin{lemma} (\cite{farb-relhyp} Proposition 4.6, \cite{bowditch-relhyp})
Let $(X,d_X)$ be a (Gromov) hyperbolic metric space and
$\mathcal{H}$ a collection of  $\epsilon$-neighborhoods of  uniformly quasiconvex mutually cobounded
uniformly separated subsets.  
Then $X$ is strongly hyperbolic relative to the collection $\mathcal{H}$. Furthermore quasigeodesics without backtracking
in $(X,d_X)$ and $(\EE (X, \HH ), d_{e})$ have similar intersection patterns with elements of $\HH$.
\label{farb2A}
\end{lemma}

Applications of Lemma \ref{farb2A} follow.

\begin{lemma}
Let $M^h$ be a hyperbolic manifold. Let $\mathcal{T}$ and $\mathcal{H}$ denote respectively a collection of 
Margulis tubes  and horoballs that are disjoint from one another. Then the elements of $\mathcal{T} \cup \mathcal{H}$
are mutually co-bounded. Hence $\widetilde{M^h}$ is strongly hyperbolic relative to the collection $\mathcal{T}  \cup  \mathcal{H}$.
\label{coboundedHor&T}
\end{lemma}

\begin{lemma}
Let $S^h$ be a hyperbolic surface, with a finite collection of disjoint simple closed geodesics
$\sigma_i $ and cusps $H_j $. Let $\SSS$ denote the collection of lifts $\widetilde{\sigma_i}$ to $\Hyp^2$  and let
$\HH$ denote the collection of lifts $\widetilde{H_j}$. Then the elements of  $\SSS \cup \HH$ 
are mutually co-bounded. Hence $\widetilde{S^h}$ is strongly hyperbolic relative to the collection $\SSS  \cup \HH$.
\label{coboundedgeods}
\end{lemma}

A closely related  theorem was proved by  McMullen
(Theorem 8.1 of \cite{ctm-locconn}). 
Let $\mathcal{H}$ be a locally finite collection of horoballs in a convex
subset $X$ of ${\mathbb{H}}^n$ 
(where the intersection of a horoball, which meets $\partial X$ in a point, 
 with $X$ is
called a horoball in $X$).

\begin{defn} {\rm The $\epsilon$-neighborhood of a bi-infinite
geodesic in ${\mathbb{H}}^n$ will be called a {\bf thickened geodesic}. } \end{defn}

\begin{theorem} \cite{ctm-locconn}
For $K, D \geq 1, \epsilon\geq 0$ there exists $R \geq 0$ such that the following holds:\\
Let $X$ be a convex subset  of ${\mathbb{H}}^n$ and let
$\mathcal{H}$  denote a uniformly $D$-separated
collection of horoballs and thickened geodesics.
Let $Y = X \setminus \bigcup_{H\in \mathcal{H}}H$ and $\gamma: I \rightarrow Y$ be a
$(K, \epsilon )$-quasigeodesic in $Y$. 
Let $\eta$ be the  geodesic in $X$ with the same endpoints as
$\gamma$. Let $\mathcal{H}({\eta})$  
be the union of all the horoballs and thickened geodesics
 in $\mathcal{H}$ meeting $\eta$. Then
$\eta\cup\mathcal{H}{({\eta})}$ is $R-$quasiconvex and $\gamma
(I) \subset  
B_R (\eta \cup \mathcal{H} ({\eta}))$. (The hyperbolic metric on ${\mathbb{H}}^n$ is understood.)
\label{ctm}
\end{theorem}

\subsection{Electric Geometry for Surfaces}
We now specialize to surfaces. 
Let $S$ be a hyperbolic surface with diameter bounded above by $K$. It follows that  injectivity radius is 
bounded below by some  $\epsilon = \epsilon (K)$. Let $\sigma$ be a finite collection of disjoint simple closed geodesics
on $S$. Component(s) of
 $S \setminus \sigma$ will be called the {\bf amalgamation component(s)}
of $(S, \sigma )$. We shall denote an amalgamation components by $S_A$. Let $(S_{Gel}, d_{Gel}) = \EE (S,S_A)$ be
 obtained from $S$ by electrocuting $S_A$'s and let the universal cover of $(S_{Gel}, d_{Gel})$ with
 the lifted pseudometric be denoted $(\til{S_{Gel}}, d_{Gel})$.
 A slightly different  path pseudometric may be constructed
on $\til S$ by declaring that\\
\noindent  
$1)$ the length of any path that lies in the interior of
 an amalgamation component is zero \\
$2)$  the length of any path  that crosses $\sigma$ once has length one \\
$3)$ the length of any other path is the sum of lengths of
pieces of the above two kinds. \\

This pseudometric differs from $(\til{S_{Gel}}, d_{Gel})$ by at most one (due to the initial and final
segments of length half). {\it We shall ignore this difference (cf Lemma \ref{hyp=elproj}).}

The fundamental group $\pi_1(S)$ may be regarded as a graph of groups with
vertex group(s) the subgroup(s)  $\pi_1(S_A)$ corresponding to amalgamation
component(s) and cyclic edge groups $\mathbb{Z}$ corresponding to $\sigma$. Then $(\widetilde{S_{Gel}}, d_{Gel})$
  is quasi-isometric to  the Bass-Serre tree of  the
splitting.

 Continuous paths in $S_{Gel}$ and $\widetilde{S_{Gel}}$ will be called
{\it electric paths}. Continuous geodesics and quasigeodesics in 
the electric metric will be called electric geodesics and electric
quasigeodesics respectively. We specialize Definition \ref{bt} to the present context, where it is slightly more
restrictive.

\begin{defn}{\rm 
An electric path $\gamma \subset \widetilde{S_{Gel}}$ is said to be an electric
$K$-quasigeodesic in $(\til{S_{Gel}}, d_{Gel})$ 
{ \bf without backtracking } if
 $\gamma$ is a $K$-quasigeodesic in  $(\til{S_{Gel}}, d_{Gel})$  and
 $\gamma$ does not return to 
any lift  $\widetilde{S_A} (\subset
\widetilde{S_{Gel}})$ of an amalgamation component
$S_A \subset S$ after leaving it.}\end{defn}

We now specialize the notion of an electro-ambient quasigeodesic to the context of surfaces.

\begin{defn}{\rm 
An electric geodesic $\lambda_e$ without backtracking in 
$(\widetilde{S_{Gel}},d_{Gel})$ is called an {\bf electro-ambient quasigeodesic} if\\
a) each segment of $\lambda_e$ lying inside a single lift $\widetilde{S_A}$ meets the boundary $\partial \widetilde{S_A}$
at most twice and is perpendicular to $\partial \widetilde{S_A}$ whenever they meet. We shall refer to these segments of $\lambda_e$ as
{\bf amalgamation segments}.  \\
b) If $a , b$ be the points of intersection of two distinct amalgamation segments of $\lambda_e$ with 
 a lift $\widetilde{\sigma}$ of
$\sigma$, then $\lambda_e \cap \widetilde{\sigma}$ is equal to $[a,b]$, the geodesic segment in  $\widetilde{\sigma}$ joining $a,
b$.  Such pieces $[a,b]$ shall be referred to
as  {\bf interpolating segments}.\\
The underlying path of  an electro-ambient quasigeodesic of the electro-ambient quasigeodesic
 in the {\em hyperbolic metric} on $\widetilde{S}$ shall be called the {\bf
  electro-ambient representative} $\lambda_q$ of $\lambda_e$. }\end{defn}

 See Figure below, where the bold line indicates the  electro-ambient quasigeodesic and the thin lines the geodesics $\widetilde{\sigma}$.

\begin{center}

\includegraphics[height=6cm]{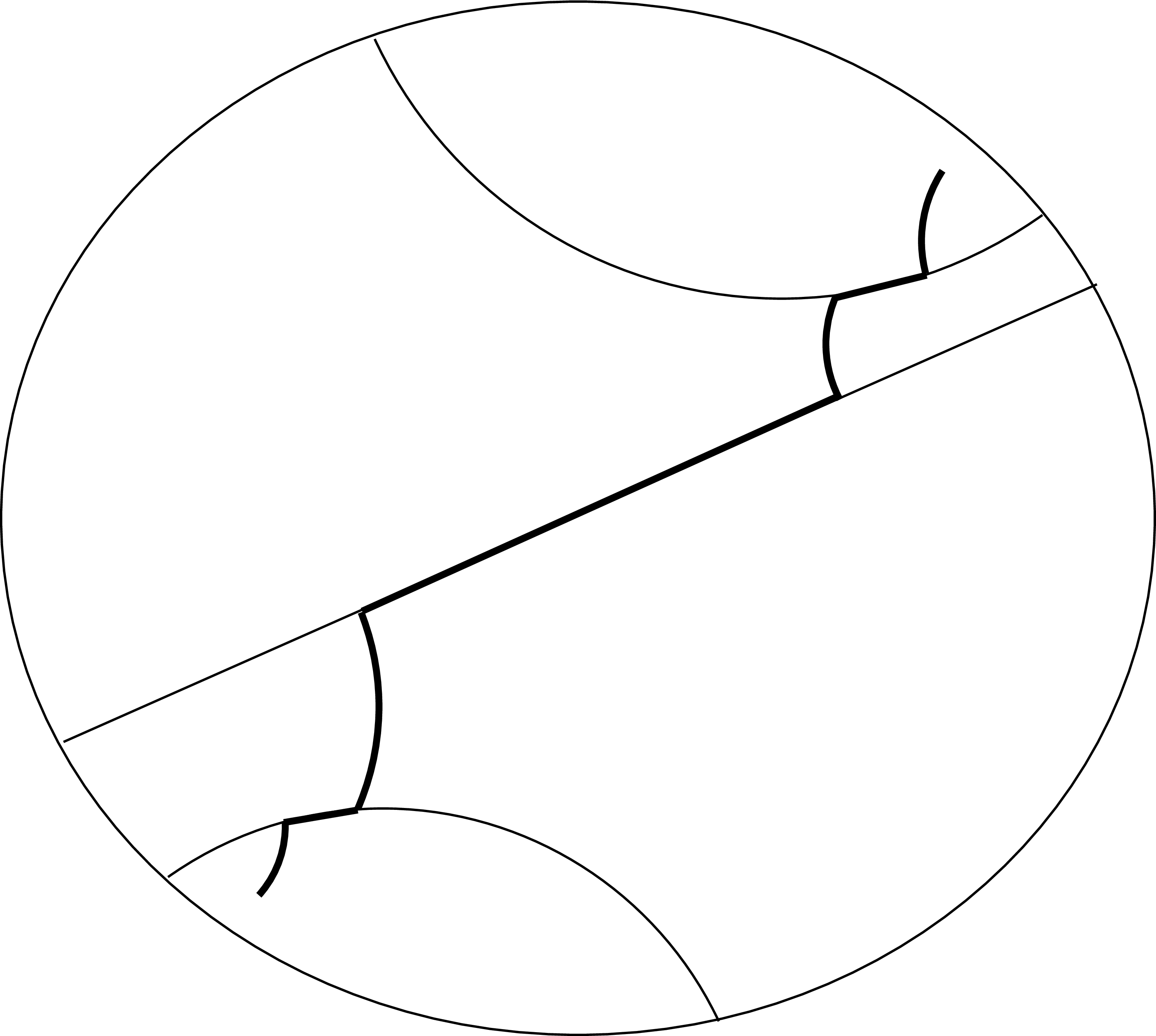}

\underline{Figure:{\it Electro-ambient quasigeodesic} }

\end{center}

The next Lemma justifies the terminology.

\begin{lemma}(See Lemma 3.7 of \cite{mahan-ibdd})
There exists $(K, \epsilon)$ such that each electro-ambient
representative $\lambda_q$ of an electric geodesic in
$(\widetilde{S_{Gel}}, d_{Gel})$
is a $(K,
\epsilon)$ hyperbolic quasigeodesic.
\label{Gea}
\end{lemma}

  \begin{proof} Let $(S_{el}, d_{el})$ denote the surface $S$ with the (collection of) geodesics $\sigma$ electrocuted.
Note that the electro-ambient quasigeodesics in $(\widetilde{S_{Gel}}, d_{Gel})$ coincide with those in the universal cover
$(\widetilde{S_{el}}, d_{el})$. Hence it suffices to show that 
electro-ambient quasigeodesics in 
$(\widetilde{S_{el}}, d_{el})$ are uniform hyperbolic quasigeodesics. 

Let $\lambda_h$ denote the hyperbolic geodesic joining the
end-points of $\lambda_e$. By Lemmas \ref{farb2A} and \ref{coboundedgeods}, $\lambda_h$ and
$\lambda_e$, and hence $\lambda_h$ and $\lambda_q$
have similar intersection patterns with $N_\epsilon (\til{\sigma}) $ for some small $\epsilon > 0$
and any lift $\til{\sigma}$ of (an element of) $\sigma$.
 Also, $\lambda_h$ and $\lambda_q$
track each other off the collection $N_\epsilon (\til{\sigma})$. Further, each {\em
  interpolating segment} of $\lambda_q$ being a {\em hyperbolic} geodesic, it
follows (from the `$K$-fellow-traveler' property of hyperbolic geodesics
starting and ending near each other, Lemma \ref{fellow})
that each interpolating segment of $\lambda_q$
lies within a $(K + 2 \epsilon )$
neighborhood of $\lambda_h$ for some fixed $K>0$. Again, since each segment of $\lambda_q$
that does not meet an electrocuted geodesic that $\lambda_h$ meets is 
of uniformly bounded length (bounded by $C$ say), we have finally that
$\lambda_q$ lies within a $(K + C + 2 \epsilon)$ neighborhood of
$\lambda_h$. Finally, since $\lambda_q$ is an electro-ambient
representative, it does not backtrack. Hence the Lemma. \end{proof}

\subsection{Electric isometries}

\begin{defn} {\rm Let $S$ be any hyperbolic surface and $\sigma$ a collection of disjoint simple closed
geodesics on $S$.
A  diffeomorphism  $\phi$ of $S$ will be called a {\bf component
  preserving diffeomorphism} if it fixes $\sigma$
pointwise and
 preserves each
amalgamation component as a set, i.e. $\phi$ sends each amalgamation
component of $(S, \sigma )$ to itself. }
\end{defn}

\begin{lemma}
Let $\phi$ denote a component preserving diffeomorphism of $S_G$. 
 Then $\phi $ induces an isometry of $(S_{Gel},d_{Gel})$.
\label{phi-isom1}
\end{lemma}

\begin{proof} In the electrocuted surface
$(S_{Gel}, d_{Gel})$, any electric geodesic $\lambda_e$ has length equal to the number of
times it crosses $\sigma$. Any component
  preserving diffeomorphism $\phi$ preserves the intersection pattern of $\lambda_e$ with amalgamation components.
Hence $\phi$ is an isometry of
$(S_{Gel}, d_{Gel})$. \end{proof}

The proof of Lemma \ref{phi-isom1} goes through verbatim after lifting to the universal cover
$(\widetilde{S_{Gel}}, d_{Gel})$. We  let $\widetilde{\phi}$ denote the lift of $\phi$ to
$(\widetilde{S_{Gel}}, d_{Gel})$.
This gives \\

\begin{lemma}
Let $\widetilde{\phi}$ denote a lift of a component preserving
 diffeomorphism $\phi$ to 
$(\widetilde{S_{Gel}},d_{Gel})$. 
 Then $\widetilde{\phi} $ induces an isometry of
$ ( \widetilde{S_{Gel}},d_{Gel})$.  
\label{phi-isom2}
\end{lemma}

\subsection{Nearest-point Projections}

The next Lemma  says nearest point projections in a $\delta$-hyperbolic
metric space do not increase distances much. This is a standard fact (See Lemma 3.1 of \cite{mitra-trees} for instance).

\begin{lemma}
Let $(Y,d_Y)$ be a $\delta$-hyperbolic metric space
 and  let $\mu\subset{Y}$ be a $C$-quasiconvex subset, 
e.g. a geodesic segment.
Let ${\pi}:Y\rightarrow\mu$ map $y\in{Y}$ to a point on
$\mu$ nearest to $y$. Then $d_Y{(\pi{(x)},\pi{(y)})}\leq{C_3}d_Y{(x,y)}$ for
all $x, y\in{Y}$ where $C_3$ depends only on $\delta, C$.
\label{easyprojnlemma}
\end{lemma}

The next lemma (from \cite{mitra-trees})  says that quasi-isometries and nearest-point projections on
(Gromov) hyperbolic metric spaces `almost commute'.

\begin{lemma}
(Lemma 3.5 of \cite{mitra-trees}) Suppose $(Y_1,d_1)$ and $(Y_2,d_2)$
are $\delta$-hyperbolic.
Let $\mu_1$ be some geodesic segment in $Y_1$ joining $a, b$ and let $p$
be any point of $Y_1$. Also let $q$ be a point on $\mu_1$ such that
${d_1}(p,q)\leq{d_2}(p,x)$ for all $x\in\mu_1$. 
Let $\phi$ be a $(K,{\epsilon})$ - quasiisometric embedding
 from $Y_1$ to $Y_2$.
Let $\mu_2$ be a geodesic segment 
in $Y_2$ joining ${\phi}(a)$ to ${\phi}(b)$ . Let
$r$ be a point on $\mu_2$ such that ${d_2}({\phi}(p),r)\leq{d_2}({\phi(p)},x)$ for $x\in\mu_2$.
Then ${d_2}(r,{\phi}(q))\leq{C_4}$ for some constant $C_4$ depending   only on
$K, \epsilon $ and $\delta$. 
\label{almost-commute}
\end{lemma}

We shall need the above Lemma for quasi-isometries
from $\widetilde{S_a}$ to  $\widetilde{S_b}$ for two different biLipschitz
 metrics on the same surface. We shall also need it for
 electrocuted surfaces.

Another property that we shall require for nearest point
projections is that nearest point projections in the hyperbolic metric on $\til S$
and that in the electric metric  $(\widetilde{S_{Gel}}, d_{Gel})$ almost agree.
 To make this precise
we make the following definition.  The hyperbolic metric on $S$ as well as $\til S$ will be denoted by $d$.

\begin{defn} {\rm Let $y \in (\til{S},d)$ and let $\mu_q$ be an electro-ambient
representative of an electric geodesic $\mu_{Gel}$ in
$(\til{S_{Gel}},d_{Gel})$. Then $\pi_e(y) = z \in \mu_q$ if the ordered pair $\{
d_{Gel}(y,\pi_e(y)), d(y, \pi_e(y) ) \}$ is minimized at $z$ in the lexicographical order on 
$({\mathbb{R}}_+ \cup \{ 0 \}) \times ({\mathbb{R}}_+ \cup \{ 0 \})$. } \end{defn}

The proof of the following Lemma shows
 that this gives us a definition of $\pi_e$ which is ambiguous up to
a finite amount of discrepancy not only in the electric metric but
also
in the hyperbolic metric.

\begin{lemma} Fix a hyperbolic surface $S$.
For all $\epsilon > 0$, there exists $C > 0$ such that if $\sigma$ is a finite collection of disjoint simple
closed geodesics such that $d(\sigma_i, \sigma_j) \geq \epsilon$ for all 
$\sigma_i \neq \sigma_j \in \sigma$, then the following holds.\\
Let $\mu$ be a hyperbolic geodesic in $(\til{S},d)$ joining $u, v \in \til{S}$. Let $({S_{Gel}},d_{Gel})$ be the electric 
space  obtained from $S$ by electrocuting the amalgamation components of $(S, \sigma )$. Let
  $\mu_{Gel}$ be an electric geodesic in $(\til{S_{Gel}},d_{Gel})$
  joining $u, v$ and let $\mu_q$ be its electro-ambient
  representative. Let $\pi_h$ denote the nearest point
  projection of $(\til{S},d)$ onto $\mu$. Then
$d(\pi_h(y) , \pi_e(y)) \leq C$.
\label{hyp=elproj}
\end{lemma}

\begin{proof} 
 Let $[u, v]$ and $[u,v]_q$ denote respectively the hyperbolic
geodesic and the electro-ambient quasigeodesic
joining $u, v \in \til{S}$.  Since $[u,v]_q$ is a hyperbolic quasigeodesic by Lemma \ref{Gea},
the nearest point projection of $y\in (\til{S},d)$ onto $[u, v]$ and $[u,v]_q$ almost
agree in the hyperbolic metric $d$. Thus, abusing notation slightly let $\pi_h$ denote nearest point projection of  
$(\til{S},d)$  onto $[u,v]_q$. Hence
it suffices to 
show that for any $y\in \til{S}$, its hyperbolic and electric  projections $\pi_h
(y), \pi_e (y)$ onto $[u,v]_q$  almost
agree. See figure below, where we denote $\pi_h
(y), \pi_e (y)$ by $p, q$ respectively.

\begin{center}

\includegraphics[height=6cm]{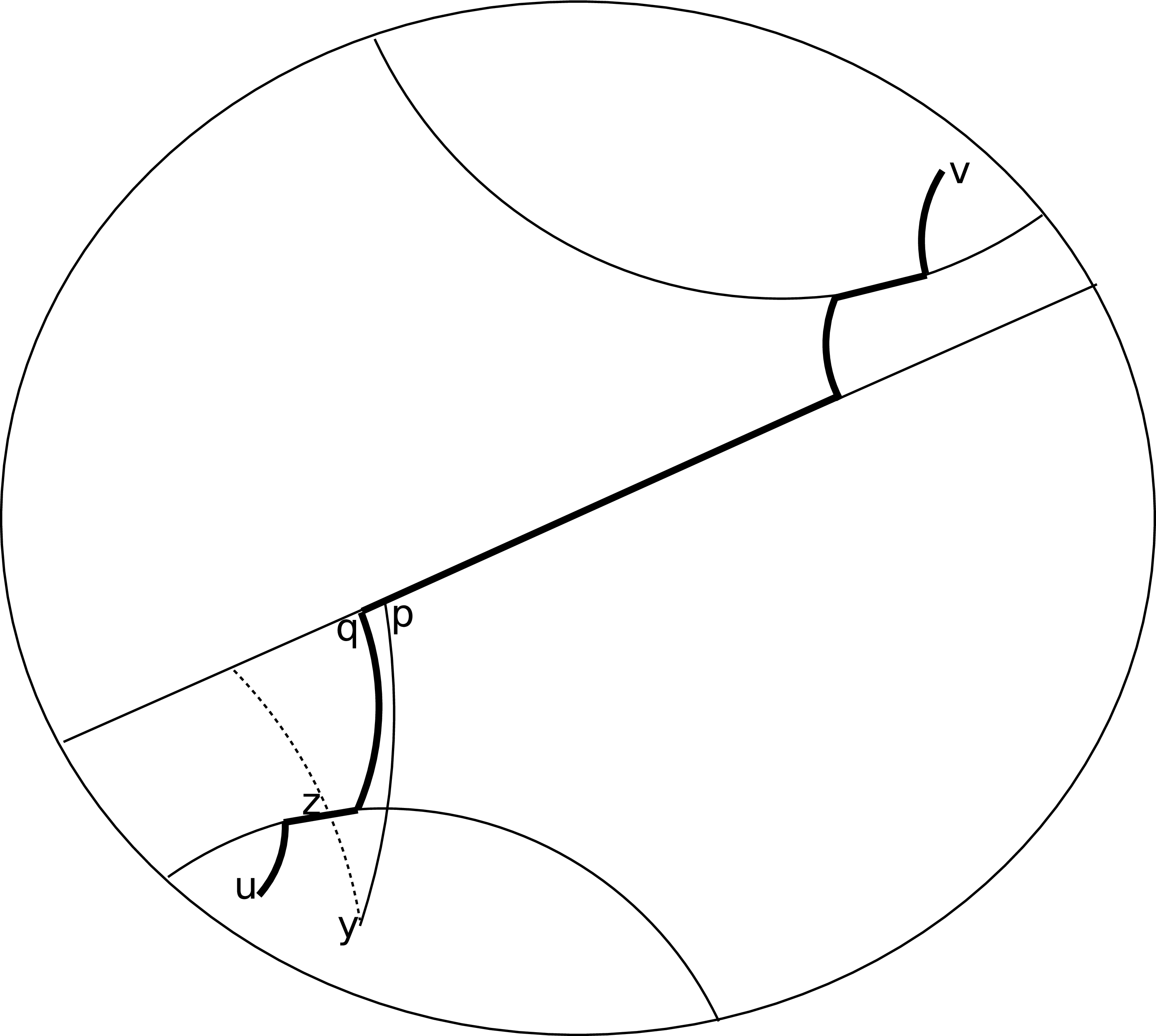}

\underline{Figure:{\it Electric and hyperbolic projections.} }

\end{center}

First note that  any  hyperbolic geodesic $\eta$ in $\widetilde{S}$
is also an 
electric geodesic in $(\til{S_{Gel}},d_{Gel})$. This follows from the fact that $(\til{S_{Gel}},d_{Gel})$
 maps to the Bass-Serre tree $T$ of the splitting of $S$ along $\sigma$, such that
 the pre-image of every vertex is a set of diameter zero in the
pseudometric $d_{Gel}$. If a path in $(\til{S_{Gel}},d_{Gel})$ projects to
a path in $T$ that is not a geodesic, then it must backtrack. Hence,
it must leave an amalgamating component and return to it. Such a path
can clearly not be a hyperbolic geodesic in $(\til{S_{Gel}},d_{Gel})$ since
each amalgamating component is convex. 

Next, it follows that  hyperbolic projections
automatically minimize electric distances. Else as in the preceding
paragraph, $[y,\pi_h (y)]$ would have to cut a lift $ \widetilde{\sigma_1}$ of
${\sigma}$
 that separates $[u,v]_q$. Further,  $[y,\pi_h (y)]$ cannot return to
 $\widetilde{\sigma_1 }$ after leaving it. 
Let $z$ be the first point at which
 $[y,\pi_h (y)]$ meets $\widetilde{\sigma_1}$ (the intersection point of the dotted line with $\widetilde{\sigma_1}$ 
in the figure above). Also let $w$ be the
point on  $[u,v]_q \cap \widetilde{\sigma_1}$ that is nearest to
$z$. Since  amalgamation segments of $[u,v]_q$ meeting
$\widetilde{\sigma_1}$ are perpendicular to the latter, it follows
that $d(w,z) < d(w,\pi_h (y) )$ and therefore  $d(y,z) < d(y,\pi_h (y)
)$ contradicting the definition of $\pi_h (y)$. Hence hyperbolic projections
automatically minimize electric distances.

Further, it follows by repeating the argument in the first paragraph
 that $[y,\pi_h (y)]$ and $[y, \pi_e (y)]$ pass through the same set
 of amalgamation components in the same order; in particular they
 cut across the same set of lifts of $\widetilde{\sigma}$. Let
 $\widetilde{\sigma_2}$ be the last such lift. Then
 $\widetilde{\sigma_2}$ forms the boundary of an amalgamation
 component $\widetilde{S_A}$ whose intersection with $[u,v]_q$ is of the
 form $[a,b] \cup [b,c] \cup [c,d]$, where $[a,b] \subset
 \widetilde{\sigma_3}$ and $[c,d] \subset \widetilde{\sigma_4}$ are
 subsegments of two lifts of $\sigma$ and $[b,c]$ is perpendicular to
 these two. Then the nearest-point 
projection of $\widetilde{\sigma_2}$ onto each of $[a,b], [b,c],
 [c,d]$ has uniformly bounded diameter. Hence the nearest point
 projection of $\widetilde{\sigma_2}$ onto the hyperbolic geodesic $[a,d]
 \subset \widetilde{S_A}$ has uniformly bounded diameter.
 The result follows. \end{proof}

\section{The Minsky Model} \label{min}

In this section we summarize the  notions and facts 
 from  \cite{minsky-elc1}, \cite{minsky-elc2} and  \cite{masur-minsky2} that we shall need.
Let $\CC (S)$ and $\PP (S)$ denote respectively the curve complex
and pants complex of a compact surface $S$, possibly with boundary, with the usual modifications for surfaces of small complexity
(see \cite{masur-minsky2} for details).

\smallskip

\noindent {\bf Split level Surfaces}\\
For our purposes, a pants decomposition of $S$ 
will be a disjoint collection of
3-holed spheres $P_1, \cdots , P_n$ embedded in $S$ such that $S \setminus \bigcup_i P_i$ is a disjoint collection of non-peripheral
annuli in $S$, no two of which are  homotopic. We shall conflate a pants decomposition of $S$ with the  collection of
(isotopy classes of) {\it non-peripheral}  boundary curves of  $P_1, \cdots , P_n$. Thus when we refer to a pair of pants in a pants decomposition 
 $P_1, \cdots , P_n$ of $S$ we are referring to one of the $P_i$'s, and when we refer to a curve in a pants decomposition 
 of $S$ we are referring to one of the  non-peripheral  boundary curves of one of the $P_i$'s.

Let $N$ be the convex core of a simply or doubly degenerate hyperbolic 3-manifold minus an open neighborhood of the cusp(s).

$N$ is homeomorphic to $S \times [0, \infty )$ or  $S \times \mathbb{R}$ according as 
$N$ is simply or doubly degenerate, where $S$ is a compact surface, possibly with boundary. 

Let $ \theta , \omega $ be positive real numbers. 
A neighborhood $N_\epsilon (\gamma )$ of a closed geodesic $\gamma (\subset N)$
is called a $(\theta, \omega)$-thin tube if the length of $\gamma$ is less than $\theta$ and the length of the shortest
geodesic on $\partial N_\epsilon (\gamma )$ is greater than $\omega$.

Let $\TT$ denote a collection of disjoint, uniformly separated $(\theta, \omega)$-thin tubes in $N$ 
such that all Margulis tubes in $N$ belong to $\TT$; in particular $\theta$ is greater than the Margulis constant. 
Let $M$ be a 3-manifold biLipschitz homeomorphic to $N$ and let $M(0)$
be the image of $N \setminus \bigcup_{T \in \TT} Int(T)$ in $M$ under  the 
 biLipschitz homeomorphism $F$. Let $\partial M(0)$ (resp. $\partial M$) denote the  boundary  of $M(0)$ (resp. $M$).

Let $(Q, \partial Q)$ be the unique hyperbolic  pair of pants such that each  component
of $\partial Q$ has length one. $Q$ will be called
the {\it standard} pair of pants.
An isometrically embedded copy of $(Q, \partial Q)$ in $(M(0), \partial M(0))$ will be said to be {\it flat}.

\begin{defn} {\rm A {\bf split level surface} associated to a pants decomposition $\{ Q_1, \cdots , Q_n \}$ of $S$ in $M(0) \subset M$
is an embedding $f : \cup_i (Q_i, \partial Q_i) \rightarrow (M(0), \partial M(0))$ such that \\
1) Each $f (Q_i, \partial Q_i)$ is flat \\
2) $f$ extends to an embedding (also denoted $f$) of $S$ into $M$ such that the interior of each annulus component of
$f(S \setminus \bigcup_i Q_i)$  lies entirely in $F(\bigcup_{T \in \TT} Int(T))$. \\
} \end{defn}

The class   of {\it all} topological embeddings from $S$ to $M$ that agree with a split level surface $f$ 
associated to a pants decomposition $\{ Q_1, \cdots , Q_n \}$ on 
$Q_1 \cup \cdots \cup Q_n$ will be denoted by $[f]$. 

We define a partial order $\leq_E$ on the collection of split level surfaces in an end $E$ of $M$ as follows: \\
$f_1 \leq_E f_2$ if there exist $g_i \in [f_i]$, $i=1,2$, such that $g_2(S)$ lies in the unbounded component of $E \setminus g_1(S)$.

\smallskip

\noindent {\bf Tight geodesics}\\
The {\it complexity} of a compact surface $S=S_{g,b}$ of genus $g$
and $b$ boundary components is defined to be $\xi ( S_{g,b} ) =
3g+b$.

For any simplex $\alpha \in \CC(Y)$, $\gamma_\alpha$ will denote a collection of disjoint simple closed curves on $S$ representing
the (homotopy classes) of vertices of $\alpha$.
A pair of simplices $\alpha,\beta$ in   $\CC(Y)$ are
said to {\em fill } an essential subsurface $Y$ of $S$ if all non-trivial non-peripheral curves in $Y$ have  essential
intersection with at least one of $\gamma_\alpha$ or $\gamma_\beta$, where we assume that representatives $\gamma_\alpha$ and $\gamma_\beta$
have been chosen to intersect each other minimally. 

Given arbitrary simplices $\alpha,\beta$ in $\CC(S)$,
form a regular neighborhood of $\gamma_\alpha\cup\gamma_\beta$, 
and fill in all
disks and one-holed disks to obtain $Y $ which is said to be {\it filled} by $\alpha,\beta$. 

For a subsurface $X\subseteq Z$ let $\boundary_Z(X)$ denote the
{\em relative boundary} of $X$ in $Z$, i.e. those boundary components
of $X$ that are non-peripheral in $Z$.

\begin{definition} {\rm
  Let $Y$ be an essential subsurface in $S$. If $\xi(Y)>4$,  a
  sequence of simplices
  $\{v_i\}_{i\in\II} \subset \CC(Y) $ (where $\II$ is a finite or
  infinite interval in $\Z$) is called } tight {\rm if\\
1) For any vertices $w_i$ of $v_i$ and $w_j$ of $v_j$ where $i\ne
  j$, $d_{\CC_1(Y)}(w_i,w_j) = |i-j|$,\\
2) Whenever $\{i-1,i,i+1\}\subset \II$, $v_i$ represents the relative
  boundary $\boundary_Y F(v_{i-1},v_{i+1})$.\\
If $\xi(Y)=4$ then a tight sequence is  the vertex sequence
of a geodesic in $\CC(Y)$.\\
A } tight geodesic {\rm $g$ in $\CC(Y)$
consists of a tight sequence
$v_0, \cdots , v_n$, and two simplices  in $\CC (Y)$, $\I=\I(g)$ and $\T=\T(g)$,
called its  initial and  terminal markings such that $v_0$ (resp. $v_n$) is a sub-simplex of $\I$ (resp. $\T$). The length of $g$ is $n$. 
$v_i$ is called a  simplex of $g$.
$Y$ is called the }  domain or  support {\rm of $g$ and is
denoted as $Y=D(g)$.  $g$ is said to be  supported in $D(g)$}.
\end{definition}

We denote the obvious linear order in $g$ as $v_i < v_j$ whenever
$i<j$.

A geodesic supported in $Y$ with $\xi(Y)=4$ is called a
$4$-geodesic.

Given a surface $W$ with $\xi(W)\ge 4$ and a simplex $v$ in $\CC(W)$
we say that $Y$ is a
{\em component domain of $(W,v)$} if  $Y$ is a component of
$W\setminus \collar(v)$, where $\collar(v)$ denotes a thin collar neighborhood of the simple closed curves.

If $g$ is a tight geodesic with domain $D(g)$,
we call $Y\subset S$ a {\em component domain of $g$} if
for some simplex $v_j$ of $g$, $Y$ is a component domain of
$(D(g), v_j)$.

\smallskip

\noindent {\bf Hierarchies}\\
The next definition is based on \cite{masur-minsky2}, describing certain special paths in $\PP (S)$ and component domains
associated to them. Paths in $\PP (S)$ will be maps $h$ from intervals $I$ in $\mathbb{Z}$ into $\PP (S)$
such that $h(i), h(i+1)$ are adjacent vertices of $\PP (S)$ for all $i, i+1 \in I$. We reverse the logic of the exposition in
\cite{masur-minsky2} slightly here by defining a hierarchy path in $\PP (S)$ first
and then associating a hierarchy of tight geodesics to it.

\begin{defn} \label{hpa} {\rm  A {\bf hierarchy path} in $\PP (S)$ 
joining pants decompositions $P_1$ and $P_2$  is a path $\rho : [0, n] \rightarrow P(S)$ joining $ \rho (0) = P_1$ to $\rho (n) = P_2$ such that \\
1) There is a collection $\{ Y \}$ of essential, non-annular subsurfaces of $S$, called
       component domains for $\rho$, such that for each component domain $Y$ there is a
       connected interval $J_Y \subset [0, n]$ with $\partial Y  \subset \rho (j)$ for each $j \in  J_Y$.\\
2) For a component domain $Y$, there exists a tight geodesic $g_Y $ supported in $Y$ such that for each
       $j \in J_Y$, there is an $\alpha \in g_Y$ with $\alpha \in \rho (j)$. \\ 
A {\bf hierarchy path} in $\PP (S)$  is a sequence $\{ P_n \}_n$ of pants decompositions of $S$ such that for any
$P_i, P_j \in \{ P_n \}_n$, $i \leq j$,  the  finite sequence  $P_i, P_{i+1}, \cdots , P_{j-1}, P_j$ is a hierarchy path
joining pants decompositions $P_i$ and $P_j$. \\
The collection $H$ of tight geodesics $g_Y $ supported in  component domains  $Y$ of $\rho$ will be called the hierarchy
of tight geodesics associated to $\rho$.} \end{defn}

The notion of hierarchy in Definition \ref{hpa} above is a special case of `hierarchies without annuli'
 described in \cite{masur-minsky2}.
The next definition allows us to associate the extra piece of data coming from
tight geodesics  supported in  component domains  of a
hierarchy path $\rho$  to the hierarchy path $\rho$.

\begin{definition}{\rm 
A {\bf slice } of a hierarchy $H$ associated to a hierarchy path $\rho$ is a set $\tau$ of pairs
$(h,v)$, where $h\in H$ and $v$ is a simplex of $h$, satisfying
the following properties:
\begin{enumerate}
\item A geodesic $h$ appears in at most one pair in $\tau$.
\item There is a distinguished pair $(h_\tau,v_\tau)$ in $\tau$,
  called the bottom pair of $\tau$. We call $h_\tau$ the bottom geodesic.
\item For every $(k,w)\in \tau$ other than the bottom pair, $D(k) $
  is a component domain of $(D(h),v)$ for some $(h,v)\in\tau$.
\end{enumerate}

A resolution of a hierarchy $H$ associated to a hierarchy path $\rho : I \rightarrow \PP (S)$ is a sequence of slices $\tau_i
=\{ (h_{i1}, v_{i1}), (h_{i2}, v_{i2}), \cdots ,  (h_{in_i}, v_{in_i}) \} $ (for $i \in I$, the same indexing set)
such that the set of vertices of the simplices  $\{  v_{i1},  v_{i2}, \cdots ,  v_{in_i} \} $ is the same as the set of the non-peripheral
boundary curves of the pairs of pants in $\rho (i) \in \PP (S)$.}
\end{definition}

\noindent {\bf Minsky Blocks} (Section 8.1 of \cite{minsky-elc1})\\ A tight geodesic in $H$ supported in a component domain of complexity $4$ is called a $4$-geodesic
and an edge of a $4$-geodesic in $H$ is called a $4$-edge.

 Given a 4-edge $e$ in $H$,
let $g$ be the $4$-geodesic containing it, and let $D(e)$ be the
domain $D(g)$. Let $e^-$ and $e^+$ denote the initial and
terminal vertices of $e$. As usual, let $\collar v$ denote a small collar neighborhood of $v$ in $D(e)$.

To each $e$  a Minsky block $B(e)$ is assigned as as follows:
\begin{center}
$B(e) =  (D(e)\times [-1,1]) \setminus ($ {\bf collar}
$(e^-)\times[-1,-1/2)\cup $  {\bf collar} $(e^+)\times(1/2,1])$.
\end{center}
That is, $B(e)$ is the product $D(e)\times [-1,1]$, with
solid-torus trenches dug out of its top and bottom boundaries,
corresponding to the two vertices $e^-$ and $e^+$ of $e$.

 The {\em gluing
boundary} of $B(e)$ is
$$
\boundary_\pm B(e) \equiv (D(e)\setminus\collar(e^\pm)) \times
\{\pm 1\}.
$$
The gluing boundary is  a union of three-holed spheres. The
rest of the boundary is a union of annuli. The top (resp. bottom) gluing boundaries of $B(e)$ are 
$(D(e)\setminus\collar(e^+)) \times \{ 1 \} $  (resp.  $(D(e)\setminus\collar(e^-)) \times \{ - 1 \} $.

\smallskip

\noindent {\bf The Model and the bi-Lipschitz Model Theorem}\\
The following Theorem summarizes 
and paraphrases what we  need in this paper from the bi-Lipschitz Model Theorem of Minsky  \cite{minsky-elc1}
and Brock-Canary-Minsky  \cite{minsky-elc2}.
(See Theorem 8.1 of \cite{minsky-elc1} in particular.)

\begin{theorem} \cite{minsky-elc1} \cite{minsky-elc2} 
Let $N$ be the convex core of a simply or doubly degenerate hyperbolic 3-manifold minus an open neighborhood of the cusp(s).
Let $S$ be a compact surface, possibly with boundary, such that
$N$ is homeomorphic to $S \times [0, \infty )$ or  $S \times \mathbb{R}$ according as 
$N$ is simply or doubly degenerate. 
 There exist $L \geq 1$, $ \theta, \omega, \epsilon, \epsilon_1 > 0$,
  a collection $\TT$ of $(\theta,\omega)$-thin tubes containing all Margulis tubes in $N$,
a  3-manifold  $M$,  and an $L$-biLipschitz homeomorphism $F: N \rightarrow M$ 
such that the following holds. \\
Let $M(0) = F(N \setminus \bigcup_{T \in \TT} Int(T))$ and let $F(\TT)$ denote the image of the collection $\TT$ under $F$.
Let  $\leq_E$ denote the partial order on the collection of split level surfaces in an end $E$ of $M$. Then there exists a sequence
$S_i$ of split level surfaces associated to pants decompositions $P_i$ exiting $E$ such that  

\begin{enumerate}
\item $S_i \leq_E S_j$ if $i \leq j$.
\item The sequence $\{ P_i \}$ is a hierarchy path in $\PP (S)$.
\item If $P_i \cap P_j = \{ Q_1, \cdots Q_l \}$ then $f_i(Q_k)=f_j(Q_k)$ for $k=1 \cdots l$, where $f_i, f_j$ are the embeddings
defining the split level surfaces $S_i, S_j$ respectively.
\item For all $i$, $P_i \cap P_{i+1} = \{ Q_{i,1}, \cdots Q_{i,l} \}$ consists of a collection of $l $ pairs of pants,
 such that $S \setminus  (Q_{i,1} \cup \cdots \cup Q_{i,l})$
has a single non-annular component of complexity $4$.
Further, there exists a Minsky block $W_i$ and an isometric map $G_i$ of $W_i$ into $M(0)$ such that $f_i(S \setminus
 (Q_{i,1} \cup \cdots \cup Q_{i,l})$ (resp. $f_{i+1}(S \setminus
 (Q_{i,1} \cup \cdots \cup Q_{i,l})$) is contained in the bottom (resp. top) gluing boundary of $W_i$. 
 \item For each flat pair of pants $Q$ in a split level  surface $S_i$ there exists an isometric embedding of $Q \times [-\epsilon, \epsilon]$
into $M(0)$ such that the embedding restricted to $Q \times \{ 0 \}$  agrees with $f_i$ restricted to $Q$.
\item For each $T\in \TT$, there exists a split level surface $S_i$ associated to pants decompositions $P_i$
such that the core curve of $T$ is isotopic to a non-peripheral boundary curve of $P_i$. The boundary $F(\partial T)$ of $F(T)$ with the induced
metric $d_T$ from  $M(0)$ is a
Euclidean torus equipped with a product structure $S^1 \times S^1_v$,
where any circle of the form $S^1 \times \{ t \} \subset S^1 \times S^1_v$
is a round circle of unit length and is called a horizontal circle; and any circle of the form 
 $\{ t \} \times  S^1_v$ is a round circle of length $l_v$ and is called a vertical circle.
\item Let $g$ be a tight geodesic other than the bottom geodesic
 in the hierarchy $H$ associated to the hierarchy path $\{ P_i \}$, let $D(g)$ 
be the support of $g$ and let $v$ be a boundary curve of $D(g)$. Let $T_v$ be the tube in $\TT$ such that
the core curve of $T_v$ is isotopic to $v$. If a vertical circle of $(F(\partial T_v), d_{T_v})$ has length $l_v$ less than
$n\epsilon_1$, then the length of $g$ is less than $n$. \end{enumerate}
\label{bilipmodel} \end{theorem}

Since the above statement is culled out of a large amount of material, particularly from \cite{minsky-elc1},
we give specific references here.\\
$M(0)$ (resp. $M$) above is denoted by $M_\nu (0)$ (resp. $M_\nu$) in Section 8 of \cite{minsky-elc1}.\\
The collection $F(\TT)$ is denoted by $\UU$ in \cite{minsky-elc1} and is called the set of tubes in $M_\nu$.\\
The hierarchy $H$ in Item (7) of Theorem \ref{bilipmodel} is constructed in Lemma 5.13 of \cite{minsky-elc1} (see also
Theorem 4.6 of \cite{masur-minsky2}) and the hierarchy path of Item (2) is obtained from it by constructing a {\it resolution}
sweeping through it in Lemma 5.8 of \cite{minsky-elc1}. (We have thus reversed the logical order of hierarchies and
hierarchy paths in our treatment.)\\ The estimate on the length of $g$ 
in Item (7) of Theorem \ref{bilipmodel} comes from Equation 9.6 of \cite{minsky-elc1}  which gives estimates on meridian coefficients.\\
The Euclidean structure of $F(T)$ for $T \in \TT$ in Item (6)
comes from gluing together the internal blocks (as well as boundary blocks) described in Section 8.1  
and in Theorem 8.1 of \cite{minsky-elc1}. \\ Theorem 8.1 of  \cite{minsky-elc1}   further describes the construction
of {\it split level surfaces} and Items (1), (3)  and (4) follow from it. \\
Item (5) simply ensures the existence of uniform product neighborhoods and follows from the fact that Minsky blocks are glued
by isometries on their $3$-holed sphere boundary components. In fact $\epsilon = \frac{1}{4}$ suffices. \\
Finally \cite{minsky-elc2} ensures that the model constructed in \cite{minsky-elc1} is indeed  biLipschitz  
homeomorphic to $N$.

We use the notation of Theorem \ref{bilipmodel} in the rest of this subsection, fixing $N, M$.

\begin{lemma} Given $l > 0$ there exists $n \in \natls$ such that the following holds. \\Let $v$ be a vertex in the hierarchy $H$ such that
the length of the core curve of the Margulis tube $T_v$ corresponding to $v$ is greater than $l$. 
Next suppose $(h,v) \in \tau_i$ for some slice $\tau_i$ of the hierarchy $H$ such that $h $
is supported on $Y$, and $D$ is a component of $Y \setminus \collar v$. Also suppose that $h_1 \in H$
such that $D$ is the support of $h_1$. Then the length of $h_1$ is at most $n$.
\label{abut}
\end{lemma}

\begin{proof} 
Let $\alpha$ be a meridian curve on $F(\partial T_v)$ such that $F^{-1}(\alpha )$ bounds a totally geodesic
disk in $T_v$.

By Item (6) of Theorem \ref{bilipmodel},
 $F(\partial T_v)$ is a metric product $S^1 \times S^1_v$. Choose horizontal and vertical curves $\alpha_h, \alpha_v$ on 
$F(\partial T_v)$. Then $\alpha$ is homologous to $(n \alpha_h + \alpha_v)$ for some integer $n$. Hence $l(\alpha ) \geq l_v$,
where $l(\alpha )$ is the length of  $\alpha$ and $l_v$ denotes the length of the vertical circle.
Since $F$ is an $L$-biLipschitz homeomorphism by Theorem \ref{bilipmodel}, it follows that $l(F^{-1}(\alpha )) \geq \frac{l_v}{L}$.
Let $\Delta$ be the totally geodesic disk bounded by $F^{-1}(\alpha )$. Then the radius $r_v$ of  $\Delta$ is bounded below by
$sinh^{-1} (\frac{l_v}{L})$. Let $l_c$ denote the length of the core curve $c_v$ of $T_v$. Then any geodesic on $\partial T_v$ homotopic 
to $c_v$ in $T_v$ has length bounded below by $\frac{l_v}{L} l_c$.

Also $l(F^{-1}(\alpha_h )) \leq L$ and $F^{-1}(\alpha_h )$ is homotopic to the core curve  $c_v$.

Hence $\frac{l_v}{L} l_c \leq L$. It follows that $l_v \leq \frac{L^2}{l_c} \leq \frac{L^2}{l}$.
The Lemma now follows from Item (7) of Theorem \ref{bilipmodel}. \end{proof}

One last fallout of the Minsky model (Theorem 8.1 of \cite{minsky-elc1} again) that we shall need is the following.

\begin{lemma}
Given $l > 0$ and $n \in \natls$, there exists $L_2 \geq 1$ such that the following holds:\\
Let $S_i, S_j$ ($i<j$) be split level surfaces associated to pants decompositions $P_i, P_j$ such that \\
a) $(j-i) \leq n$ \\
b) $P_i \cap P_j$ is a (possibly empty) pants
decomposition of $S \setminus W$, where $W$ is an essential (possibly disconnected) subsurface of $S$ such that
each component $W_k$ of $W$ has complexity $\xi (W_k) \geq 4$.\\
c)For any $k$ with $i < k < j$, and 
$(g_D, v) \in \tau_k$ for $D \subset W_i$ for some $i$, no curve in $v$ has a geodesic realization in $N$ of length less than $l$. \\
Then there exists an $L_2$-biLipschitz embedding  $G: W \times [-1,1] \rightarrow M$, such that \\
1) $W$ admits a hyperbolic metric given by $W = Q_1 \cup \cdots \cup Q_m$ where each $Q_i$ is a flat pair of pants. \\
2) $ W \times [-1,1]$ is given the product metric. \\
3) $f_i(P_i \setminus P_i \cap P_j) \subset W \times \{ -1 \}$ and $f_j(P_j \setminus P_i \cap P_i) \subset W \times \{ 1 \}$.
\label{thick} \end{lemma}

\noindent {\bf Idea of Proof:}
What the Lemma above says is that if a `thick' piece of the manifold $N$ is trapped between split level surfaces $S_i, S_j$,
then it is biLipschitz to a product region on the support of the hierarchy path between $S_i, S_j$.
This follows from the construction of the model in Section 8.2 of \cite{minsky-elc1} along with Equation 9.6 of  \cite{minsky-elc1}.
The lower bound on lengths of hierarchy curves in hypothesis (c) ensures an {\it upper bound} on the twist coefficient
($[h_v] $ in Equation 9.6 of  \cite{minsky-elc1}) exactly as in the proof of Lemma \ref{abut}. Hence the `full hierarchy' path
(including annuli in the sense of \cite{masur-minsky2})
 between $S_i, S_j$ equipped with markings is of length bounded  in terms of $l, n$. This guarantees the existence
of a biLipschitz product region as required. Since this is the only place where we shall require full hierarchies and twist coefficients
in this paper, and since the rest of the proof of Lemma \ref{thick} follows Lemma \ref{abut}, we omit the details,
referring the interested reader to Section 8.2 of \cite{minsky-elc1}. (See also \cite{masur-minsky2} where a quasi-isometry
is constructed between the mapping class group and the full marking complex. Interpreted in these latter terms there is 
a
bounded length element in the mapping class group $MCG(W)$  taking the marking on $S_i \cap W$ to 
the marking on $S_j \cap W$.)
$\Box$

\section{Split Geometry}\label{sec:split}

\subsection{Constructing Split Level Surfaces} \label{splitgeometry}
The aim of this subsection is to extract a special sequence of split level surfaces from
 the sequence of split level surfaces constructed in Theorem \ref{bilipmodel}. The main point is to ensure that successive
split level surfaces are separated by a definite amount in $M(0)$.
We continue with the notation of Theorem \ref{bilipmodel} in this subsection.

Fix an $l > 0$. The precise value of $l$ will be less than the Margulis constant for
hyperbolic 3-manifolds and will be determined by the  Drilling Theorem \ref{drill} to
be used in the next subsection. We shall henceforth refer to Margulis
tubes that have core curve of length $\leq l$ as {\bf thin Margulis
  tubes} and the corresponding vertex $v$ as a {\bf thin vertex}.

For convenience start with a  doubly degenerate surface group. Let $\rho (i) = \{ P_i \}$ be a hierarchy path provided by Item (2) of
Theorem \ref{bilipmodel}. Let  $H$ be the hierarchy of tight geodesics associated to $\{ P_i \}$
and $\cdots , \tau_{i-1}, \tau_i ,
\tau_{i+1}, \cdots $ be a resolution. 
Let $S_i$ be the split level surface corresponding to $P_i$ and let $\tau_i$ be the slice
whose vertices comprise the curves in $P_i$. Let $S_i^s$ denote the collection of flat pairs of pants occurring in the image of $S_i$ in
$M(0)$. The metric on the model manifold and the induced path metric on
$M(0)$ will be denoted by $d_M$ and will be called the {\it model metric}. Thus $S_i$ is an embedding and $S_i^s$ is the image in $M(0)$
of a collection of pairs of pants.

\begin{definition}  A curve $v$ in $H$ is {\bf $l$-thin} if the core curve of the Margulis tube $T_v$ has length less than or equal to $l$. \\
A curve $v$ is said to split a pair of split level surfaces $S_i$ and $S_j$ ($i<j$) if $v$ occurs as a vertex
in both $\tau_i$ and $\tau_{j-1}$.\\
A pair of split level surfaces $S_i$ and $S_j$ ($i<j$) is said to be an {\bf $l$-thin pair} if there exists an   $l$-thin curve $v$ such that
$v$ occurs as a vertex
in  both $\tau_i$ and $\tau_{j-1}$.\\
A pair of split level surfaces $S_i$ and $S_j$ ($i<j$) is said to be an {\bf $l$-thin pair on a component
domain $D$}  if  \\
a) $P_i \cap P_j$ is a pants decomposition of 
$S \setminus D$, none of whose curves are $l$-thin.  \\
b) There exists a tight geodesic $g_D \in H$ supported on $D$ such that $(g_D, u) \in \tau_k$
for all $i <  k < j$, where the multicurve $u$ contains an $l$-thin curve. (Here $D$ could be $S$ itself.) Further we demand that
the initial and final vertices  of $g_D$ consist of curves contained in (the boundary curves of) 
$P_i, P_j$ respectively. 

A pair of split level surfaces $S_i$ and $S_j$ ($i<j$) is said to be an {\bf $l$-thick pair}
(or an $l$-thick pair on $S$) if no curve $v \in \tau_k$ is  $l$-thin for $i < k < j$.
\end{definition}

In fact in criterion (b) of the definition of an $l$-thin pair on a component
domain $D$, we might as well have assumed that the initial and final markings of $g_D$, ${\bf I}(g_D)$ and ${\bf T}(g_D)$
respectively, are precisely $P_i, P_j$.
This is the case when  the markings are {\it complete} in the sense of \cite{minsky-elc1}.

 By Definition \ref{hpa},  Item (1), the set   $ J(v) = \{ i:  v \in \rho (i) \}$
 is an interval. Consider the family of
intervals $\{ J(v): v \in g_H \}$, where $g_H$ is the distinguished
main geodesic (bottom geodesic) for the hierarchy $H$. Then
$\bigcup_v \{ J(v) : v\in g_H  \} = \mathbb{Z}$. This follows from

the fact that each $\tau_i$ has a simple closed curve corresponding to some
vertex in $g_H$.

 Any pair $v_i, v_{i+1}$ of simplices (multicurves) which
form  successive vertices of the base
geodesic $g_H$ are at a distance of $1$ from each other by
tightness of $g_H$.

\smallskip

\noindent {\bf Selecting Split Level Surfaces}\\
We shall now construct a subset $\mathbb{I}$ of $\mathbb Z$ by selecting a subsequence of the elements  $\{ P_i \}$
of the hierarchy path.
Let $\tau_{m_i}$ be the first slice in the resolution such that $(g_H,v_i) \in \tau_{m_i}$.
Let ${\mathbb{I}}_1 = \{ m_i : i \in \mathbb{Z} \}$. We shall now expand the set ${\mathbb{I}}_1$ if necessary as follows.

If some curve in $v_i$ is $l$-thin, then we declare that $[m_i, m_{i+1}] \cap \mathbb{I} = \{ m_i, m_{i+1} \}$, i.e.
no integer strictly between $m_i, m_{i+1}$ is added to ${\mathbb{I}}_1$.

More generally, for any $j \in \mathbb Z \setminus {\mathbb{I}}_1 $, choose $i$ such that $m_i <j < m_{i+1}$.

Then $j \in {\mathbb{I}}_2$
if and only if there exists $k$ \\
a) either with $j < k \leq m_{i+1}$ such that $S_j, S_k$ form an $l$-thin pair on some component domain $D$.\\
b) or with $k < j \leq m_{i+1}$ such that $S_k, S_j$ form an $l$-thin pair on some component domain $D$.

Finally set $\mathbb{I} = \mathbb{I}_1 \cup {\mathbb{I}}_2$. Then $\mathbb{I} = \{ \cdots , n_{i-1}, n_i, n_{i+1}, \cdots \}$
inherits  a linear order  from $\mathbb Z$ such that $j < k$ implies $n_j < n_k$.

Note that the same construction works for simply degenerate groups if we replace $\mathbb Z$ by
$\mathbb N$.

The next few Propositions identify some of the features of the selection $\mathbb{I}$. The main point is to show that the sequence
of split level surfaces $S_{n_i}$ makes definite progress out an end.

\begin{prop} Let $\mathbb{I} = \{ \cdots , n_{i-1}, n_i, n_{i+1}, \cdots \}$ be as above.
There exists  a positive integer $N_0$ such that for all $i$,\\
a) either  $(S_{n_i},S_{n_{i+1}})$ is an $l$-thin pair on some component domain $D$ \\
b) or $(S_{n_i},S_{n_{i+1}})$ is an $l$-thick pair and $n_{i+1} - n_i \leq N_0$.
\label{selection1} \end{prop}

\begin{proof} Suppose that $(S_{n_i},S_{n_{i+1}})$ is not an $l$-thin pair on some component domain.

Then, by the construction of $\mathbb{I} $ and Lemma \ref{abut}, there exists $N_1 (= N_1(l))$ such that
for all $k$ with $n_i < k < n_{i+1}$, if $\tau_k 
=\{ (h_{k1}, v_{k1}), (h_{k2}, v_{k2}), \cdots ,  (h_{km_k}, v_{km_k}) \} $, the length of the tight geodesic $h_{ki}$
satisfies $l(h_{ki}) \leq N_1$. Further, none of the curves in $v_{ki}$ are $l$-thin, ensuring $l$-thickness
of the pair $(S_{n_i},S_{n_{i+1}})$.

Note that $m_k \leq \xi (S)$ where $\xi (S)$ is the complexity of $S$. Also the number of component domains in $W \setminus \collar (v)$
for $W = D(h_{ki})$  is certainly bounded above by the number of pairs of pants in a pants decomposition of $S$ and hence
by $\xi (S)$. Therefore $(n_{i+1} - n_i)$ is bounded above by $N_1^{\xi (S)} \times \cdots \times N_1^{\xi (S)}$ (${\xi (S)}$ times).
Choosing $N_0 = N_1^{{\xi (S)}^2}$ we are done. \end{proof}

The next Proposition asserts that between two successive split level
surfaces $S_{m_i}$ and $S_{m_{i+1}}$ selected from the base geodesic, our selection process
 `interpolates' a uniformly bounded number of new split level
surfaces. Equivalently the cardinality of the set $({\mathbb{I}}_2 \cap [m_i, m_{i+1}])$ is uniformly bounded.

\begin{prop} Let $\mathbb{I} = \{ \cdots , n_{i-1}, n_i, n_{i+1}, \cdots \}$ and
${\mathbb{I}}_1 = \{ \cdots , m_{i-1}, m_i, m_{i+1}, \cdots \}$  be as above. 
There exist  a positive integer $N_2$ such that for all $i$, if $n_j = m_i$ and $n_k = m_{i+1}$
then $k - j \leq N_2$.
\label{selection2} \end{prop}

\begin{proof}  Let $k \in \mathbb{I}$. So
$S_k$ is a split level surface  interpolated between $S_{m_i}$ and $S_{m_{i+1}}$ for some
$k$ with $m_i < k < m_{i+1}$. Let the corresponding slice $\tau_k 
=\{ (h_{k1}, v_{k1}), (h_{k2}, v_{k2}), \cdots ,  (h_{km_k}, v_{km_k}) \} $. Then there exists
a unique `subslice' $\tau_k^0
=\{ (h_{k1}, v_{k1}), (h_{k2}, v_{k2}), \cdots ,  (h_{kr_k}, v_{kr_k}) \} $, with $r_k \leq m_k$ such that
the length of the tight geodesic $h_{ki}$
satisfies $l(h_{ki}) \leq N_1$ for all $i \leq r_k$ and $S_k$ is a split level surface.

Since the total number of such choices is bounded above by $N_0$ by the proof of Proposition \ref{selection1},
and for each such choice at most two (by the construction of ${\mathbb{I}}_2$ above) split level surfaces are introduced,
it follows that the total number of $l$-thin split level surfaces $S_k$ with $m_i < k < m_{i+1}$ is bounded above by $2N_0
=2N_1^{{\xi (S)}^2}$. Choosing $N_2= 2N_0$ we are done. \end{proof}

\begin{lemma}
There exists $n$ such that each thin curve splits at most
$n$ split level surfaces in the  sequence $\{ S_{n_i} : i \in \mathbb{I} \}$. \label{Margulis-gqc}
\end{lemma}

\begin{proof} Since, for any $i$ the number of split level 
surfaces $S_{n_i} $ between $S_{m_i}$ and  $S_{m_{i+1}}$ ($m_i,  m_{i+1} \in {\mathbb{I}}_1$)
 is at most $N_2$  by Proposition \ref{selection2},
it suffices to prove that any thin curve splits a uniformly bounded
number of   $S_{m_i}$'s.  

If a curve $v$  splits both $S_{m_i}$ and  $S_{m_{j}}$,
then  $v$ belongs to both the pants
decomposition $P_{m_i}$ and  $P_{m_{j}-1}$. 

Suppose $(g_S, v_{m_k}) \in \tau_{m_k}$ for $k=i, j$, where $g_S$ denotes the bottom geodesic of the hierarchy $H$.
Then the distance
between  $v_{m_i}$ and $v_{m_{j}}$ in $\CC (S)$ 
is at most 3 by tightness, i.e. $|i-j| \leq 3$. 
  Taking $n = 3 N_2$, we are through. \end{proof}

\noindent {\bf Pushing split level Surfaces Apart } \\
We shall now use Item (5) of Theorem \ref{bilipmodel} to `thicken' each of the $S_{n_i}$'s if necessary, so that successive split 
level surfaces can be arranged to be
uniformly separated.  Recall that $S_i^s$ is the collection of flat embedded pairs of pants in $M(0)$ corresponding to the split level surface $S_i$.

\begin{defn}
A pair of split level surfaces $S_i$ and $S_j$ ($i<j$) is said to be {\bf $k$-separated} if \\
a) for all $x \in S_i^s$, 
$d(x,S_j^s) \geq k$\\
b)Similarly, for all $x \in S_j^s$, $d(x,S_i^s) \geq k$. \end{defn}

\begin{lemma} \label{selection3} Let $\mathbb{I} = \{ \cdots , n_{i-1}, n_i, n_{i+1}, \cdots \}$ be as above.
There exist  $k_0 > 0$ and a sequence of split level surfaces $\Sigma_{i}$ and a  positive integer $N_0$ such that for all $i$,
$(\Sigma_i,\Sigma_{i+1})$ is $k_0$-separated and\\
a) either  $(\Sigma_i,\Sigma_{i+1})$ is an $l$-thin pair on some component domain $D$ \\
b) or $(\Sigma_i,\Sigma_{i+1})$ is an $l$-thick pair and $n_{i+1} - n_i \leq N_0$.
 \end{lemma}

\begin{proof} By Proposition \ref{selection1}, the sequence  $\{ S_{n_i} \}_i$ satisfies one of the alternatives (a) or (b).
It remains to modify  $\{ S_{n_i} \}_i$ such that $(S_{n_i},S_{n_{i+1}})$ are $k_0$-separated for some $k_0 >0$ and all $i$.

By Item (5) of Theorem \ref{bilipmodel}, there exists $\epsilon > 0$ such that for all flat pairs of pants $Q_i$ in $S^s_{n_i}$ there exists
an isometric embedding  $H_i: Q_i \times [-\epsilon,\epsilon]$ into $M(0)$.

Also, by Lemma \ref{Margulis-gqc}, there exists $n \in \natls$ such that if $Q_k \in P_{n_i} \cap P_{n_j}$, then $|i-j| \leq n$.
Further the collection $\{i\in \mathbb{I} : Q_k \in P_{n_i} \}={\mathbb{I}}_{Q_k}$ is an interval in $\mathbb{Z}$.
Let $\epsilon_2 = \frac{\epsilon}{n}$. If ${\mathbb{I}}_{Q_k} = [a_k, b_k] \subset \Z$, define
$Q_{ks} = H_k|_{Q_k \times s\epsilon_2}$. Note that $(b_k-a_k) \leq n$.

For each embedding $f_{n_i}$ defining the split level surface $S_{n_i} $, and $Q_k \in P_{n_i} \cap P_{n_j}$ for some $j \neq i$,
let ${\mathbb{I}}_{Q_k} = [a_k, b_k]$ and $s = (n_i - a_k)$. Then
define $f_{n_i}^{\prime}|_{Q_k} = H_k|_{Q_k \times s\epsilon_2}$.

Now, let $\Sigma_i$ be the split level surface defined by $f_{n_i}^{\prime}|_{Q_k}$, whenever $Q_k \in P_{n_i}$.
Choosing $k_0=\epsilon_2$, it follows that successive split level surfaces are $k_0$-separated. \end{proof}

Let $\TT_l$ denote the collection of tubes in $\TT$ whose core curves have length less than $l$.
Also let $M(l) = M(0) \bigcup_{T \in \TT \setminus \TT_l} F(T)$ denote the union of $M(0)$ and all $l$-thick tubes.

\begin{defn} {\rm An $L$-biLipschitz {\bf split  surface} in $M(l)$ associated to a pants decomposition $\{ Q_1, \cdots , Q_n \}$ of $S$
and a collection $\{ A_1, \cdots , A_m \}$ of complementary annuli  in $S$ 
is an embedding $f : \cup_i Q_i \bigcup  \cup_i A_i \rightarrow M(l)$ such that\\
1) the restriction  $f: \cup_i (Q_i, \partial Q_i) \rightarrow (M(0), \partial M(0))$ is a split level surface. \\
2) the restriction $f: A_i \rightarrow M(l)$ is an $L$-biLipschitz embedding.\\
3)  $f$ extends to an embedding (also denoted $f$) of $S$ into $M$ such that the interior of each annulus component of
$f(S \setminus (\cup_i Q_i \bigcup  \cup_i A_i))$  lies entirely in $F(\bigcup_{T \in \TT_l} Int(T))$.}\end{defn}

\noindent {\bf Note:} The difference between a split level surface and a split surface is that the latter may contain
biLipschitz annuli in addition to flat pairs of pants.

\smallskip

Let
$\Sigma_{i}^{s}$ denote the union of the collection of flat pairs of pants
and biLipschitz annuli in the image of the embedding  $\Sigma_{i}$.

\begin{theorem} 
Let $N, M, M(0), S, F$ be as in Theorem \ref{bilipmodel} and $E$ an end of $M$. For any $l$ less than the Margulis constant,
let $M(l) = \{ F(x) : {\rm injrad_x} (N) \geq l \}$. Fix a hyperbolic metric on $S$ such that each component of $\partial S$ is 
totally geodesic of length one (this is a normalization condition).
 There exist $ L_1 \geq 1$, $  \epsilon_1 > 0$, $n \in \natls$, 
 and a sequence $\Sigma_i$ of $L_1$-biLipschitz, $  \epsilon_1$-separated split  surfaces exiting the end $E$ of $M$
such that for all $i$, one of the following occurs: 
\begin{enumerate}
\item An $l$-thin curve splits the pair $(\Sigma_i ,\Sigma_{i+1})$, i.e. the associated split level surfaces form
an $l$-thin pair. 
\item there exists an $L_1$-biLipschitz embedding  $$G_i: (S\times [0,1], (\partial S)\times [0,1]) \rightarrow (M, \partial M)$$
such that $\Sigma_i^s = G_i (S\times \{ 0\})$ and $\Sigma_{i+1}^s = G_i (S\times \{ 1\})$
\end{enumerate}
Finally, each $l$-thin curve in $S$ splits at most
$n$ split level surfaces in the  sequence $\{ \Sigma_{i} \}$. \label{wsplit}
\end{theorem}

\begin{proof} By Lemma \ref{selection3}, there exists $k>0$, a  positive integer $N_0$ and
a sequence of {\it split level surfaces} $\Sigma_{i}^0$  such that for all $i$,
$(\Sigma_i^0,\Sigma_{i+1}^0)$ is $k$-separated and\\
a) either  $(\Sigma_i^0,\Sigma_{i+1}^0)$ is an $l$-thin pair on some component domain $D$ \\
b) or $(\Sigma_i^0,\Sigma_{i+1}^0)$ is an $l$-thick pair and $n_{i+1} - n_i \leq N_0$.\\
(We add the superscript $0$ to indicate that we are still dealing with split level surfaces and not split surfaces.)

In Case (a), there exists an $l$-thin curve splitting the pair of split level surfaces $(\Sigma_i^0,\Sigma_{i+1}^0)$.

In Case (b),  
let $P_{n_i}, P_{n_{i+1}}$  be  the pants decompositions  associated to $\Sigma_i^0,\Sigma_{i+1}^0$ and let $P_{n_i} \cap P_{n_{i+1}}$ be 
a (possibly empty) pants
decomposition of $S \setminus W$, where $W$ is an essential (possibly disconnected) subsurface of $S$ such that
each component $W_k$ of $W$ has complexity $\xi (W_k) \geq 4$. Hence
by Lemma
\ref{thick}, 
 there exists $L_2 \geq 1$ and an $L_2$-biLipschitz embedding  $G: W \times [-1,1] \rightarrow M$, such that \\
1) $W$ admits a hyperbolic metric given by $W = Q_1 \cup \cdots \cup Q_m$ where each $Q_i$ is a flat pair of pants. \\
2) $ W \times [-1,1]$ is given the product metric. \\
3) $f_{n_i}(P_{n_i} \setminus P_{n_i} \cap P_{n_{i+1}}) \subset W \times \{ -1 \}$ and 
$f_{n_{i+1}}(P_{n_{i+1}} \setminus P_{n_i} \cap P_i) \subset W \times \{ 1 \}$.

Also, from the proof of Lemma \ref{selection3}, there exists $\epsilon > 0$
such that for all $i$, there exits an isometric embedding $H_{n_i}: 
(P_{n_i} \cap P_{n_{i+1}}) \times [0, \epsilon] \rightarrow M$ such that 
$H_{n_i} (P_{n_i} \cap P_{n_{i+1}}) \times \{ 0\} \subset f_{n_i}(P_{n_i} )$ and 
$H_{n_i} (P_{n_{i+1}} \cap P_{n_{i+1}}) \times \{ \epsilon\} \subset f_{n_{i+1}}(P_{n_{i+1}})$.

Finally since $(\Sigma_i^0,\Sigma_{i+1}^0)$ is an $l$-thick pair, there exists standard annuli $A_1, \cdots A_p$,
$L_3 = L_3(l)\geq 1$, $\epsilon_1 >0$ and $L_3-$biLipschitz 
embeddings $\Gamma_j: A_j\times [-1,1] \rightarrow \bigcup_{T \in \TT \setminus \TT_l} F(Int(T)) $ such that \\
a) $S=\cup_k P_k \bigcup \cup_j A_j$ is the union of the pairs of pants above along with the annuli $A_j$.\\
b) $f_{n_i}$ restricted to  $A_j$ agrees with $\Gamma_j$ restricted to $A_j\times \{-1\}$.\\
c) $f_{n_{i+1}}$ restricted to  $A_j$ agrees with $\Gamma_j$ restricted to $A_j\times \{1\}$.\\

Pasting these  maps \\
i) $G: W \times [-1,1] \rightarrow M$, \\

ii) $f_{n_i}(P_{n_i} \setminus P_{n_i} \cap P_{n_{i+1}}) \subset W \times \{ -1 \}$,\\
iii) $H_{n_i}: 
(P_{n_i} \cap P_{n_{i+1}}) \times [0, \epsilon] \rightarrow M$, and\\
iv) $\Gamma_j: A_j\times [-1,1] \rightarrow \bigcup_{T \in \TT \setminus \TT_l} F(Int(T)) $ \\
along the common boundaries, we
obtain an $L_1$-biLipschitz embedding  $$G_i: (S\times [0,1], (\partial S)\times [0,1]) \rightarrow (M, \partial M)$$
such that the {\it split surfaces} $\Sigma_i^s = G_i (S\times \{ 0\})$ and $\Sigma_{i+1}^s = G_i (S\times \{ 1\})$.

Lemma \ref{Margulis-gqc} now proves the final assertion. \end{proof}

Pairs of split surfaces satisfying Alternative (1) of Theorem \ref{wsplit} will be called an $l$-thin pair of split surfaces (or simply a thin
pair if $l$ is understood). Similarly, pairs of split surfaces satisfying Alternative (2) of Theorem \ref{wsplit} will be called an $l$-thick pair

(or simply a thick
pair) of split surfaces.

\begin{rmk} \label{wsplitrmk} {\rm 
The notion of split surface could be made a bit more general. We might as well  require a split surface to be
a (uniformly) biLipschitz embedding of a bounded geometry subsurface of $S$ containing a pants decomposition.
Theorem \ref{wsplit} then summarizes the consequences of the Minsky model that  we shall need in this paper.
We have thus constructed  the following from  the Minsky model:\\
1) A sequence of split surfaces $S^s_i$ exiting the
end(s) of $M$, where $M$ is marked with a homeomorphism to $S \times J$ ($J$ is $\mathbb{R}$ or $[0, \infty )$ according as $M$ is totally or simply degenerate). 
$S_i^s \subset S \times \{ i \}$. \\
2) A collection of Margulis tubes $\mathcal{T}$ in $N$ with image  $F(\mathcal{T})$ in $M$
(under the biLipschitz homeomorphism between $N$ and $M$). We refer to the elements of 
$F(\mathcal{T})$ also as Margulis tubes.\\
3) For each  complementary annulus of $S^s_i$ with core $\sigma$,
  there is a Margulis tube $T\in \mathcal{T}$ whose core is freely homotopic to $\sigma$
   such that $F(T)$ intersects $S^s_i$ at the boundary. (What this roughly
  means is that there is an $F(T)$ that contains the complementary
  annulus.) We say that $F(T)$ splits $S^s_i$.\\
4) There exist constants $\epsilon_0 >0, K_0 >1$ such that
for all $i$, either there exists a Margulis tube splitting both $S^s_i$
  and $S^s_{i+1}$, or else $S_i (=S^s_i)$ and $S_{i+1}  (=S^s_{i+1})$ have injectivity radius bounded
below by $\epsilon_0$ and bound a {\bf thick block} $B_i$, where a thick block is defined to
be a $K_0-$biLipschitz homeomorphic image of $S \times I$. \\
5) $F(T) \cap S^s_i$ is either empty or consists of a pair of
  boundary components of $S^s_i$ that are parallel in $S_i$. \\
6) There is a uniform upper bound $n = n(M)$ on the number of surfaces that
 $F(T)$ splits. \\
For easy reference later on, a model manifold satisfying conditions (1)-(6) above is said
to have {\bf weak split geometry}.}
\end{rmk}

We have isolated the features of weak split geometry in Remark \ref{wsplitrmk} above so as to emphasize the point that it
is possible to make a definition
independent of the Minsky model and the hierarchy machinery.\footnote{This was our original approach to the
main Theorem of this paper: Prove it for more and more general model geometries, e.g. bounded geometry \cite{mahan-bddgeo},
i-bounded geometry \cite{mahan-ibdd}, amalgamation geometry and split geometry \cite{mahan-amalgeo}. Finally prove that the Minsky
model satisfies split geometry \cite{mahan-split}.} This will be useful for easy referencing in
 \cite{mahan-elct} and \cite{mahan-series2}. In fact a strengthening of weak split geometry will be enough to guarantee
the existence of Cannon-Thurston maps as we shall see below.

\subsection{Split Blocks}

\begin{defn} Let $(\Sigma_i^s, \Sigma_{i+1}^s)$ be a thick pair of split surfaces in $ M$. 
The closure of the bounded component of
$M \setminus (\Sigma_i^s \cup \Sigma_{i+1}^s)$ will be called a thick block.\end{defn}

Note that a thick block is uniformly biLipschitz to the product $S \times [0,1]$ and that its boundary components are 
$\Sigma_i^s, \Sigma_{i+1}^s$.

\begin{defn} Let $(\Sigma_i^s, \Sigma_{i+1}^s)$ be an $l$-thin pair of split surfaces in $M$
and $F(\TT_i)$ be the collection of $l$-thin Margulis tubes that split both $\Sigma_i^s, \Sigma_{i+1}^s$. The closure of the union of the
bounded components of
$M \setminus ((\Sigma_i^s \cup \Sigma_{i+1}^s)\bigcup_{T\in \TT_i} F(T))$ will be called a split block.
Equivalently, the closure of the union of the
bounded components of
$M(l) \setminus (\Sigma_i^s \cup \Sigma_{i+1}^s)$ is a split block.\end{defn}

Topologically, a  split block $B^s $
is a topological product $S^s \times I$ for some {\em not necessarily connected 
  $S^s$}. However, its upper and lower boundaries need not be $S^s
\times 1$ and $S^s \times 0$. We only require that the upper and
lower boundaries
be {\it  split subsurfaces} of $S^s$. This is to allow for Margulis tubes
starting (or ending) within the split block. Such tubes would split
one of the horizontal boundaries but not both. We shall call such
tubes {\bf hanging tubes}. Connected components of split blocks
are called {\bf split components}. By $l$-thinness,
 there is a {\it non-empty} collection of $l$-thin Margulis tubes, called {\bf splitting tubes},
splitting a split block. For each splitting tube $F(T)$ of a split block $B^s $, the intersection $(B^s \cap F(T)) \subset M$ 
is called the {\bf vertical boundary} of the splitting tube. Note that the 
vertical boundary of a splitting tube is the union of two disjoint annuli.

See figure below, where the left split component has four hanging tubes
and the right  split component has two hanging tubes. The vertical space between the components is the place where an 
$l$-thin Margulis tube splits the
split block into two split components. 

\begin{center}

\includegraphics[height=4cm]{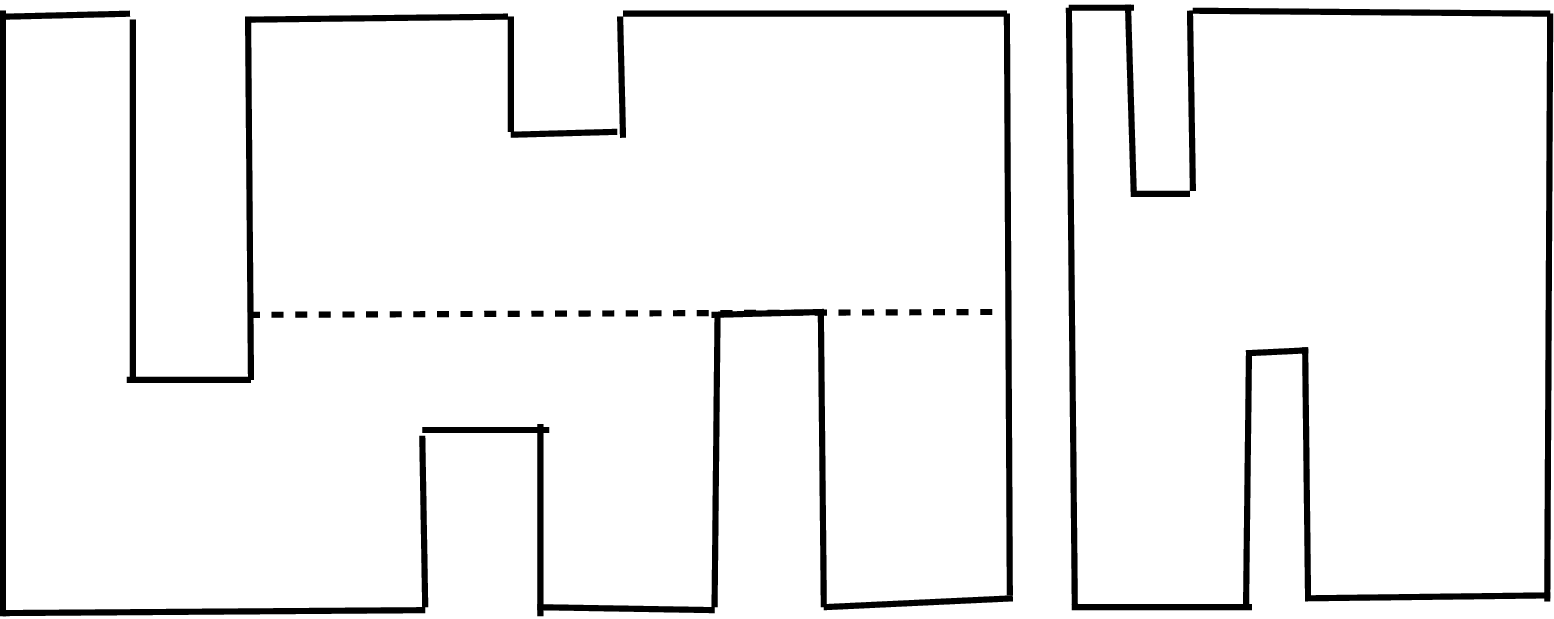}
\underline{Figure:  {\it Split Components of Split
Block with hanging tubes} }

\end{center}

Observe further that for each hanging tube $F(T)$, there exists a split surface $S^s$ (marked with a dotted line in the figure)
which intersects the boundary $F(\partial T)$ non-trivially and such that $S^s$ contains an annulus whose core-curve is homotopic 
(in $M$) to the core curve of $F(T)$. Also, the closure of $(F(\partial T) \setminus S^s)$ consists of precisely two annuli
called the {\bf vertical boundary} of the hanging tube. We can assume further that\\
a) $S^s \cap F(\partial T)$ is a biLipschitz
annulus called the {\bf horizontal boundary} of $F(T)$ in the split block $B^s$.\\
b) the union of the vertical and horizontal boundaries of an $l-$thin hanging tube $F(T)$ in $B^s$  is precisely equal to
$F(T) \cap B^s$.

Note that the whole manifold $M$ is  the union of \\
a) Thick blocks (biLipschitz homeomorphic to $S \times I$) \\
b) Split blocks (homeomorphic to $S^s \times I$ for some split surfaces)\\
c) $l$-thin Margulis tubes. \\

Also note that the union of thick and split blocks is $M(l)$, which is the complement (in $M $) of the
union  of $l$-thin Margulis tubes. Each of these Margulis tubes
splits a uniformly bounded number of split blocks and might end in a hanging tube.

\subsection{Electrocutions} \label{electcns}

For any hanging tube or splitting tube $F(T_j)$ in a split block $B^s$ (with top and bottom split surfaces
$\Sigma^s_k, \Sigma^s_{k+1}$, say), 
let $A_{ji}=S^1 \times [0,l_{ji}]$, $(i=1,2)$, be the vertical boundaries ($i=1,2$ correspond to the left and right vertical annuli
in the previous Figure). 
Let the metric product
$S^1 \times [0,1]$ be called the {\bf standard annulus}
 for splitting tubes. For hanging tubes the standard annulus
 will be  $S^1 \times [0,1/2]$. 

 We shall define a {\bf welded split block} ${B_{wel}}$ (homeomorphic to $S \times [0,1]$) to be a split block with
identifications on vertical boundaries of splitting tubes and hanging tubes.
Let $H_{ji}^0: [0,l_{ji}] \rightarrow [0,x]$ be the unique linear surjective map (scaling) taking $0$ to $0$ and $l_{ji}$ to $x$,
where $x$ is $1$ or $1/2$ according as $F(T_j)$ is a splitting tube or a hanging tube.
Now define $H_{ji}: S^1 \times [0,l_{ji}] \rightarrow S^1 \times [0,x]$ by $H_{ji}(y,z) = (y, H_{ji}^0(z))$. Finally
extend $H_{j1} \cup H_{j2}$ continuously to the horizontal boundaries $S^1 \times [-\epsilon,\epsilon]$ of hanging tubes $F(T_j)$ 
as Lipschitz maps to 
$S^1 \times \{ p \}$ by $H_{j}(y,z) = (y,p)$ where $p$ is either $0$ or $1/2$ according as the horizontal boundary of 
$F(T_j)$ lies at the bottom or the top of the hanging tube (for instance in the figure, the horizontal boundary marked with
a dotted line lies at the top of a hanging tube).
Now glue the mapping
cylinders of $H_{j1} \cup H_{j2} \cup H_j$ (for hanging tubes) and $H_{j1} \cup H_{j2}$ 
 (for splitting tubes)  to  $F(\partial T_j)\cap B^s$
to obtain the  {\bf welded split block} ${B_{wel}}$. Note that ${B_{wel}}$ is homeomorphic to $S \times [0,1]$.
The images of the standard annuli in ${B_{wel}}$ after the identification shall simply be called {\it standard annuli in ${B_{wel}}$}.

For each hanging tube, there exists one distinguished curve on either $\Sigma^s_k$ or $ \Sigma^s_{k+1}$.
When the hanging tube intersects $ \Sigma^s_{k+1}$, this is the image of $S^1 \times \{ 1/2 \}$ contained in the standard
annulus after identification. Similarly, 
when the hanging tube intersects $ \Sigma^s_{k}$, this is the image of $S^1 \times \{ 0 \}$ contained in the standard
annulus after identification. Again, for each splitting tube, there exists two
 distinguished curves, one each on  $\Sigma^s_k$ and $ \Sigma^s_{k+1}$ -
the images of $S^1 \times \{ 0, 1 \}$ contained in the standard
annulus after identification. Such
simple closed curves shall be called {\bf weld curves}.
The resulting metric on ${B_{wel}}$ will be denoted by ${d_{wel}}$.

We shall equip ${B_{wel}}$ with a new pseudometric.
Equip the standard annulus $S^1 \times [0,x]$ (where $x$ is $1$ or $1/2$) with the product of the zero metric on the $S^1$-factor
and the Euclidean metric on the $[0,x]$ factor.  Let $(S^1 \times [0,x], d_0)$ denote the resulting pseudometric.
The {\bf tube-electrocuted metric} $d_{tel}$   is defined to be the pseudometric metric that
 agrees with ${d_{wel}}$ away from the standard annuli in ${B_{wel}}$ and with $d_0$ on 
the standard annuli in ${B_{wel}}$. To distinguish it from $({B_{wel}}, {d_{wel}})$
we shall represent the new space and the pseudometric on it by $(B_{tel},d_{tel})$.
Note that the underlying topological spaces 
${B_{wel}}$ and $B_{tel}$ are the same and homeomorphic to $S \times [0,1]$.

Recall that in defining thick blocks, $S$ was equipped with a fixed hyperbolic metric.
If $\Sigma^s_k$ is the  bottom  split surface of the split block $B^s_k$ and also 
the  top  split surface of a (thick or split) block block $B^s_{k-1}$, then the common  split surface $\Sigma^s_k$
can be easily extended over complementary annuli to a common uniformly biLipschitz embedding of $S$ into welded blocks $B_{wel,k}$
and $B_{wel,k-1}$, where we define $B_{wel,m} = B_m$ for thick blocks. When $B_{k-1}$ is thick, this follows from the fact
that the complementary annuli are uniformly biLipschitz embeddings of $S^1 \times [0,1]$. When
$B^s_{k-1}$ is split, the mapping cylinder construction above restricted to $\Sigma^s_k$ is the same whether $\Sigma^s_k$
is regarded as the  bottom  split surface of  $B^s_k$ or
the  top  split surface of  $B^s_{k-1}$.
We shall continue to denote the extended split surface by $\Sigma^s_k$ and call it a split surface in $B_{wel,k}$.
From now on shall drop the suffix $wel$ from (thick or split) blocks $B_{wel,k}$ and denote them simply as $B_k$.
Note that {\it all such extended split surfaces} are homeomorphic to $S$ via uniformly biLipschitz homeomorphisms.

\smallskip

\noindent {\bf The Welded Model Manifold}\\
Gluing successive welded blocks along common split surfaces we obtain the {\bf welded model manifold} 
$({M_{wel}}, {d_{wel}})$ homeomorphic to
$S \times J$, where $J = \mathbb{R}$ or $[0, \infty )$ according as the original manifold $N$ is doubly or simply degenerate.

It remains to construct the tube-electrocuted pseudometric $d_{tel}$ on ${M_{wel}}$.
The tube electrocuted metrics on successive welded split blocks  coincide on the common split surface.
The same is clearly true if the successive blocks are thick.

If a weld curve lies in $\Sigma^s_k= B_{wel,k} \cap B_{wel,k-1}$ and precisely one of $B_{wel,k}, B_{wel,k-1}$ is a thick block,
we   fix the convention that for the tube electrocuted metric $d_{tel}$ on  ${M_{wel}}$: {\it All weld curves have length zero.}

Gluing successive tube electrocuted blocks using the convention above, we obtain the {\bf tube electrocuted manifold} $(M_{tel}, d_{tel})$.
Observe that the underlying topological manifolds ${M_{wel}}$ and $M_{tel}$ are the same. (The notation $({M_{wel}}, {d_{wel}})$ and 
$(M_{tel}, d_{tel})$ is used to distinguish the metrics.)

The union of the images of the contiguous mapping cylinders of maps $H_{j1} \cup H_{j2} \cup H_j$ (or
$H_{j1} \cup H_{j2}$) in $(M_{tel}, d_{tel})$ 
associated to a particular 
$l$-thin Margulis tube $T$ (and hence $F(T)$) is topologically a solid torus $T^t$. Equipped with the tube electrocuted metric,
$(T^t, d_{tel})$
is of diameter at most $n$ by Theorem \ref{wsplit}. The collection of all $T^t$'s in  
$(M_{tel}, d_{tel})$  is denoted $\TT^t$. (We shall continue to use the same notation $\TT^t$ for the collection of 
$T^t$'s in  $(M_{wel}, d_G)$ to be defined below.)

The images of split components $K$ of $B^s$ in $B_{tel}$ will continue to be called split components of $B_{tel}$.
A lift $\til{K}$ of a split component $K$ of $(B_{tel}, d_{tel})$ to the universal cover 
$(\til{B_{tel}}, d_{tel})  $
shall be termed a {\bf split component}  of $\til{B_{tel}}$.

Let $d_G$ be the (pseudo)-metric
obtained by electrocuting the collection $\KK$ of split components $\til{K}$ in  $(\til{B_{tel}}, d_{tel})
\subset (\til{M_{tel}}, d_{tel})$ as $({B_{tel}}, d_{tel})$ ranges over all split blocks.
$d_G$ will be called the the {\bf graph metric} on $\til{M_{tel}} (=\til{M_{wel}}$). Thus $(\til{M_{wel}}, d_G)$ is isometric to
$\EE (\til{M_{wel}},\KK)$ with the electric metric.

\begin{rmk} {\bf Alternate Description:} {\rm There is an alternate description of a pseudometric on $\til M$ which makes it quasi-isometric to 
$(\til{M_{wel}}, d_G)$.  For each lift $\til{K} \subset
\til{M}$ of a split component $K$ of a split block of $M(l) \subset M$, 
there are lifts of $l$-thin Margulis tubes that share the boundary of  $\til{K}$ in $\til{M}$. Adjoining these lifts to 
$\til{K}$ we obtain {\bf extended split components}. Let $\KK^\prime$ denote the collection of extended split components  in $\til{M}$.
We continue to denote  the collection of split components in $\til{M(l)} \subset \til{M}$ by $\KK$.
Let $\til{M(l)}$ denote the lift of $M(l)$ 
to $\til M$. 
Then the inclusion of $\til{M(l)}$ into $\til{M}$ gives a quasi-isometry between $\EE (\til{M(l)}, \KK)$ and
$\EE (\til{M}, \KK^\prime)$ equipped with the respective electric metrics. This  follows from the last assertion of Theorem \ref{wsplit}.

 Again, by the last assertion of Theorem \ref{wsplit},
the inclusion of $\til{M(l)}$ into $\til{M_{wel}}$ gives a quasi-isometry between $\EE (\til{M(l)}, \KK)$ and 
$\EE (\til{M_{wel}}, \KK) (= (\til{M_{wel}}, d_G))$.

Therefore $(\til{M_{wel}}, d_G)$ is quasi-isometric to 
$\EE (\til{M}, \KK^\prime)$. We shall henceforth identify $\EE (\til{M}, \KK^\prime)$ with $(\til{M_{wel}}, d_G)$ via this
quasi-isometry without explicitly mentioning the quasi-isometry. The electric metric on $\EE (\til{M}, \KK^\prime)$
shall therefore be denoted by $d_G$ also.
We shall find it
easier to use $\EE (\til{M}, \KK^\prime)$ when dealing with all of $\til M$, whereas $(\til{M_{wel}}, d_G)$ will be more useful
when dealing with the block structure of $\til{M_{wel}}$.} \label{alt} \end{rmk}

\begin{rmk} {\rm Here is the raison d'etre for the two closely related but different electric spaces.
In the ladder construction of Section \ref{ladder} below, it is important that a split surface goes `all the way
across', i.e. is an embedded copy of $S$. There is no canonical way to do this in the model manifold $M$. In fact
for the ladder construction of Section \ref{ladder} to work, it is important that a split surface
in  
$({M_{wel}}, {d_{wel}})$ is an embedded copy of $S$ having uniformly bounded geometry.
 This is simply not possible in $M$ as Margulis tubes may be arbitrarily thin. 
On the other hand, we finally need to control hyperbolic geodesics in $\til{N}$ by means of the ladder.
Since $M$ is biLipschitz to $N$, we can equivalently control them in $\til M$. 
The {\it Alternate Description}
above establishes a way of transferring this control from $(\til{{M_{wel}}}, {d_{G}})$ to $\til M$, which is where we really want the
control on geodesics.}
\end{rmk}

The following definition illustrates this passing back and forth between these two quasi-isometric electric spaces.

\begin{definition} Let $Y \subset \til{N}$ and $X=F(Y)$.  $X \subset \til{M}$ is said to 
be $\Delta$-graph quasiconvex if for any hyperbolic geodesic $\mu$ joining $a, b \in Y$,
$F(\mu )$ lies inside $N_\Delta (X, d_G) \subset \EE (\til{M}, \KK^\prime)$. \end{definition}

For $X$ a split component, define $CH(X) = F(CH(Y))$, where $CH(Y)$
  is the convex hull of $Y$ in $\til{N}$.
Then $\Delta$-graph quasiconvexity of $X$
is equivalent to the condition that  $dia_G (CH(X))$  is bounded by $\Delta^{\prime} = \Delta^{\prime}(\Delta )$ as any split component has diameter one in $(\til{M_{tel}}, d_G)$.

\subsection{Quasiconvexity of Split Components} \label{sec-gqc}

We now proceed to show further that split components are
   quasiconvex (not necessarily uniformly) in the hyperbolic metric,
  and uniformly quasiconvex in the graph metric, i.e. we need to
  show {\em hyperbolic quasiconvexity} and {\em uniform graph
  quasiconvexity} of split components.

\smallskip

\noindent {\bf Hyperbolic Quasiconvexity:}

Let $N = {\Hyp}^3/\Gamma$ 
be a complete hyperbolic 3-manifold. Then \cite{Thurstonnotes} there exists a geometrically finite hyperbolic manifold with 
compact convex core
$N_{gf}$ and a strictly type-preserving embedding $i$ of $N_{gf}$ into $N$, which is a homotopy equivalence.
Then for any boundary component $S^h$ of $N_{gf}$, $i_\ast (\pi_1(S^h)) \subset \pi_1(N)$ is called a peripheral subgroup. 
In the Theorem below $\pi_1(N)$ will be identified with a Kleinian group $\Gamma$ and the peripheral subgroup $i_\ast (\pi_1(S^h))$
with a Kleinian subgroup of $\Gamma$.

\begin{theorem} {\bf Covering Theorem} \cite{Thurstonnotes}
\cite{canary-cover} Let $N = {\Hyp}^3/\Gamma$ 
be a complete hyperbolic 3-manifold. A finitely generated subgroup
${\Gamma}'$ is geometrically infinite if and only if it  contains
a finite index subgroup of a geometrically infinite peripheral
subgroup.
\label{cover}
\end{theorem}

We shall now specialize
the Thurston-Canary covering theorem  \ref{cover} to the case under
consideration, namely, infinite index free subgroups of surface
Kleinian groups. 

\begin{lemma} Let $N$ be a simply or doubly degenerate hyperbolic 3-manifold homotopy equivalent to a surface equipped with
a weak split geometry model $M$.
For $K$ a split component, let $\tilde{K}$ be a lift to $\til N$.  Then there exists $C_0
= C_0(K)$ such that  the convex hull of $\tilde K$ minus cusps lies in a
$C_0$-neighborhood of $\tilde K$ in $\til N$. \label{hypqc}
\end{lemma}

\begin{proof} 
 Let $\Gamma = \pi_1(N)$, and ${\Gamma}^{\prime}= i_\ast (\pi_1(K)) (\subset \Gamma )$.

Then $\Gamma$ itself is the unique peripheral subgroup.
Since ${\Gamma}^{\prime}$ has infinite index in $\Gamma$, it follows from Theorem \ref{cover} that 
 ${\Gamma}^{\prime}$ is geometrically finite. The result follows. (Cusps need to be excised because the model manifold is biLipschitz
homeomorphic to $N$ minus cusps.) \end{proof}

\noindent {\bf Graph Quasiconvexity:} \\
Next, we shall prove that each split component  is uniformly graph
quasiconvex. We begin with the following Lemma. Recall that we are
dealing with simply or totally degenerate groups without accidental
parabolics.

\begin{lemma} Let $\Sigma$ be a component of a proper extended split subsurface $S^s_i$
 of $S$.
Any (non-peripheral) simple closed curve in $S$ appearing
in the hierarchy whose free homotopy class 
has a representative lying in $\Sigma$ must have a geodesic
representative in $M$ lying  within a uniformly
bounded distance of $S^s_i$ in the graph metric $d_G$. \label{hier-qc}
\end{lemma}

\begin{proof} Suppose a  curve $v$ in the hierarchy is
homotopic into $\Sigma$.  Then $v$ is at a distance of at most $1$ in the curve
complex from each of   the boundary components of $\Sigma$. Since $\Sigma$ is a proper subsurface of $S$, the relative
boundary $\partial_S(\Sigma) \neq \emptyset$.  Let
$\alpha$ be such a boundary component. Next, suppose that the geodesic realization of $v$ in $N$
intersects some block $B^s_j$ (via the correspondence in Alternate Description \ref{alt}). Then $v$ must be at a distance of
at most one from some
 curve $\sigma $ in the base geodesic $g_H$ forming an
  element of the pants decomposition of the split surface $S^s_j$.
  
By
 tightness, the distance from $\alpha$ to $\sigma$ in the curve
 complex is at most $2$. Hence the distance ($\leq |i-j|$) of $S^s_j$ from
 $S^s_i$ (in the $d_G$ metric) is
$\leq 2n$ from Lemma \ref{Margulis-gqc}. Therefore $v$ is realized
within a distance $2n$ of $S^s_i$ in the graph metric $d_G$.
\end{proof}

Recall (Definition 8.8.1 of \cite{Thurstonnotes}) that a {\bf pleated surface}  in a hyperbolic three-manifold $N$ is
a complete hyperbolic surface $S$ of finite area, together with an isometric map $h: S \rightarrow N$ 
such that every $x \in S$ is in the interior of some geodesic segment in $S$ which
is mapped by $h$ to a geodesic  segment in $N$. Also, $h$ maps cusps to cusps. We refer the reader to Section 8.8
of \cite{Thurstonnotes} for further details. A pleated surface is said to be incompressible if
$h_\ast : \pi_1(S) \rightarrow \pi_1(N)$ is injective. A standard fact about hyperbolic surfaces and pleated surfaces is Lemma \ref{pleated-bdd}
below.
(See the proof of Proposition 8.8.5 of \cite{Thurstonnotes} for instance.)

An {\bf $l$-thin annulus} on a hyperbolic surface $S^h$ is a maximal connected component of the set $\{ x \in S^h: injrad_x(S^h) < \frac{l}{2} \}$.
This is the 2-dimensional analogue of an $l$-thin Margulis tube. Note that an  $l$-thin annulus 
may also be 
a neighborhood of a cusp in $S^h$.

\begin{lemma} \label{pleated-bdd}  \cite{Thurstonnotes} \cite{bonahon-bouts}
For all $l > 0$ and $g, n \in \natls$ there exists $\Delta=\Delta (l, g, n) >0$ such that the following holds.\\
Let $S^h$ be any hyperbolic surface of genus $g$ and $n$ boundary components and/or cusps.
Let $\AAA_l$ be the collection of $l$-thin annuli. Then $\EE (S^h , \AAA_l )$ has diameter
less than $\Delta$ in the electric metric. \\
Again, let $N$ be a hyperbolic 3-manifold and let $\TT_l$ be the collection of $l$-thin Margulis tubes and cusps in it. Let $h: S \rightarrow N$ 
be an incompressible  pleated surface. Then $h(S)$ has diameter less than $\Delta$ in the electric metric on $\EE (N, \TT_l )$.\end{lemma}

Next, we show that any (non-peripheral) simple closed curve $v_i$ in $S^s_i$
(not just hierarchy curves as in Lemma \ref{hier-qc})
 must be realized within a
uniformly bounded distance of $S^s_i$ in the graph metric. In fact we shall show further
that any pleated surface which contains at least one boundary geodesic
of $\Sigma$ in its pleating locus lies within a
uniformly bounded distance of $S^s_i$ in the graph metric.

\begin{lemma}
There exists $B > 0$ such that the following holds:\\
Let $\Sigma$ be a proper split subsurface of $S^s_i$. Then any pleated
surface with at least one boundary component coinciding with a
geodesic representative of a non-peripheral component of $\partial \Sigma$ must lie
within a $B$-neighborhood of $S^s_i$ in $(M, d_G) = (\EE(\til{M}, \KK^{\prime}))/\Gamma$, where $\KK^{\prime}$ denotes the collection of extended split components in $\til M$
and $\Gamma$ is the fundamental group of $M$ regarded as the group of deck transformations of $\til M$. In
particular, every simple closed curve in $S$ homotopic into $\Sigma$
has a geodesic representative
within a $B$-neighborhood of $S^s_i$ in $(M, d_G) $.
\label{scc-qc}
\end{lemma}

\begin{proof}
Choose a curve $v_i$ homotopic to a simple closed curve on
 $\Sigma$. Let $\alpha$ denote its
geodesic
realization in $N$. 

Let $\Sigma_p$ be any pleated (sub)surface 
 whose boundary
coincides with the geodesics representing the boundary components of
$\Sigma$. (See \cite{Thurstonnotes} for the construction of such pleated surfaces.)
In particular, we may choose $\Sigma_p$ such that its
pleating locus contains  $v_i$. 
 Since $S^s_i$ is a split surface in $M$, 
the topological type of $\Sigma_p$ has finitely many possibilities.
By Lemma \ref{pleated-bdd} the diameter of 
$\Sigma_p$ is bounded by $\Delta=\Delta (l)$ in the electric metric on $\EE (N, \TT_l )$.
Since the $l$-thin components of the  boundary of $\Sigma_p$ are contained in $l$-thin tubes bounding $S^s_i$, it follows that 
$\Sigma_p$ (and $\alpha$ in particular)
 lies in a $\Delta$ neighborhood of $S^s_i$ in the electric metric on $\EE (N, \TT_l )$.
Since each $T \in \TT_l$ is contained in the image of some $K \in \KK^{\prime}$ under the quotient map $(\EE(\til{M}, \KK^{\prime}))
\rightarrow (\EE(\til{M}, \KK^{\prime}))/\Gamma = (M, d_G)$ the result follows. \end{proof}

\begin{rmk}  In \cite{bowditch-model}, Bowditch indicates a
method to obtain a related (stronger) result that given $B_1 > 0$,
there exists $B_2 > 0$ such that any two simple
closed curves realized within a Hausdorff distance $B_1$ of each
other in $M$ are within a  distance $B_2$ of each other in the curve
complex.\end{rmk}

\subsection{Drilling and Filling} In this subsection we summarize some material that will be needed in Section \ref{sec-gqc1}
to prove uniform graph quasiconvexity of split components.

The {\it
 Drilling Theorem} of Brock and Bromberg
 \cite{brock-bromberg-density}, which built on work of Hodgson and
 Kerckhoff \cite{HK-cone} \cite{HK-univlbds} is given below. 
We shall invoke a
 version of this theorem which is closely related to one used by
 Brock and Souto in \cite{brock-souto}.

\begin{theorem} \cite{brock-bromberg-density}
 For each
$L
>1$, and $n $ a positive integer,
there is an $\ell >0$ so that if $N_{gf}$ is a geometrically finite
hyperbolic 3-manifold and $c_1, \cdots c_n$ are geodesics in $N_{gf}$
with length $\ell_{N_{gf}}(c_i) < \ell$ for all $c_i$, then there is an
$L$-biLipschitz diffeomorphism of pairs
$$h \colon (N_{gf} \setminus \cup_i \tube(c_i), \cup_i \boundary \tube(c))  \to (N_{gf}^0 \setminus
\cup_i \cusp(c_i), \cup_i \boundary \cusp(c_i))$$ where
$N_{gf}\setminus \cup_i \tube(c_i)$ denotes
 the complement of a standard tubular neighborhood of $\cup_i c_i$ in $N_{gf}$,
$N_{gf}^0$ denotes the complete hyperbolic structure on $N_{gf} \setminus
\cup_i c_i$, and $\cusp(c_i)$ denotes a standard rank-2 cusp
corresponding to $c_i$. \label{drill}
\end{theorem}

$N_{gf}^0$ is said to be obtained from $N_{gf}$ by {\bf drilling}.
We remark here (following 
\cite{brock-bromberg-density}) that the drilled manifold is the
unique hyperbolic manifold which has the same conformal structure
on its domain of discontinuity, but has core curves of $\T$
drilled out to give rank $2$ parabolics.

The {\it Filling Theorem}
 of Thurston \cite{Thurstonnotes} (generalized by Canary
\cite{canary-cover})  we shall require is stated below.

\begin{theorem}  \cite{Thurstonnotes} \cite{canary-cover}
Given any quasifuchsian surface group $\Gamma$ and $N = {\Hyp}^3/\Gamma$
there exists $\delta > 0$ depending only on the Euler characteristic of the surface
such that for all $x \in CC(N)$, the convex core of $N$, there exists
a pleated
 surface $\Sigma$ such that $d(x, \Sigma ) \leq \delta$.
\label{filling}
\end{theorem}

\subsection{Proof of Uniform Graph-Quasiconvexity}\label{sec-gqc1}

We
need to prove the uniform graph quasiconvexity of split components.

Let $B^s$ be a split block with a splitting $l$-thin Margulis tube $T$. We aim at

showing:

\begin{prop} {\bf Uniform Graph Quasiconvexity of Split Components:}
Each (extended) split component  $\til K$ is uniformly 
graph-quasiconvex in  $(\til{M}, d_G)$.
\label{gr-qc}
\end{prop}

The proof of Proposition \ref{gr-qc} will occupy  this entire
subsection.

Let
 $B^s \subset B = S \times I$ be a split block with horizontal boundary consisting of split surfaces $S^s_j, S^s_{j+1}$. 
Let $\bigcup_i T_i$ be the union of $l$-thin
Margulis tubes splitting $B^s$ (we suppress the dependence on the index $j$ for the time being). Let $K$ be a split component. Then $ K = (S_1 \times
I)$ topologically for a subsurface $S_1$ of $S$. Also, let $\partial_s K
= \partial S_1 \times I$ denote the collection of boundary annuli of $K$ that abut the splitting tubes. Let $\partial S_1 = \bigcup_i \sigma_i =
\sigma$ be the finite collection of boundary curves. $\sigma $ is
thus a {\bf multicurve}. Each $\sigma_i$ is homotopic to the core
curve of an $l$-thin splitting Margulis tube $T_i$. Let $\T = \bigcup_i T_i$. $\T$
will be referred to as a {\bf multi-Margulis tube}.

We have already shown in Lemma \ref{hypqc}  that $\pi_1 (S_1 )
\subset \pi_1 (S) $ includes into $\pi_1(N)$ as a geometrically finite subgroup
of $PSl_2 (C)$. Let $N_1$ be the cover of $N$ corresponding to
$\pi_1 (K) = \pi_1 (S_1 )$. Then $N_1$ is geometrically finite.
Let $\T^1$ be the multi-Margulis tube in $N_1$ that consists of
tubes that are (individually) isometric to individual components
of the multi-Margulis tube $\T$.

Let $N_{1d}$ be the hyperbolic
manifold obtained from $N_1$ by drilling out the core curves of
$\T^1$.  Since $N_1$ is geometrically
finite, so is $N_{1d}$.

We first observe that the boundary of the augmented Scott core $X$
of $N_{1d}$ is incompressible away from cusps. To see this, note
that $X$ is double covered by a copy of $D \times I$ with solid
tori drilled out of it, where $D$ is the {\it double} of $S_1$
(obtained by doubling $S_1$ along its boundary circles).

Identify $X$ with the convex core $CC(N_{1d})$ of
$N_{1d}$. We also identify $D$ with the convex core boundary.
Since $D$ is incompressible away from cusps, we have the following.

\begin{lemma} (Chapter 8, \cite{Thurstonnotes})
$D$ is a pleated surface.
\label{D-pleated}
\end{lemma}

Since $N_1$ is the cover of $N$ corresponding to $\pi_1 (K)
\subset \pi_1 (N)$, $K$ lifts to an embedding into $N_1$. Adjoin
the multi-Margulis tube $\T^1$ to (the lifted) $K$ to get an
augmented split component $K_1$. Let $K_{1d} \subset N_{1d}$
denote $K_1$ with the components of $\T^1$ drilled. We want to
show that $D$ lies within a uniformly bounded distance of $K_{1d}$
in the lifted graph metric on $N_{1d}$. This would be enough to
prove  a version of Proposition \ref{gr-qc} for the drilled
manifold $N_{1d}$ as the split geometry structure gives rise to a
graph metric on $N$, hence a graph metric on $N_1$ and hence
again, a graph metric on $N_{1d}$. Finally, we shall use the
 Theorem \ref{drill}
to complete the proof of Proposition \ref{gr-qc}.

\begin{lemma}
There exists $C_1$ such that for any split component $K$, $D$
lies within a  uniformly bounded neighborhood of $K_{1d}$
in $N_{1d}$. \label{D-gr-qc}
\end{lemma}

\begin{proof} {\underline \em Case 1}:
 $D \cap K_{1d} \neq \emptyset$ \\
If $D$ intersects $K_{1d}$, then the Lemma follows directly from Lemma \ref{pleated-bdd}: Incompressible pleated surfaces have bounded diameter in the graph
metric $d_G$.

\smallskip

\noindent  {\underline \em Case 2}: 
$D \cap K_{1d} = \emptyset$ \\
This is the more difficult case because  a priori $D$ might
lie far from $K_{1d}$. Recall that $F: N \rightarrow M$ is a biLipschitz homeomorphism between the hyperbolic manifold and the model
manifold. Let $M_1 = F(N_1)$.
Let $B$ denote the block (split or thick) in the {\it model manifold} $M$ 
 containing $F(K)$. Let $B_1 \subset M_1$
denote its lift to  $M_1$. Let $B_{1d}$ denote $B_1$ with $\T^1$
drilled.

Then $B_{1d} - F(K_{1d})$ is topologically a disjoint union of
`vertically thickened flaring annuli' $F(A_i)$, say. Each $A_i \subset F^{-1}(B_1)$ is of the
form $S^1 \times [0,\infty )$ where $S^1 \times \{ 0 \}$ lies on 
 $T_i$. 

More elaborately, what this means is the following.
Identifying $B$ with $S \times I$, we may identify $B_1$ with
$S_1^a \times I$, where $S_1^a$ is the cover of $S$ corresponding
to the subgroup $\pi_1 (S_1) \subset \pi_1 (S)$. Then $S_1^a$ may
be regarded as $S_1$ union a finite collection of flaring annuli
$F(A_i)$ (one for each boundary component of $S_1$). Thus $B_1$ is
the union of a core   $F(K_1)$ and a collection of {\it vertically
thickened flaring annuli} of the form $F(A_i) \times I$. Hence
$B_{1d}$ is the union of a core   $F(K_{1d})$ and the collection of
{\it vertically thickened flaring annuli}  $F(A_i) \times I$. Also
the boundary $\partial A_i = A_i \cap T_i$ is a curve of fixed
length $\epsilon_0$.  Let us fix one such annulus $A_1$. 
 Refer figure below (where we have removed subscripts for convenience): \\

\begin{center}

\includegraphics[height=4cm]{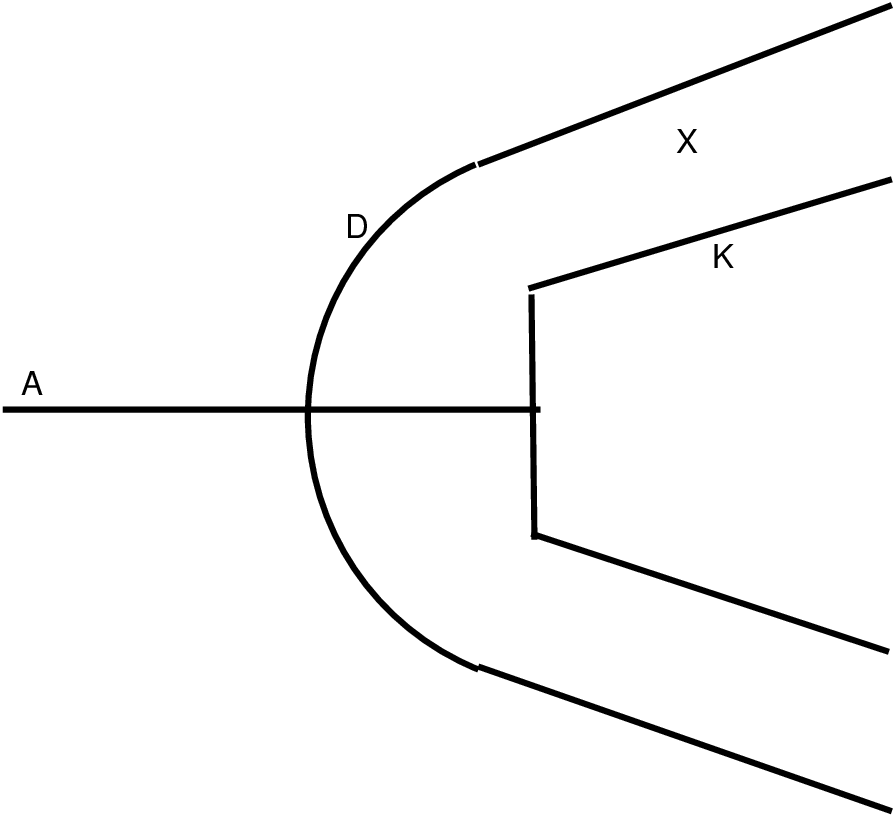}

\underline{Figure:  {\it Graph Quasiconvexity} }

\end{center}

Recall that $D$ bounds $X$ and $X$ contains $K_{1d}$. Thicken the convex core slightly to $N_\epsilon (X)$ such that
its boundary, $D_\epsilon$ is a smooth surface. 

Let $\til{M_1}$ denote the cover of $M_1$ corresponding to $i_\ast \pi_1(A_1)$, where $i$ denotes the inclusion map. Let $\til{D_\epsilon}$
denote the lift of $D_\epsilon$ to $\til{M_1}$. Then each lift $A_1 \times \{t\}$ separates $\til{D_\epsilon}$ since $i_\ast (\pi_1(A_1)) \subset \pi_1(D)$ is a subgroup
such that the cover $\til D$ has two ends.

Hence, by a small homotopy of $A_1$, we can assume that \\
a) $F^{-1}(F(A_1) \times I)$
is a smooth manifold (with boundary) biLipschitz homeomorphic to $(F(A_1) \times I)$.\\
b) $F(D_\epsilon)$ is transverse to each $(F(A_1) \times \{ t \})$ for $t$ belonging to an interval $I_1$ of some fixed length $h_0 > 0$
(equal to the uniform lower bound on the height of
split blocks $B^s$)
contained in $I$.\\

Since $A_1 \times \{t\}$ separates $\til{D_\epsilon}$ in $\til{N_1}$, it follows that for each $t \in I_1$, 
$F(D_\epsilon)$ intersects $(F(A_1) \times \{ t \})$
 in an essential loop $\alpha_t$ parallel to
$\partial A_1$. Hence $D_\epsilon$ must contain an annulus of the form
$\alpha \times I_1 \subset F(A_1) \times I_1$. Also the
 the length of $I_1$  is at least
 $\frac{h_0}{L}$, where $L$ is the biLipschitz constant for $F$. Since this is true for all $\epsilon$, it is also true for the pleated surface $D$. Hence for at least some $t \in I$, the length
of $\alpha \times \{ t \}$ is uniformly bounded (by $\frac{2L\pi
(4g-4)}{h_0}$) by the Gauss-Bonnet Theorem applied to $D$. Much more is true in fact, but this is enough for
our purposes.

Since $\alpha \times \{ t \} \subset A_1 \times \{ t \}$ and the
latter is an exponentially  flaring annulus, it follows that 
there exist uniform constants
$C_0 > 0, \eta > 1$ such that if $d(\alpha \times \{ t \} , \partial (A_1 \times \{ t \}) \geq d_0$
then the length
of $\alpha \times \{ t \}$ is bounded below by $C_0 \eta^{d_0}$. 

These two estimates imply that
there is some point
$ p \in \alpha \times \{ t \} \subset D$ such that $d(p, T_1
)$ is uniformly bounded (in terms of the genus of $S$ and the
minimal height of split blocks $h_0$), where $T_1$ is the drilled Margulis tube intersecting $A_1 \times \{ t \}$ non-trivially.

By Lemma \ref{pleated-bdd} and using the biLipschitz homeomorphism $F$ between $N_{1d}$ and $M_{1d}$, the diameter of $D$
is uniformly bounded in the graph metric lifted to $M_{1d}$.
Hence, by the triangle inequality, $D$ lies in a uniformly bounded
neighborhood of $K$ in the graph metric (using either of the descriptions of the graph metric in Remark \ref{alt}). \end{proof}

\noindent {\bf An Alternate Proof of Lemma \ref{D-gr-qc}:} 
A simpler proof of the fact that $D$
lies in a uniformly bounded neighborhood of $K_1$ in the graph
metric may alternately be obtained directly as follows. First,
that $M_1$ is geometrically finite by the {\bf Covering Theorem }
of Thurston \cite{Thurstonnotes} and Canary \cite{canary-cover}
(See Lemma \ref{hypqc}). Next, by a theorem of Canary and Minsky
\cite{can-min}, it follows that the convex hull boundary $D$ of
$M_1$ can be approximated by simplicial hyperbolic surfaces (see \cite{can-min} for details)
homotopic to $D$ with short tracks. Thus any simplicial
hyperbolic approximant $D_a$ would have to have bounded area and
hence bounded diameter modulo Margulis tubes (as in Lemma
\ref{scc-qc}). Thus so would $D$. Now, we repeat the argument in
the proof of Lemma \ref{drill-gr-qc}, to conclude that $D$ and
hence the convex core $CC(M_1)$ of $M_1$ lies in a uniformly
bounded neighborhood of $K_1$ in the graph metric. $\Box$

\smallskip

This approach
would circumvent the use of the Drilling Theorem at this
stage. However, since we shall again need it below, we retain our
approach here. 

\smallskip

Since $D$ bounds $X$, we would like to claim that the conclusion of
Lemma \ref{D-gr-qc}  follows
with $X$ in place of $D$. 
Though this  does not a priori follow in the hyperbolic
metric, it does follow for the graph metric.
This is because the double cover of $X$ is a 'drilled quasifuchsian'
manifold
(i.e. it is essentially $(D \times I)$ with some short curves
drilled). Further, any point in the convex core of
a quasifuchsian
$(D \times I)$ is close to a pleated
 surface by Theorem \ref{filling}.  Essentially the same argument
as in Lemma \ref{D-gr-qc} applies now. Details will be given below.

\begin{lemma}
There exists $C_1$ such that for any split component $K$, $K_{1d}$
is uniformly graph-quasiconvex in $M_{1d}$. \label{drill-gr-qc}
\end{lemma}

\begin{proof}
$X$ is double covered by $D \times I$ with cores of some Margulis
tubes drilled. Let $X_1$ denote this double cover. Note that $X_1$ is
convex, being a double cover of the convex compact $X$. 
By Theorem \ref{drill}, there exists $l > 0$ such that the drilled and undrilled manifolds are $2$-biLipschitz homeomorphic
away from Margulis tubes and cusps provided the Margulis tubes are $l$-thin.

Perform $(1,m)$ Dehn
filling on $X_1$ with sufficiently large $m=m(l)$ to ensure that the  resulting Margulis tube is $l$-thin. 
Let $X_{1f}$ be the resulting Dehn-filled manifold. By Theorem \ref{drill}, $X_{1f}$ is
uniformly quasiconvex in $M_{1f} = {\Hyp}^3/\Gamma$ where $\Gamma$ is a
quasiFuchsian surface group obtained by the above Dehn filling. 
(Theorem \ref{drill} gives a uniform
biLipschitz map outside Margulis tubes.)

Next, by Theorem \ref{filling}, for all $x \in X_{1f}$ there exists a
pleated surface $\Sigma \subset X_{1f}$ such that
$d(x,\Sigma ) \leq \delta$ where $\delta$ depends only on the genus of
$D$. 

Returning to $X_1$ via the Drilling Theorem \ref{drill} we see that
for all $x \in X_1$,

\begin{enumerate}
\item Either there exists a uniformly biLipschitz image of a
  hyperbolic surface $\Sigma_1 \subset X_1$ such that $d(x,\Sigma_1 )
  \leq \delta$. This is the case that the pleated surface $\Sigma$ misses all
  filled Margulis tubes.
\item Or, there exists a uniformly biLipschitz image of a
 subsurface $\Sigma_1$ of a
 hyperbolic surface  such that $d(x,\Sigma_1 )
  \leq \delta$  and such that the boundary of $\Sigma_1$ lies in a
 Margulis tube.  This is  the case that the pleated surface $\Sigma$ meets some
  filled Margulis tubes. Here,
 we can take $\Sigma_1$ to be the image of the
component of  {\em ($\Sigma$ minus Margulis tubes)} that lies near
 $x$.
\end{enumerate}

Again, passing down to $X$ under the double cover (from $X_1$ to $X$), 
we have, for all $x \in X$,

\begin{enumerate}
\item Either there exists a uniformly biLipschitz image of a
  hyperbolic surface $\Sigma_1 \subset X$ parallel to $D$ such that  $d(x,\Sigma_1 )
  \leq \delta$.
\item Or, there exists a uniformly biLipschitz image of a
 subsurface $\Sigma_1$ of a
 hyperbolic surface  such that $d(x, \Sigma_1 )
  \leq \delta$  and such that the boundary of $\Sigma_1$ lies on a
 Margulis tube. Further, $\Sigma_1$ is incompressible in the complement of $l$-thin Margulis tubes.
\end{enumerate}

In either case, the argument for Lemma \ref{D-gr-qc} now shows that for
 all $x \in X$ the distance $d_G(x,K_{1d})$ is uniformly bounded in
 the graph-metric $d_G$.
 Thus, we have shown that
$K_{1d}$ is uniformly graph-quasiconvex in $M_{1d}$. \end{proof}

 To
complete the proof
 of Proposition \ref{gr-qc} it is necessary to translate the
 content of Lemma \ref{drill-gr-qc} to the `undrilled' manifold
 $N_1$. We shall need to invoke
the Drilling Theorem \ref{drill} again.

\smallskip

\noindent {\bf Concluding the Proof of Proposition \ref{gr-qc}}:\\
While recovering data about $N_1$, it is  slightly
easier to handle the case where $D \cap K_{1d} = \emptyset$. Since we shall use the convex core boundary for both the drilled
as well as the undrilled manifolds in the rest of the proof, we change notation slightly and use\\
a) $D_d$ for the convex core boundary of the drilled manifold $N_d$.\\
b)  $D$ for the convex core boundary of the undrilled manifold $N$.\\

\smallskip

\noindent {\bf Case 1: $D_d \cap K_{1d} = \emptyset$ } \\
Filling $N_{1d}$ along the (drilled) $\T^1$, we get back $N_1$.
Since $D_d$ misses $K_{1d}$, the filled image of $X$ in $N_1$ is
$C_1$-quasiconvex for some $C_1$, depending on the biLipschitz
constant of Theorem \ref{drill} above. (One can see this easily
for instance from the fact that there is a uniform Lipschitz
retract of $N_{1d} - X$ onto $D_d$).

\smallskip

\noindent {\bf Case 2: $D_d \cap K_{1d} \neq \emptyset$ } \\
If $D_d$ meets some Margulis tubes $\T^1$, we enlarge $D$ to $D^{\prime}$ in $X_1$  by
letting $D^{\prime}$ be the boundary of $X_1 = X \cup \T^1$. The annular
intersections of $D$ with Margulis tubes are replaced by boundary
annuli contained in the boundary of $\T^1$.

It is easy enough to check that  $X_1$ is uniformly quasiconvex in the hyperbolic metric:   look at a universal cover $\til{X_1}$ of $X_1$ in
$\til{N_{1}}$. Then $\til{X_1}$ is a union of $\til{X}$ and
the lifts of $T$ that intersect it. All these lifts of $T$ are
disjoint. Hence $\til{X_1}$ is a `star' of convex sets all of
which intersect the convex set $\til X$. By (Gromov)
$\delta$-hyperbolicity, such a set is uniformly quasiconvex.

 Then as before, there
is a uniform Lipschitz retract of $N_{1d} - X_1$ onto
$D^{\prime}$. But now $D^{\prime}$ misses the interior of
$K_{1d}$ and we can apply the previous argument.

By Theorem \ref{drill} above, the diameter of the convex core boundary $D$ (or $D^{\prime}$ if $D$
intersects some Margulis tubes) in $N_1$ is bounded in terms of
the diameter of the convex core boundary $D_d$ in $N_{1d}$ and the {\bf uniform} biLipschitz
constant $L$ obtained from Theorem \ref{drill} above. Further, the
distance of $D$ from $K_1 \cup \T^1$ in $N_1$ is bounded in terms
of the distance of $D_d$ from $K_{1d} \cup \partial \T^1$ in
$N_{1d}$ and the  biLipschitz constant $L$.

Hence we can translate the content of Lemma \ref{drill-gr-qc} to
the `undrilled' manifold $N_1$. This
concludes the proof of Proposition \ref{gr-qc}: 
{\em Split components are uniformly graph-quasiconvex.} $\Box$

\begin{rmk} Our proof above
 uses the fact that the convex core $X$ of
$N_{1d}$ is a rather well-understood object, namely, a manifold
double covered by a  drilled convex hull of a quasi-Fuchsian
group. Hence,  it follows that the convex core $X$ is
uniformly congested, i.e. it has a uniform upper bound on its
injectivity radius. This is an approach to a conjecture of
McMullen \cite{bielefeld} (See also Fan \cite{fan1} \cite{fan2}). \\
A further point to be noted  is that we have  implicitly used here the idea of
drilling {\em disk-busting curves} introduced by Canary in
\cite{canary} and used again by Agol in his resolution of the
tameness conjecture \cite{agol-tameness}. \end{rmk}

\begin{rmk} Recall that {\it extended split components} were defined in $\til N$ by adjoining Margulis tubes abutting lifts of split components
to $\til N$.
The proof of Proposition \ref{gr-qc} establishes also the uniform graph-quasiconvexity of {\it extended split components} in $\til N$.
The metric obtained by electrocuting the family of convex hulls of extended split components in $\til N$ will be denoted as $d_{CH}$.
\label{extddsplit}
\end{rmk}

\subsection{Hyperbolicity in the graph metric} First a word about the modifications necessary for Simply Degenerate Groups.\\

\noindent {\bf Simply Degenerate Groups} We have so far mostly assumed,
for simplicity, that we  are dealing with totally
degenerate groups.  In a
simply degenerate $N$, the Minsky model is {\em uniformly}
biLipschitz to $N$ only in a neighborhood $E$ of the end. In this
case $(N \setminus E)$ is homeomorphic to $S \times I$. We declare
$(N \setminus E)$ to be the first block - a `thick block' in the
split geometry model. Thus the boundary blocks of Minsky are put
together to form one initial thick block. This changes the
biLipschitz constant, but the rest of the discussion, including
Proposition \ref{gr-qc} go through as before.

\bigskip

Construct
  a second auxiliary metric ${\til{N}}_2=(\til{N}, d_{CH})$ 
by electrocuting
  the elements $CH(\til{K})$ of convex hulls of {\it extended split components}. We show that the spaces ${\til{N}}_1 = (\til{N}, d_{G})$ and
${\tilde{N}}_2=(\til{N}, d_{CH})$ are quasi-isometric. In fact we show that the identity map from $\til N$ to itself induces this quasi-isometry after
the two different electrocutions.

\begin{lemma}
 The identity map from $\til N$ to itself
 induces a quasi-isometry of 
${\widetilde{N}}_1$ and ${\widetilde{N}}_2$. 
\label{qi12}
\end{lemma}

\begin{proof} We use $d_1$, $d_2$ as shorthand for the electric
metrics $d_G$ and $d_{CH}$ on ${\widetilde{N}}_1$ and ${\widetilde{N}}_2$. Since $\til{K}
\subset CH( \til{K} )$ for every split component, we have straightaway\\
\begin{center}
$d_1 (x,y) \leq d_2 (x,y)$ for all $x, y \in \til{M}$
\end{center}

To prove a reverse inequality with appropriate constants, it is enough
to show that each set $CH(\til{K} )$ (of diameter one in ${\til N}_2$) has uniformly
bounded diameter in ${\widetilde{N}}_1$. To see this, note that by definition of
graph-quasiconvexity, there exists $n$ such that for all $\til K$
and each point $a$ in $CH(\til{K})$, there exists a point $b \in \til{K}$
with $d_1(x,y) \leq n$. Hence by the triangle inequality, \\
\begin{center}
$d_2 (x,y) \leq 2n+1$ for all $x, y \in {CH(\til{K})}$
\end{center}
  Therefore, \\
\begin{center}
$d_2 (x,y) \leq (2n+1)(d_1 (x,y)+1)$ for all $x, y \in \til{N}$
\end{center}

This proves the Lemma. \end{proof}

\begin{cor} ${\til{N}}_1 = (\til{N}, d_{G})$ is Gromov-hyperbolic.
\label{dGhyp}
\end{cor}

\begin{proof} By Lemma \ref{farb1A}, ${\til{N}}_2=(\til{N}, d_{CH})$ is a $\delta$-hyperbolic metric space
for some $\delta \geq 0$. By quasi-isometry invariance of Gromov hyperbolicity,
so is ${\til{N}}_1 = (\til{N}, d_{G})$. \end{proof}

We have thus constructed a sequence of split surfaces
that satisfy the following two conditions in addition to Conditions (1)-(6)
of Remark \ref{wsplitrmk} for the Minsky model of a simply or totally degenerate surface group:\\

\begin{defn} \label{gqc} A model manifold of weak split geometry is said to be of {\bf split geometry} if \\
7) Each split component $\til{K}$ is 
quasiconvex (not necessarily uniformly) in the hyperbolic metric on $\til{N}$. \\
8) Equip $\til{N}$ with the {\it graph-metric} $d_G$ obtained by
electrocuting  (extended) split components $\til{K}$. Then the convex hull
$CH( \til{K})$ of any split component $\til{K}$ has uniformly bounded
diameter in the metric $d_G$. 
\end{defn}

Hence by Lemma \ref{hypqc} and Proposition \ref{gr-qc} we have the following.

\begin{theorem}
Any simply or doubly degenerate
surface group without accidental parabolics is biLipschitz homeomorphic to a model of split geometry. \label{minsky-split}
\end{theorem}

\section{Constructing Quasiconvex Ladders and Quasigeodesics} \label{ladder}

To avoid confusion we summarize the various metrics on $\til{M} , \til{N}$ and related
models that will be used: \\
1) The hyperbolic metric $d$ on $\til N$.\\
2) The weld-metric $d_{wel}$ obtained after welding the boundaries of Margulis tubes of $\til M$
to standard annuli (and before tube electrocution) where each horizontal circle
of a Margulis tube $T$ 
has a fixed non-zero length. This gives the welded model manifold $({M_{wel}}, d_{wel})$.\\
3) The tube-electrocuted metric $(M_{tel},d_{tel})$. We remind the reader that the underlying manifolds
${M_{wel}}, {M_{tel}}$ are the same.\\
4) The graph metric $d_G$. This is the notation for the 
electric metric on  $\EE(\til{M_{wel}}, \KK)$, where $\KK$ denotes the collection of
split components. We shall also use it for 
the 
electric metric on  $\EE (\til{N} , {\KK}^{\prime})$,
 where $\KK^{\prime}$ denotes the collection of extended
split components in $\til N$. The two electric metrics are quasi-isometric by  Remark \ref{alt}.\\

There will be two (families of) metrics on the universal cover $\til S$ of $S$:\\
1) The graph-electrocuted metric $d_{Gel}$ obtained by electrocuting the amalgamation components of $\til S$ that the lift of a weld-curve
cuts $\til S$ into.\\
2) The (Gromov) $\delta$-hyperbolic metric $d$ 
on $\til S$ obtained by lifting the metric on the {\it welded} surface. Recall that the metric $d$ on 
$\til S$ is the lift to the universal cover of a metric 
 on $S$  obtained by cutting out thin annuli and then welding the boundaries of the resulting extended split surface together.  The latter is uniformly biLipschitz to a fixed
hyperbolic structure on $S$. Hence we shall use $d$ to denote both the hyperbolic metric as well as those uniformly biLipschitz to it.\\

Note that the path metric induced on $\til{ S} \subset \til{B}$ by the graph metric $d_G$ 
on  $\EE(\til{M_{wel}}, \KK)$ is precisely $d_{Gel}$.

\subsection{Construction of Quasiconvex Sets for Building Blocks}

In this subsection, we describe the construction of a {\bf hyperbolic ladder} ${\LL}_\lambda$ 
restricted to  building blocks $B$. Putting these together we will show later that ${\LL}_\lambda$  is quasiconvex in
$(\til{M_{wel}}, d_G)$.

\smallskip

\noindent {\bf Construction of ${\LL}_\lambda (B)$ - Thick Block}\\
Let $B$ be a thick block. By definition $B$ is a uniformly biLipschitz homeomorphic
image of $S \times I$. Let $F_B: S \times I \rightarrow B$ denote the biLipschitz homeomorphism.

Let $\lambda = [a,b]$ be a geodesic segment in
$\widetilde{S}$. Let $\lambda_{Bi}$ denote $F_B(\lambda \times \{ i
\})$ for $i=0,1$. 

Equivalently, let  $\phi : F_B(\til{S} \times \{0\}) \rightarrow F_B(\til{S} \times \{1 \})$ be 
given by $\phi (F_B(x,0)) = F_B(x,1)$. 
 The induced map on geodesics will be denote by  $\Phi$, which can be described as follows. Let $\lambda$
be a geodesic joining $a, b \in F_B(\til{S} \times \{0\})$ and let $\Phi ( \lambda )$ denote the
a geodesic 
joining $\phi (a), \phi (b)$. Let $\lambda_{B1}$ denote $\Phi (\lambda
) \times \{ 1 \}$.

For the universal cover $\widetilde{B}$ of the thick block $B$, define

\begin{center}

${\LL}_\lambda (B) = \bigcup_{i=0,1} \lambda_{Bi}$

\end{center}

\begin{defn} Each $\widetilde{S} \times i$ for $i = 0, 1$
 will be called a {\bf horizontal sheet} of
$\widetilde{B}$ when $B$ is a thick block.\end{defn}

\noindent {\bf Construction of ${\LL}_\lambda (B)$ - Split Block}\\
As above, let $\lambda = [a,b]$ be a geodesic segment in
$\widetilde{S}$, where $S$ is regarded as the base surface of $B$ in the tube electrocuted model. Let $\lambda_{B0}$ denote $\lambda \times \{ 0
\}$. Then for each split component $K$, $K\cap (S \times{i})$ ($i=0,1$) is an {\it amalgamation component} of $\til S$.
Also, $S \times{i}$, ($i=0,1$), are the boundary {\it welded split surfaces} forming the horizontal boundary of $B$, uniformly biLipschitz to $S$ with
a fixed hyperbolic metric.   Note further that the induced 
path metric $d_{Gel}$ on $\til{ S} \times{i}$ ($i=0,1$) is the electric pseudo-metric on $\widetilde{S}$ obtained by
electrocuting amalgamation components of $\til S$.

Let $\lambda_{Gel}$ denote the electro-ambient quasigeodesic (Lemma \ref{Gea}) joining $a,
b$ in $(\til{S}, d_{Gel})$.  Let $\lambda_{B0}$ denote
$\lambda_{Gel} \times \{ 0 \}$. 

Then the map $\phi : S \times \{ 0 \} \rightarrow S \times \{ 1 \}$ taking $(x,0)$ to $(x,1)$ is 
a component
preserving diffeomorphism. Let $\tilde{\phi }$ be the lift of  $\phi$ to
$\widetilde{S}$ equipped with  the electric metric $d_{Gel}$. Then $\til \phi$ is an isometry by Lemma \ref{phi-isom2}. Let
${\Phi }$ denote
the induced map on electro-ambient quasigeodesics, i.e. if $\mu = [x,y] \subset (
\widetilde{S} , d_{Gel} )$, then $\Phi ( \mu ) = [ \phi (x), \phi (y)
]$ is the electro-ambient quasigeodesic joining $\phi (x), \phi (y)$. Let $\lambda_{B1}$
denote $\Phi ( \lambda_{Gel} ) \times \{ 1 \}$.

For the universal cover $\widetilde{B}$ of the split block $B$, define:

\begin{center}

${\LL}_\lambda (B) = \bigcup_{i=0, 1} \lambda_{Bi}$

\end{center}

\begin{defn}  Each $\widetilde{S} \times i$ for $i = 0
, 1$ will be called a {\bf horizontal sheet} of
$\widetilde{B}$ when $B$ is a split block.\end{defn}

\noindent {\bf Construction of $\Pi_{\lambda ,B}$ - Thick Block}\\
For $i = 0, 1$, let  $\Pi_{Bi}$ denote nearest point
projection of
$\widetilde{S} \times \{ i \}$
onto $\lambda_{Bi}$ in the path metric on $\widetilde{S} \times \{ i \}$.

For the universal cover $\widetilde{B}$ of the thick block $B$, define:

\begin{center}

$\Pi_{\lambda ,B}(x) = \Pi_{Bi}(x) , x \in \widetilde{S} \times \{ i
  \} , i=0,1$

\end{center}

\smallskip

\noindent {\bf Construction of $\Pi_{\lambda ,B}$ - Split Block}

For $i = 0, 1$, let  $\Pi_{Bi}$ denote nearest point
projection of
$\widetilde{S} \times \{ i \}$
onto $\lambda_{Bi}$.

Here the
nearest point projection is taken in the sense of the definition
preceding
Lemma \ref{hyp=elproj}, i.e.  minimizing the ordered pair $(d_{Gel},
d)$ in the lexicographic order on $\mathbb{R} \times \mathbb{R}$
(where $d_{Gel}, d$ refer to electric and (biLipschitz)-hyperbolic
metrics respectively.)

For the universal cover $\widetilde{B}$ of the split block $B$, define:

\begin{center}

$\Pi_{\lambda ,B}(x) = \Pi_{Bi}(x) , x \in \widetilde{S} \times \{ i
  \} , i=0,1$

\end{center}

{\bf $\Pi_{\lambda , B}$ is a coarse Lipschitz retract - Thick Block}

The proof for a thick block is exactly as in \cite{mitra-trees} and
\cite{mahan-bddgeo}. We omit it here.

\begin{lemma} (Theorem 3.1 of \cite{mahan-bddgeo})
There exists $C > 0$ such that the following holds: \\
Let $x, y \in \widetilde{S} \times \{ 0, 1\} \subset \widetilde{B}$
for some thick block $B$. 
Then $d( \Pi_{\lambda , B} (x), \Pi_{\lambda , B} (y)) \leq C d(x,y)$.
\label{retract-thick}
\end{lemma}

\noindent {\bf $\Pi_{\lambda , B}$ is a retract - Split Block}
\begin{lemma}
There exists $C > 0$ such that the following holds: \\
Let $x, y \in \widetilde{S} \times \{ 0,1\} \subset \widetilde{B}$
for some split block $B$. 
Then $d_{G}( \Pi_{\lambda , B} (x), \Pi_{\lambda , B} (y)) \leq C d_{G}(x,y)$.
\label{retract-thin}
\end{lemma}

\begin{proof} It is enough to show this for the  following cases: \\
Case 1) $x, y \in \widetilde{S} \times \{ 0 \} $ OR 

 $x, y \in \widetilde{S} \times \{ 1 \} $. \\ This follows directly from Lemma \ref{easyprojnlemma}.\\
Case 2) $x = (p,0)$ and $y = (p,1)$ for some $p \in \til S$. \\ First note that $(\til{ S} , d_{Gel})$ is uniformly $\delta$-hyperbolic as a metric space (in fact 
uniformly quasi-isometric to a tree) and
$\til{\phi} : \til{S} \times \{ 0 \} \rightarrow \til{S} \times \{ 1 \}$ induces an isometry of the $d_{Gel}$ metric by Lemma \ref{phi-isom2} as $\phi$ is a component
preserving diffeomorphism. Case 2 now follows from the fact that that quasi-isometries and nearest-point projections almost commute (Lemma \ref{almost-commute} ).\end{proof}

In the next section, we shall come across the situation where one horizontal surface $\til S \times \{ i \}$ can occur
as the bottom surface of a split  block $B_2$ and as the top surface of a thick block $B_1$, or vice versa. 
Alternately it could occur as the bottom surface of a split  block 
and as the top surface of a different split block where the collection of splitting tubes differ.
In either situation we shall denote the bottom block by $B_1$ and the top block by $B_2$.
In this case, the nearest point projection
could be in any  of the following senses: \\
a) Projection onto a (biLipschitz)-hyperbolic geodesic $[a,b]$ in the (biLipschitz)-hyperbolic metric $d$ on $\til S$. \\
b) Projection onto an electro-ambient quasigeodesic $[a,b]_{ea}$ minimizing the ordered pair $(d_{Gel1},
d)$, where $d_{Gel1}$ denotes the electric metric on $S$ induced by the split block $B_1$. \\
c) Projection onto an electro-ambient quasigeodesic $[a,b]_{ea}$ minimizing the ordered pair $(d_{Gel2},
d)$, where $d_{Gel2}$ denotes the electric metric on $S$ induced by the split block $B_2$. 

\begin{lemma} {\bf $\Pi_{\lambda , B}$ is coarsely well-defined:} There exists $C_0 > 0$ such that the following
holds. \\
Suppose that $\Pi^1_{\lambda , B}$
and $\Pi^2_{\lambda , B}$ are projections defined in any two of the above senses.
Then $$d( \Pi^1_{\lambda , B}(p), \Pi^2_{\lambda , B}(p)) \leq C_0$$ for all $p \in \til{S}$.
\label{disamb0}
\end{lemma}

\begin{proof} By Lemma \ref{hyp=elproj},
 hyperbolic and electric projections of $p$ onto the
(Gromov) $\delta$-hyperbolic geodesic $[a,b]$ and the electro-ambient geodesic
$[a,b]_{ea}$ respectively `almost agree':  If $\pi_h$ and $\pi_e$
denote the hyperbolic and electric projections, then
there exists (uniform) $C_1 > 0$ such that
$d ( \pi_h(p), \pi_e(p)) \leq C_1$.
The Lemma follows if one of the blocks are thick.

If both blocks are split blocks, then $d ( \pi_h(p), \Pi^i_{\lambda , B}(p)) \leq C_1$, for $i=1,2$
by the above argument. Taking $C_0 = 2 C_1$, we are through. \end{proof}

\subsection{Construction of ${\LL}_\lambda$ and $\Pi_\lambda$}

A subset $Z \subset (X,d)$ shall be called a coarse $k$-net in $X$ if $X=N_k(Z,d)$. A subset 
$Z \subset (X,d)$ shall be called a coarse net if it is a coarse $k$-net in $X$ for some $k$.

Given a manifold $M$ of split geometry, we know that $M$ is
homeomorphic to $S \times J$ for $J = [0, \infty ) $ or $( - { \infty
  }, {\infty })$. By definition of split geometry, 
there exists a sequence of  blocks $B_i$ (thick or split) such that $M_{wel}=\cup_i B_i$.
 Denote: \\
$\bullet$ $\LL_{\mu , B_i} = \LL_{i \mu }$ \\
$\bullet$ $\Pi_{\mu , B_i} = \Pi_{i \mu }$ \\

Now for a block $B = S \times I$ (thick or amalgamated),  a natural
map $\Phi_B$ may be defined taking
 $ \mu = \LL_{\mu}(B) \cap F_B(\widetilde{S} \times \{ 0 \}) $ to a
geodesic $\LL_{\mu}(B) \cap F_B(\widetilde{S} \times \{ 1 \}) = \Phi_B ( \mu
)$. Similarly $\Phi_B^{-1}$ may be defined taking

 $ \mu = \LL_{\mu}(B) \cap F_B(\widetilde{S} \times \{ 1\}) $ to 
$\LL_{\mu}(B) \cap F_B(\widetilde{S} \times \{ 0 \}) = \Phi_B^{-1} ( \mu
)$. 

Let the map
$\Phi_{B_i}$   (resp.  $\Phi_{B_i}^{-1}$ ) be denoted as $\Phi_i$  (resp.  $\Phi_{i}^{-1}$ ).

We start with a reference block $B_0$ and a reference geodesic segment
$\lambda = \lambda_0$ on the `lower surface' $\widetilde{S} \times \{
0 \}$.
Now inductively define: \\
$\bullet$ $\lambda_{i+1}$ = $\Phi_i ( \lambda_i )$ for $i \geq 0$\\
$\bullet$ $\lambda_{i-1}$ = $\Phi_i^{-1} ( \lambda_i )$ for $i \leq 0$ \\

Finally define

$${\LL}_\lambda = \bigcup_i \lambda_i.$$

${\LL}_\lambda$ is the {\it hyperbolic ladder} promised.

Recall that 
each $\widetilde{S} \times i$ for $i = 0, 1$ is called a {\bf horizontal sheet} of
$\widetilde{B}$. We will restrict our attention to the union of the horizontal
sheets $\widetilde{M_H}\subset \widetilde{M_{wel}}$ with the 
metric induced from the graph model. Since
 $\widetilde{M_H}$ is a {\it coarse $1$-net} in $(\widetilde{M_{wel}}, d_G)$, we will be
able to get all the coarse information we need by restricting
ourselves to  $\widetilde{M_H}$.

Clearly, 
${\LL}_\lambda  \subset \widetilde{M_H} \subset \widetilde{M_{wel}}$.

Let the bottom horizontal sheet of $\til{B_i}$ be denoted as $\til{S_i}$.
$\Pi_{i \lambda }$ is defined to be the nearest point projection of $\til{S_i}$ onto $\lambda_i$.

\begin{rmk} \label{disamb} {\rm
As noted earlier, the nearest point projection $\Pi_{i \lambda }$
could be in any  of the following senses: \\
a) Projection onto a (biLipschitz)-hyperbolic geodesic $[a,b]$ in the (biLipschitz)-hyperbolic metric $d$ on 
$\til S$. \\
b) Projection onto an electro-ambient quasigeodesic $[a,b]_{ea}$ minimizing the ordered pair $(d_{Gel1},
d)$, where $d_{Gel1}$ denotes the electric metric on $S$ induced by the split block $B_1$ whose top boundary is $S$. \\
c) Projection onto an electro-ambient quasigeodesic $[a,b]_{ea}$ minimizing the ordered pair $(d_{Gel2},
d)$, where $d_{Gel2}$ denotes the electric metric on $S$ induced by the split block $B_2$   whose bottom boundary is $S$. \\
By Lemma \ref{disamb0},} {\bf $\Pi_{i \lambda }$ is coarsely well-defined, i.e. any two choices are a uniformly bounded $d-$distance apart.} \end{rmk}

Hence we define the projection
$$\Pi_\lambda = \bigcup_i \Pi_{i \lambda }.$$

$\Pi_\lambda$ is defined from  $\widetilde{M_H}$ to ${\LL}_\lambda$.

\begin{theorem}
There exists $C > 0$ such that 
for any geodesic $\lambda = \lambda_0 \subset \widetilde{S} \times \{
0 \} \subset \widetilde{B_0}$, the retraction $\Pi_\lambda :
\widetilde{M_H} \rightarrow {\LL}_\lambda $ satisfies
$$d_{G}( \Pi_{\lambda } (x), \Pi_{\lambda } (y)) \leq C
d_{G}(x,y) + 
C.$$
\label{retract}
\end{theorem}

\begin{proof} This is now
a direct consequence of Lemmas \ref{retract-thick} and \ref{retract-thin} and  Remark \ref{disamb}. 
 \end{proof}

 For Theorem \ref{retract} above, note that all that we
really
require is that the universal cover $\widetilde{S}$ is a  Gromov-hyperbolic
metric space. There is {\em no restriction on $\widetilde{M_H}$.} In fact,
Theorem \ref{retract} would hold for general stacks of (Gromov) hyperbolic
metric spaces with  blocks of split geometry. However, in the present situation we have more

\begin{cor} \label{ladder-qc}
${\LL}_\lambda $ is quasiconvex in $(\til{M_{wel}}, d_G)$.
\end{cor}

\begin{proof} By Corollary \ref{dGhyp} $(\til{M_{wel}}, d_G)$ is (Gromov)-hyperbolic.
Hence ${\LL}_\lambda $  is a coarse Lipschitz retract in a (Gromov)-hyperbolic space
by Theorem \ref{retract}. Therefore
${\LL}_\lambda $ is quasiconvex in $(\til{M_{wel}}, d_G)$. \end{proof}

\subsection{Heights of Blocks}

Recall that each thick or split block $ B_i$ is identified with $S \times I$
where each fiber $\{ x \} \times I$ has length $\leq l_i$ for some
$l_i$, called the {\it thickness} of the block $B_i$.

\noindent {\bf Observation:} $\widetilde{M_{H}}$ is a `coarse net' in $(\widetilde{M_{wel}}, d_G)$
in the  graph metric, but not in the weld metric $d_{wel}$, the tube-electrocuted metric $d_{tel}$, nor the model metric $d_M$ (cf. Remark \ref{alt} for  $d_M$).
In the graph model, any point can be connected by a vertical segment
of length $\leq 1$ to one of the boundary horizontal sheets.

However,  there are points within
split components which are at a  $d_{wel}$-distance of the order of $l_i$ from the boundary horizontal sheets. Since $l_i$ could be arbitrary,
$\widetilde{M_{H}}$ is no longer necessarily a `coarse net' in $(\widetilde{M} , d_{wel})$ or $(\widetilde{M} , d_{tel})$.

\begin{lemma}
There exists a function $g: \mathbb{Z} \rightarrow \mathbb{N}$ such that
for any block $B_i$ (resp. $B_{i-1}$), and $x \in \lambda_i$, there exists $x^{\prime} 
\in \lambda_{i+1}$ (resp. $\lambda_{i-1}$) for $i \geq 0$ (resp. $i \leq
0$), satisfying:

\begin{center}

$d_{wel}(x, x^{\prime}) \leq g(i)$,  \\  $d_{M}(x, x^{\prime}) \leq g(i)$

\end{center}
\label{bddheight}
\end{lemma}

\begin{proof}
Let $\mu \subset \widetilde{S} \times \{ 0 \} \subset \widetilde{B_i}$ be a
geodesic in a (thick or split) block. 
Then  from the product structure on the block $B_i$, there exists a $(K_i, \epsilon_i )$- quasi-isometry $\psi_i$ 
from $\widetilde{S} \times \{ 0 \}$ to $\widetilde{S} \times \{ 1 \}$
and $\Psi_i$ is the induced map on geodesics. Hence, for any $x \in
\mu$, $\psi_i (x)$ lies within some bounded distance $C_i$ of $\Psi_i
( \mu )$. But $x$ is connected to $\psi_i (x)$ by \\
\noindent {\it Case 1 - Thick Blocks:} a vertical segment of uniformly bounded length ($\leq C$ say).  \\

\noindent{\it Case 2 - Split Blocks:}

\smallskip

Thus $x$ can be connected to a point $x^{\prime} \in \Psi_i ( \mu
)$ by a path of length less than $g(i) =   l_i +C_i + C$.
 Recall that $\lambda_i$ is the geodesic on the lower
horizontal surface of the block $\widetilde{B_i}$. The same can be
done for blocks $\widetilde{B_{i-1}}$ and {\em going down} from
$\lambda_i$ to $\lambda_{i-1}$.

By Remark \ref{alt}, the same argument works for the model manifold $(\til{M}, d_M)$.
\end{proof}

\section{Recovery} The previous Section was devoted to constructing a quasiconvex ladder in the graph metric which is an electric metric. In this section we shall
be concerned with recovering information about hyperbolic geodesics from electric ones.
Since a host of  metrics will make their
appearance in this section, we shall refer to (quasi)geodesics in $(\til{M_{wel}},d_G)$, $(\til{M_{wel}},d_{wel})$, 
$(\til{M_{tel}}, d_{tel})$  and
$(\til{M},d_{CH})$ as $d_G$-(quasi)geodesics, $d_{wel}$-(quasi)geodesics, $d_{tel}$-(quasi)geodesics and $d_{CH}$-(quasi)geodesics respectively.
Recall that the union of the horizontal sheets $\til{S_i} \subset \til{M_{wel}}$ is denoted as $\til{M_H}$ and that
the projection $\Pi_\lambda$ occurring in Theorem \ref{retract} is defined only on $\til{M_H}$ and {\it not}
all of $\til{M_{wel}}$.

\subsection{Scheme of Recovery}\label{scheme}
The recovery is in several stages. We sketch the scheme of recovery in some detail in this subsection for the convenience of the reader. 
A first problem in recovering data about hyperbolic geodesics from $d_G$-geodesics
is the absence of canonical representatives in $(\til{M_{wel}},d_{wel})$ of $d_G$-geodesics.
In Section \ref{adm-sxn}, we address this problem by making a choice of paths in  $(\til{M_{wel}},d_{wel})$ 
representing $d_G$-geodesics. We call these {\it admissible paths}. Roughly speaking, admissible paths
are built up of \\
a) vertical segments of the form $\{ x \} \times [0,1]\subset \til{B} = \til{S}\times [0,1] $,
where $B$ is a block (thick or split) and $x \in \til{S}$.\\
b) horizontal segments consisting of geodesics in the horizontal sheets of $\til{M_H}$.\\
Let $\lambda \subset \til{S} (\subset \til{M_{wel}})$ be a geodesic in the intrinsic metric on $\til{S}$, where
$S$ is identified with the base surface $S\times \{ 0 \}$ of the first block in $\til{M_H}$.
Let $\beta_e$ denote an admissible path representing a $d_G$-geodesic joining the endpoints of   $\lambda$ in $(\til{M_{wel}},d_G)$.

We would like to project $\beta_e$ using ${\Pi}_\lambda$  onto the ladder ${\mathcal{L}}_\lambda$ to obtain a quasigeodesic contained
in ${\mathcal{L}}_\lambda$. Unfortunately, ${\Pi}_\lambda$ is defined only on $\til{M_H}$ and there is no natural
way to extend it to all of $\til{M_{wel}}$. To circumvent this problem we first define in Section \ref{adm-sxn}
a subcollection of the family of
admissible paths, called ${\mathcal{L}}_\lambda$-admissible paths.  Roughly speaking, ${\mathcal{L}}_\lambda$-admissible paths
are those admissible paths whose horizontal segments lie on or near ${\mathcal{L}}_\lambda$. 

Then in Section \ref{join-sxn} we project $\beta_e\cap \til{M_H}$ using ${\Pi}_\lambda$  onto the ladder ${\mathcal{L}}_\lambda$.
Since $\beta_e$ is itself an admissible path, there is a sequence of points $a_1, b_1, a_2, b_2, \cdots , a_k, b_k$ such that
the piece of $\beta_e$ joining $a_i$ to $b_i$ is horizontal, whereas the piece of $\beta_e$ joining $b_i$ to $a_{i+1}$ is vertical.
In particular, $b_i$ and $a_{i+1}$ must lie in the same split component if they lie in (the universal cover of) a split block. 
In this case ${\Pi}_\lambda (b_i)$ and ${\Pi}_\lambda (a_{i+1})$ must also lie in the same split component. This allows us to join 
the sequence of points ${\Pi}_\lambda (a_1), {\Pi}_\lambda (b_1), {\Pi}_\lambda (a_2), {\Pi}_\lambda (b_2), \cdots , 
{\Pi}_\lambda (a_k), {\Pi}_\lambda (b_k)$ by alternating horizontal and vertical segments to obtain an ${\mathcal{L}}_\lambda$-admissible path
$\beta_{adm}$ representing a (uniform) $d_G$-quasigeodesic joining the endpoints of   $\lambda$ in $(\til{M_{wel}},d_G)$.
Lemma \ref{far-nopunct} now establishes that if  $\lambda$ lies outside a large ball about a reference point in $\til S$, then 
$\beta_{adm}$ also lies outside a large ball about a reference point in $(\til{M_{wel}},d_{wel})$.

In Section \ref{recovery1} we construct an electro-ambient quasigeodesic  $\beta_{ea}$ in $(\til{M_{wel}},d_G)$ from the
${\mathcal{L}}_\lambda$-admissible path
$\beta_{adm}$,  constructed in Section \ref{join-sxn}. The idea is simple. Denote by $\til{K_{ij}}
\subset \til{B_i}$ the  split components in the universal cover of a split block $B_i$. Replace the intersection
$\beta_{adm} \cap \til{K_{ij}}$
(of $\beta_{adm}$ with any such split component $\til{K_{ij}}$) by a geodesic in $\til{K_{ij}}$ joining the
 end-points of $\beta_{adm} \cap \til{K_{ij}}$.
 Then
$\beta_{ea}$ continues to satisfy the conclusions of Lemma \ref{far-nopunct}, i.e. if  $\lambda$ lies outside a large ball
 in $\til S$, then 
$\beta_{ea}$  lies outside a large ball  in $(\til{M_{wel}},d_{wel})$.
What is crucial
at this stage of the recovery is the quasiconvexity of $\til{S_i}$ and $\til{S_{i+1}}$ in $\til{B_i}$, where
the quasiconvexity constant depends only
on $i$.

Finally in Section \ref{recovery2}, we construct an electro-ambient quasigeodesic  $\beta_{ea2}$ in $(\til{M},d_{CH})$
from $\beta_{ea}$. To do this, we first replace $\beta_{ea}$ by a path $\beta_{ea1}$ in $\til{M}$ such that 
$\beta_{ea1}$ coincides with $\beta_{ea}$
outside Margulis tubes and consists of hyperbolic geodesic segments within Margulis tubes. Then as above, we  
replace the intersection
$\beta_{ea1} \cap CH(\til{K_{ij}})$
(of $\beta_{ea1}$ with the convex hull $CH(\til{K_{ij}})$
of a split component $\til{K_{ij}}$) by a geodesic in $CH(\til{K_{ij}})$ joining the
 end-points of $\beta_{ea1} \cap CH(\til{K_{ij}})$. This gives us the required electro-ambient  quasigeodesic  $\beta_{ea2}$ in $(\til{M},d_{CH})$.
Again, $\beta_{ea2}$ continues to satisfy the conclusions of Lemma \ref{far-nopunct}, i.e. if  $\lambda$ lies outside a large ball
 in $\til S$, then 
$\beta_{ea2}$  lies outside a large ball  in $\til{M}$ (where the latter is equipped with the model metric).
The last statement follows from the uniform graph  quasiconvexity of split components (Proposition \ref{gr-qc}).

It is a small step from here to the main Theorem \ref{crucial} in Section \ref{ct-nopunct}; so we mention it here. 
Lemma \ref{ea-strong} ensures that the geodesic $\beta^h$ in $\til{M}$ joining the end-points of $\beta_{ea2}$ 
lies in a uniformly bounded neighborhood
of $\beta_{ea2}$ (see the figure just after Lemma \ref{ea-strong}). Note that it is at this stage that 
 we use explicitly the {\it weak relative hyperbolicity} of $\til M$ relative to the collection
of convex hulls of split components. Though  $\beta_{ea2}$ could be very far from a hyperbolic geodesic,
 Lemma \ref{ea-strong} forces  $\beta^h$ to lie in a bounded neighborhood of it.
 Hence if $\lambda$ lies outside a large ball
 in $\til S$, then $\beta^h$  lies outside a large ball  in $\til{M}$. Lemma
 \ref{contlemma} now furnishes the Cannon-Thurston map we want.

\subsection{Admissible Paths}\label{adm-sxn}
We want to first define a collection of  paths
 lying in a bounded
neighborhood of  ${\LL}_\lambda$ in $(\til{M_{wel}},d_G)$. Since ${\LL}_\lambda$ is not connected, it
 does not make  sense to speak of the path-metric on
 ${\LL}_\lambda$. To remedy this we shall introduce in this subsection the class of {\it ${\LL}_\lambda$-elementary admissible paths}
whose horizontal pieces are contained  in a neighborhood of  ${\LL}_\lambda$ in $\til{M_H}$. Further the distance of 
${\LL}_\lambda$-elementary admissible paths from ${\LL}_\lambda$
 will be controlled. An {\bf  ${\LL}_\lambda$-admissible path} will be a composition
 of  ${\LL}_\lambda$-elementary admissible paths.

We first define admissible paths in general. Let $B$ be a  thick or split block  in $M_{wel}$. We shall
identify $B$ with a product $S \times I$ as usual. In particular for $B$ a split block, and  any $x \in \widetilde{S}$, a
vertical segment of the form $x \times [0,1]$ will be assumed to be contained in some split component $\til{K} \subset \til{B}$. 

\begin{defn}{\rm
An {\bf admissible path in $\til{B}$} $ (\subset \til{M_{wel}})$
 is a path that can be decomposed into subpaths  of the following two types:\\
1) Horizontal segments along some $\widetilde{S} \times \{ i \}$
  for  $i = \{ 0, 1 \}$.\\
2)  Vertical segments of the form $x \times [0,1]$   where $x \in \widetilde{S}$.\\

An {\bf admissible path} $\sigma$ in $\til{M_{wel}}$ is a path such that for every (thick or split)
block $B$, any connected component of 
$\sigma \cap \til{B}$ is an admissible path in $\til{B}$.

An {\bf admissible $K-$ quasigeodesic} is an
admissible path that is a $K-$quasigeodesic in
$(\widetilde{M_{wel}}, d_{G})$. }
\end{defn}

\begin{lemma} \label{admqg} Given $K \geq 1$ there exists $K_1 \geq 1$ such that the following holds:\\
Let  $\beta_e$ be a ($d_G-$) $K-$quasigeodesic  in
$(\widetilde{M_{wel}}, d_{G})$. Then 
there exists an  admissible $K_1-$quasigeodesic $\beta_e^{\prime}$ joining the end-points of $\beta_e$.
\end{lemma}

\begin{proof}
Without
loss of generality, we can assume that $\beta_e$ does not
 back-track relative 
to the collection of split components, as any back-tracking can be removed without increasing the
$d_G-$length of $\beta_e$ (see \cite{farb-relhyp} for instance).
 We shall now convert
$\beta_e$ into an  admissible electric quasigeodesic without backtracking
joining the same pair of points as $\beta_e$. To do this we shall look at 
connected components of 
$\beta_e \cap \til{B}$ for any block $B$ and replace them with admissible paths. We identify $B$ with $S \times [0,1]$ and we call 
$\til{S} \times \{ 0 \}$ and $\til{S} \times \{ 1 \}$  the lower and upper boundary components of $\til B$.
Also let $P_0, P_1$ denote the natural projections from $\til{S} \times [0,1]$ to $\til{S} \times \{ 0 \}$ and $\til{S} \times \{ 1 \}$
respectively given by $P_0 (x,t) = (x, 0)$ and $P_1 (x,t) = (x, 1)$.

Now let $B$ be a  block (thick or split) and let $\beta_e \cap \til{B} \neq \emptyset$. Let $\beta_1$
be a connected component of $\beta_e \cap \til{B}$. Let $b_1, b_2$ be the end-points of $\beta_1$.
Two cases arise.

If both $b_1, b_2$ belong to the same boundary component, then we replace $\beta_1$ by $\beta_1^{\prime}=
P_i(\beta_1)$, where $i = 0$ or $1$ according as
$b_1, b_2 \in \til{S} \times \{ 0 \}$ or $\til{S} \times \{ 1 \}$. 

If  $b_1, b_2$ belong to different boundary components, then assume without loss of generality that 
$b_1 \in \til{S} \times \{ 0 \}$ and let $b_2 = (z,1) \in \til{S} \times \{ 1 \}$. Then replace $\beta_1$ by 
$\beta_1^{\prime}=P_0(\beta_1) \cup \{ z \} \times [0,1]$.

Performing this replacement for every block $B$ and every connected component of $\beta_e \cap \til{B}$ we obtain the 
required admissible quasigeodesic  $\beta_e^{\prime}$  joining the end-points of $\beta_e$.

It remains to show that if $\beta_e$ is a ($d_G-$) $K-$quasigeodesic, then
 $\beta_e^{\prime}$ is indeed an admissible $K_1-$quasigeodesic, where $K_1$ depends only on $K$. 

 For $B$ a split block, the $d_G$ length of any $\beta_1 \subset \til{B}$ is the same as  the $d_G$ length of the corresponding $\beta_1^{\prime}$
constructed to replace it as above. This is because the $d_G$ length of  $\beta_1$ is equal to 
the number of
split blocks that $\beta_1$ cuts. 

For $B$ a thick block, the inclusion of $\til{S} \times \{ 0 \}$ (or $\til{S} \times \{ 1 \}$) into $\til B$ is a uniform quasi-isometry as
the thickness of thick blocks is uniformly bounded. Hence $\beta_1^{\prime}$ is a $K_1-$quasigeodesic where $K_1$ depends only on $K$.
The Lemma follows.
\end{proof}

We shall now choose a subclass of these admissible paths to define
${\LL}_\lambda$-elementary admissible paths. The constants $C, C(B), K(B)$ etc. below will be independent of the geodesic $\lambda$, the initial
geodesic in the ladder ${\LL}_\lambda$.

\smallskip

\noindent {\bf  ${\LL}_\lambda$-elementary admissible paths in the thick block}\\
Let $B = S \times [i,i+1]$ be a thick block, where each $(x,i)$ is
connected by a vertical segment  to $( x,i+1)$. Let
$\phi$ be the map that takes $(x,i)$ to $(x,i+1)$.
 Also let $\Phi$ be the map on geodesics induced by
$\phi$. Let ${\LL}_\lambda \cap \widetilde{B} = \lambda_i \cup
\lambda_{i+1}$ where $\lambda_i$ lies on $\widetilde{S} \times \{ i
\}$ and $\lambda_{i+1}$ lies on  $\widetilde{S} \times \{ i+1
\}$. Let $\pi_j$, for $j = i, i+1$  denote nearest-point projections of
$\widetilde{S} \times \{ j \}$ onto $\lambda_j$. Since $\phi$
is a quasi-isometry, there exists $C > 0$ such that\\
a)  for all $(x,i) \in
\lambda_i$, $(x,i+1)$ lies in a $C$-neighborhood of $\Phi
(\lambda_i ) = \lambda_{i+1}$. \\
b) for all $z \in \til{S}$, $d_{wel}(\pi_i(z,i), \pi_{i+1}(z,i+1))\leq C$ (by Lemma \ref{retract-thick} or Theorem \ref{retract}).\\
We emphasize here that $C$ is independent of both the thick block $B$ and the geodesic $\lambda$ (and hence the ladder ${\LL}_\lambda$).
It depends only on the model manifold $M$.

The same conclusions hold for $\phi^{-1}$ and
points in $\lambda_{i+1}$, where $\phi^{-1}$ denotes the {\em
  quasi-isometric inverse} of $\phi$ from
$\widetilde{S} \times \{ i + 1 \}$ to $\widetilde{S} \times \{ i \}
$. The {\bf
  ${\LL}_\lambda$-elementary admissible paths} in $\widetilde{B}$ are defined to be paths consisting
of the following:\\
1) Horizontal geodesic subsegments of $\lambda_j$,  $j = \{ i, i+ 1 \}$.\\
2) Vertical segments of $d_G$ length $1$ joining $x  \times \{ 0 \}$ to
$x  \times \{ 1 \}$.  Note that for thick blocks,  $d_G=d_{wel}$.\\
3) Horizontal geodesic segments lying in a $C$-neighborhood of
  $\lambda_j$, $j = i, i+1$.

\smallskip

\noindent {\bf  ${\LL}_\lambda$-elementary admissible paths in the split block}\\
Let $B = S \times [i,i+1]$ be a split block, where each $(x,i)$ is
connected by a  segment of $d_G$ length one and $d_{wel}$-length
 $\leq C(B)$ (due to bounded thickness of $B$, Lemma \ref{bddheight}) to $( 
x,i+1)$. As before we regard $\phi$ as the  map from
$\widetilde{S} \times \{ i \}$ to $\widetilde{S} \times \{ i+1 \}$
that is the identity on the first component.
 Also let $\Phi$ be the map  on electro-ambient quasigeodesics induced by
$\phi$. Let ${\LL}_\lambda \cap \widetilde{B} = \bigcup_{j=i,  i+1}
 \lambda_j$ where $\lambda_j$ lies on $\widetilde{S} \times \{ j \}$.  $\pi_j$  denotes nearest-point projection of
$\widetilde{S} \times \{ j \}$ onto $\lambda_j$ (in the appropriate
sense - minimizing the ordered pair of electric and hyperbolic distances). Since $\phi$ 
is an electric isometry, but a hyperbolic quasi-isometry, there exist
$C > 0$ (uniform constant) and $K=K(B)$  such that \\
a) for all $x \in
\lambda_i$, $\phi(x)$ lies in a ($d_G-$) $C$-neighborhood and
a $d_{wel}-$ $K$-neighborhood of $\Phi
(\lambda_{i} ) = \lambda_{i+1}$.\\
b)   for all $z \in \til{S}$, $d_{G}(\pi_i(z,i), \pi_{i+1}(z,i+1))\leq C$ (by Lemma \ref{retract-thin} or Theorem \ref{retract})
and $d_{wel}(\pi_i(z,i), \pi_{i+1}(z,i+1))\leq K$ (by Lemma \ref{retract-thick}).\\
The last statement follows from the fact that the block $B$ is topologically a product and hence the map $\phi$ is a quasi-isometry,
with quasi-isometry constants depending on $B$.

We re-emphasize here that $C$ is independent of both the split block $B$ and the geodesic $\lambda$ (and hence the ladder ${\LL}_\lambda$),
whereas $K=K(B)$ depends on the split block $B$ but is independent of the geodesic $\lambda$.

 The same holds for $\phi^{-1}$ and
points in $\lambda_{i+1}$, where $\phi^{-1}$ denotes the {\em
  quasi-isometric inverse} of $\phi$ from
$\widetilde{S} \times \{ i + 1 \}$ to $\widetilde{S} \times \{ i \}
$. It is worth pointing out here that Remark \ref{alt} will be used later to pull back information from the graph metric in $(M_{wel},d_G)$
to the model manifold $(\til{M},d_M)$ and hence via the biLipschitz homeomorphism $F^{-1}$ to $\til{N}$
to give information in the hyperbolic metric.

Again, since $\lambda_{i}$ and $\lambda_{i+1}$ are electro-ambient
quasigeodesics, we further note that  for all $(x,i) \in \lambda_i$,
$(x,i+1) \in N_K(\lambda_{i+1},d)$, where $d$ is the (biLipschitz) hyperbolic metric on $\til S$.

The {\bf
  ${\LL}_\lambda$-elementary admissible paths} in $\widetilde{B}$ consist
of the following:\\
1) Horizontal  subsegments of $\lambda_j$,  $j = \{ i, i+ 1 \}$.\\
2) Vertical segments  joining $x  \times \{ i  \}$ to
$x  \times \{ i+1 \}$. These have $d_{wel}-$ `thickness' $l=l(B)$ and $d_G-$ thickness one, by Lemma
\ref{bddheight}.\\
3) Horizontal geodesic segments lying in a {\em (biLipschitz) hyperbolic}
 $K (=K(B))$-neighborhood of
  $\lambda_j$, $j = i, i+1$.\\
4) Horizontal (biLipschitz) hyperbolic segments of {\em electric length $\leq C$}
  and {\em  (biLipschitz) hyperbolic length $\leq K(B)$} joining points of the form
  $( \phi (x), i+1)$ to a point on $\lambda_{i+1}$ for $x \in
  \lambda_{i}$. \\
5) Horizontal (biLipschitz) hyperbolic segments of {\em electric length $\leq C$}
  and {\em (biLipschitz) hyperbolic length $\leq K(B)$} joining points of the form
  $( \phi^{-1} (x), i)$ to a point on $\lambda_{i}$ for $x \in
  \lambda_{i+1}$. \\

\noindent {\bf Definition:} An   ${\LL}_\lambda$-admissible path is a continuous path that can be decomposed as a union of a sequence of
${\LL}_\lambda$-elementary admissible paths with disjoint interiors.

The next lemma follows from the above definition and Lemma
\ref{bddheight}.

\begin{lemma}
There exists a function $g: \mathbb{Z} \rightarrow \mathbb{N}$ such that
for any block $B_i$, and $x$ lying on an ${\LL}_\lambda$-admissible
path in $\widetilde{B_i}$, there exist $y
\in \lambda_{i}$ and $z \in \lambda_{i+1}$ such that

\begin{center}
$d_{wel}(x, y) \leq g(i)$ \\
$d_{wel}(x, z) \leq g(i)$ \\
$d_{M}(x, y) \leq g(i)$ \\
$d_{M}(x, z) \leq g(i)$ \\
\end{center}
\label{bddheight-adm}
\end{lemma}

 The following is an easy Corollary of Lemma \ref{bddheight-adm} above,

\begin{cor}
There exists a function $h: \mathbb{Z} \rightarrow \mathbb{N}$ such that
for any block $B_i$, and $x$ lying on a ${\LL}_\lambda$-admissible
path in $\widetilde{B_i}$, there exist $y
\in \lambda_0 = \lambda$ such that:

\begin{center}

$d_{wel}(x, y) \leq h(i)$ \\
$d_M(x, y) \leq h(i)$

\end{center}
\label{bddht-final}
\end{cor}

\begin{proof}
Let $h(i) = \Sigma_{j = 0 \cdots i} g(j)$ be the sum of the values of
$g(j)$ as $j$ ranges from $0$ to $i$ (with the assumption that
increments are by $+1$ for $i \geq 0$ and by $-1$ for $i \leq
0$). \end{proof}

\noindent {\bf  Note:} In  Lemma \ref{bddheight-adm} and
Corollary \ref{bddht-final}, it is important to note
that the distance $d_{wel}$ (resp. $d_M)$ is the weld (resp. model)  metric, not the graph metric. 
This is because the lengths occurring in  ${\LL}_\lambda$-elementary admissible paths  of types (4) and (5) above are (biLipschitz) hyperbolic lengths depending only on $i$ (in $B_i$).

\begin{lemma}
There exists a function $M(N) : \mathbb{N} \rightarrow \mathbb{N}$
such that  $M (N) \rightarrow \infty$
 as $N \rightarrow \infty$ for which the following holds:\\
For any geodesic $\lambda \subset \widetilde{S} \times \{ 0 \} \subset
\widetilde{B_0}$, a fixed reference point 
 $p \in \widetilde{S} \times \{ 0 \} \subset
\widetilde{B_0}$ and any $x$ on an ${\LL}_\lambda$-admissible path, 

\begin{center}

$d(\lambda , p) \geq N \Rightarrow d_{wel}(x,p) \geq M(N)$ and $ d_{M}(x,p) \geq M(N)$.

\end{center}

\label{far-nopunct}
\end{lemma}

\begin{proof} Suppose that $\lambda$ lies outside $B_N(p)$, the $N$-ball about
a fixed reference point $p$ on the boundary horizontal surface
$\widetilde{S} \times \{ 0 \} \subset \widetilde{B_0}$. Then by
Corollary \ref{bddht-final}, any $x$ lying on an ${\LL}_\lambda$-admissible
path in $\widetilde{B_i}$ satisfies

\begin{center}

$d_{wel}(x, p) \geq N - h(i).$ 

\end{center}

Also, since the electric, and $d_{wel}-$`thickness' (the
shortest distance between its boundary horizontal sheets) is $\geq k_0$ (by uniform $k_0-$ separatedness
of  horizontal sheets),
we get,

\begin{center}

$d_{wel}(x, p) \geq |i|k_0$ 

\end{center}

Assume for convenience that $i \geq 0$ (a similar argument works,
reversing signs for $i < 0$). Then,

\begin{center}

$d_{wel}(x, p) \geq $ ${\rm min}_i$ $max \{ ik_0,  N - h(i) \}$ 

\end{center}

Let $h_1 (i) = h(i) + ik_0$.  Then $h_1$ is a
monotonically increasing function on the integers.
If
 $M(N)=h_1^{-1} (N)$ denote the largest positive
integer $n$ such that $h_1(n) \leq N$, then clearly, 
 $M (N) \rightarrow \infty$
 as $N \rightarrow \infty$. Also, $d_{wel}(x,p) \geq k_0 M(N)$ and the first conclusion of the Lemma follows.

The same arguments work for $(\til{M}, d_M)$. \end{proof}

\subsection{Projecting to ${\mathcal{L}}_\lambda$ and Joining the Dots} \label{join-sxn}
\begin{defn}
An {\bf ${\LL}_\lambda$ admissible ($d_G$) $K-$ quasigeodesic} is an
${\LL}_\lambda$ admissible path that is a $K-$quasigeodesic in
$(\widetilde{M_{wel}}, d_{G})$. 
\end{defn}

Our strategy in this subsection is to project the intersection of an
 admissible quasigeodesic (Lemma \ref{admqg}) with the horizontal sheets $\til{M_H}$ onto
 ${\LL}_\lambda$ and then obtain a {\it connected} ${\LL}_\lambda$-admissible quasigeodesic from it
by interpolating ${\LL}_\lambda$-admissible paths. We think of this last step as {\it "joining the dots"}. 
The end product is thus a {\it connected} $d_G$-quasigeodesic built up of ${\LL}_\lambda$
admissible paths.

\begin{lemma}
There exists $K \geq 1$ and a function $M(N) : \mathbb{N} \rightarrow \mathbb{N}$
with  $M (N) \rightarrow \infty$
 as $N \rightarrow \infty$ such that the following holds:\\
Let $B_0$ denote the first block (thick or split) in $M_{wel}$ and let $S \times \{0\}$ denote its lower  boundary.
For a fixed reference point 
 $p \in \widetilde{S} \times \{ 0 \} \subset
\widetilde{B_0}$, and any geodesic $\lambda \subset \widetilde{S} \times \{ 0 \} \subset
\widetilde{B_0}$, 
there exists an  ${\LL}_\lambda$ admissible ($d_G$) $K-$ quasigeodesic $\beta_{adm} \subset \til{M_{wel}}$  without 
backtracking, such that \\
$(1)$ $\beta_{adm}$ joins
the end-points of $\lambda$. \\
$(2)$ 
$d(\lambda , p) \geq N \Rightarrow d_{wel}(\beta_{adm},p) \geq M(N)$. 
\label{adm-qgeod-props}
\end{lemma}

\begin{proof} Let $a, b$ denote the end-points of $\lambda$.
First, by Lemma \ref{admqg} there exists an admissible $d_G$-geodesic $\beta_e  \subset \til{M_{wel}}$ joining  $a, b$.
We now look at $\Pi_\lambda(\beta_e \cap \widetilde{M_H})$ obtained by acting
 on $\beta_e \cap \widetilde{M_H}$ by $\Pi_\lambda$. From Theorem \ref{retract}, we
shall  conclude that the image $\Pi_\lambda ( \beta_e \cap \widetilde{M_H})$  is a
$d_G$ quasigeodesic carried by ${\LL}_\lambda$ in an appropriate sense as explicated below.


Since $\beta_e$ is itself an admissible path, there is a sequence of points $a=a_1, b_1, a_2, b_2, \cdots , a_k, b_k=b$ such that
the piece of $\beta_e$ joining $a_i$ to $b_i$ is horizontal, whereas the piece of $\beta_e$ joining $b_i$ to $a_{i+1}$ is vertical.
In particular, $b_i$ and $a_{i+1}$ must lie in the same split component if they lie in (the universal cover of) a split block. 
In this case ${\Pi}_\lambda (b_i)$ and ${\Pi}_\lambda (a_{i+1})$ must also lie in the same split component. We shall now join 
the sequence of points ${\Pi}_\lambda (a_1), {\Pi}_\lambda (b_1), {\Pi}_\lambda (a_2), {\Pi}_\lambda (b_2), \cdots , 
{\Pi}_\lambda (a_k), {\Pi}_\lambda (b_k)$ by  horizontal and vertical segments
to obtain an ${\mathcal{L}}_\lambda$-admissible path
$\beta_{adm}$ as follows:\\
 For all $i$, $[{\Pi}_\lambda (a_i), {\Pi}_\lambda (b_i)] $ will be  a geodesic in the horizontal sheet $\til{S_i}$,
joining ${\Pi}_\lambda (a_i), {\Pi}_\lambda (b_i)$. 

The ${\LL}_\lambda$-admissible path joining 
$ {\Pi}_\lambda (b_i), {\Pi}_\lambda (a_{i+1})$ requires more care to define. For notational simplicity, let 
$ b_i =p$ and $ a_{i+1}=q$.\\
1) Let $[p,q]$ be a  vertical segment in a thick block joining $p, q$. Then $\Pi_\lambda
(p), \Pi_\lambda (q)$ are a uniformly bounded $d_{wel}-$ distance apart by Theorem \ref{retract}. Hence, by
 Lemma \ref{retract-thick}, we can join 
 $\Pi_\lambda
(p), \Pi_\lambda (q)$ by an ${\LL}_\lambda$-admissible path of length
bounded by some $C_0$ (independent of $B$, $\lambda$). 

For a thick block, we define the ${\LL}_\lambda$-admissible path joining $\Pi_\lambda
(p), \Pi_\lambda (q)$ to be any such ${\LL}_\lambda$-admissible path of uniformly bounded  $d_{wel}-$length. \\
2) Let $[p,q]$ be a vertical segment in a split block $\til{B_i} $ 
of $d_G$ length one and  $d_{wel}-$ length $\leq l_i$  joining $p, q$, where $ p \in \til{S_i}$, the lower horizontal
boundary of $\til{B_i} $  and $ q \in \til{S_{i+1}}$, the upper horizontal
boundary of $\til{B_i} $. Since
 $p,q$  lie within a 
split component,
 $d_G(\Pi_\lambda(p), \Pi_\lambda(q)) = 1$, that is to say $\Pi_\lambda(p), \Pi_\lambda(q)$ also
lie within
a split component.
This is because the projection of a split component lies within
a single split component. Hence there exists an admissible path
$[\Pi_\lambda(p), \Pi_\lambda(q)]$ of $d_G$ length one joining $\Pi_\lambda(p), \Pi_\lambda(q)$. 
Further, by Lemma \ref{retract-thick} again, we can join 
 $\Pi_\lambda
(p), \Pi_\lambda (q)$ by an ${\LL}_\lambda$-admissible path of $d_{wel}-$length
bounded by some $C_i$ (dependent on $B_i$ but independent of $\lambda$). Note  that,
since $C_i$ depends on $B_i$, it depends on $l_i$ in particular.\\
3) By Remark \ref{disamb} the two images under nearest projection of a point in $\til{S_i}$ onto  
respectively a hyperbolic geodesic and an electro-ambient quasigeodesic 
in  $\til{S_i}$  (joining
any pair of points) 
 are a uniformly bounded (biLipschitz)-hyperbolic distance apart. Hence, by
  Lemma \ref{disamb0}, we can join them
  by an ${\LL}_\lambda$-admissible path of length
bounded by some uniform $C_1$ (independent of $B_i$, $\lambda$). 

A clarificatory remark as to why segments of type (3) are necessary: In defining ${\LL}_\lambda$, we have had
to make a choice. Suppose $\lambda_i \subset {\til{S}}_i$. Then $S_i$ is the common boundary of two blocks. In case
both are split blocks then there is a choice of $\lambda_i$ out
of two electro-ambient quasigeodesics involved. If one is a split block and the other a thick block, 
then there is a choice of $\lambda_i$  involved out of an  electro-ambient quasigeodesic and a geodesic. The  different
nearest point projections corresponding to the different choices of  $\lambda_i$  differ by a uniformly bounded
amount (Remark \ref{disamb}). Segments of type (3) take care of this bounded discrepancy.\\

For a split block, we define the ${\LL}_\lambda$-admissible path joining $\Pi_\lambda
(p), \Pi_\lambda (q)$ to consist of one ${\LL}_\lambda$-admissible path constructed
in Step (2) above and (at most) two segments of uniformly bounded  $d_{wel}-$length as in Step (3). Thus
an ${\LL}_\lambda$-admissible path joining $\Pi_\lambda
(p), \Pi_\lambda (q)$ contains one vertical segment of type (2) typically sandwiched between two segments of type (3).\\

Joining ${\Pi}_\lambda (a_i), {\Pi}_\lambda (b_i) $ by $[{\Pi}_\lambda (a_i), {\Pi}_\lambda (b_i)] $
and $ {\Pi}_\lambda (b_i), {\Pi}_\lambda (a_{i+1})$ by ${\LL}_\lambda$-admissible paths as above, we obtain the required 
${\LL}_\lambda$ admissible ($d_G$) $K-$ quasigeodesic $\beta_{adm} \subset \til{M_{wel}}$.

  By Theorem \ref{retract},
there exists $K \geq 1$ such that $\beta_{adm}$
 represents a  ($d_G$)-$K $-quasigeodesic. This proves statement (1) of the Lemma.

After "joining the dots" by ${\LL}_\lambda$-admissible paths as above, we can assume further that the ${\LL}_\lambda$-admissible  quasigeodesic
$\beta_{adm}$
thus obtained does not backtrack relative to split components.
 Conclusion (2) of the Lemma now follows from Lemma \ref{far-nopunct}
since we have obtained an admissible quasigeodesic built up out of  ${\LL}_\lambda$-admissible paths. 
\end{proof}

\subsection{Recovering Electro-ambient Quasigeodesics I}\label{recovery1}

This subsection is devoted to extracting an electro-ambient
quasigeodesic $\beta_{ea}$ in $(\til{M_{wel}}, d_G)$
from an ${\LL}_\lambda$-admissible quasigeodesic $\beta_{adm}$. 
 $\beta_{ea}$ shall satisfy  the property indicated by Lemma
\ref{adm-qgeod-props} above. 

\begin{lemma}
There exists $\kappa \geq 1$ and a function $M^{\prime}(N) : \mathbb{N} \rightarrow \mathbb{N}$
with  $M^{\prime} (N) \rightarrow \infty$
 as $N \rightarrow \infty$ such that the following holds:\\
Let $B_0$ denote the first block (thick or split) in $M_{wel}$ and let $S \times \{0\}$ denote its lower  boundary.
For a fixed reference point 
 $p \in \widetilde{S} \times \{ 0 \} \subset
\widetilde{B_0}$, and any geodesic $\lambda \subset \widetilde{S} \times \{ 0 \} \subset
\widetilde{B_0}$, 
there exists an 
  electro-ambient $\kappa$-quasigeodesic $\beta_{ea}$  without 
backtracking in $(\til{M_{wel}}, d_G)$, such that \\
$\bullet$ $\beta_{ea}$ joins
the end-points of $\lambda$. \\
$\bullet$ 
$d(\lambda , p) \geq N \Rightarrow d_{wel}(\beta_{ea},p) \geq M^{\prime}(N)$. \\
\label{adm-ea-props}
\end{lemma}

\begin{proof} From Lemma \ref{adm-qgeod-props}, we have an
${\LL}_\lambda$ - admissible $\kappa$-quasigeodesic $\beta_{adm}$ without backtracking (with respect to the collection $\KK$
of  split components  $\til{K} $) and a function $M(N)$
 satisfying
the conclusions of the Lemma.
Since $\beta_{adm}$ does not backtrack, we can decompose it as a  union
of non-overlapping segments $\beta_1, \cdots \beta_k$, such that only successive $\beta_i$'s intersect at one common end-point and each
$\beta_i$ is \\
a) either an ${\LL}_\lambda$- admissible quasigeodesic lying outside 
split components, \\
b) or 
an ${\LL}_\lambda$-admissible quasigeodesic lying entirely
within some split component ${\widetilde{K}}_{n(i)}$. Further, since
$\beta_{adm}$ does not backtrack  relative to split components, we can assume that all ${\widetilde{K}}_{n(i)}$'s are
distinct, i.e. $i\neq j \Rightarrow {\widetilde{K}}_{n(i)} \neq {\widetilde{K}}_{n(j)}$.

We  modify $\beta_{adm}$ to an electro-ambient quasigeodesic
$\beta_{ea}$ in $(\til{M_{wel}},d_G)$ as per the following recipe:\\
1) $\beta_{ea}$ coincides with $\beta_{adm}$ outside split
  components. \\
2) If some $\beta_i$ lies within a split component
  ${\widetilde{K}}_{n(i)}$, then
we replace it by a geodesic $\beta_{i}^{ea}$ {\it in the intrinsic metric on ${\widetilde{K}}_{n(i)}$}
joining the end-points of $\beta_i$. Of course  $\beta_{i}^{ea}$ lies
 within   ${\widetilde{K}}_{n(i)}$.

 Since $\beta_{ea}$ coincides with $\beta_{adm}$ outside split
  components  and since $\beta_{adm}$ is  a ($d_G$) $\kappa$-quasigeodesic, therefore $\beta_{ea}$ represents a ($d_G$) $\kappa$-quasigeodesic. 
Hence, the resultant path $\beta_{ea}$ is  an  electro-ambient
$\kappa$-quasigeodesic without backtracking.

Next, since any amalgamation component of $\til{S}$ is
   quasiconvex   in the split component $\til{K}$ containing it, each segment
$\beta_i^{ea}$ lies in a $C_i$ neighborhood of $\beta_i$. Here $C_i$
depends on the quasiconvexity constants of the amalgamation components in split components and hence only on
  the thickness $l_i$ of the split component $K_{n(i)}$.

We let $C(m)$ denote the maximum of the (finitely many) values of $C_i$ for the split components 
of  
$\til{B_m}$, where we take $C(m)=0$ if $B_m$ is thick (this makes sense as $\beta_{ea}$ coincides with
 $\beta_{adm}$ outside split
  components). 
Then, as in the proof of Lemma \ref{far-nopunct}, we have for any $z
\in \beta_{ea} \cap \til{B_m}$, 

\begin{center}
$ d(z, p) \geq $ max $(mk_0, M(N) - C(m))$
\end{center}

Again, as in Lemma \ref{far-nopunct}, this gives us a (new) function 
$M^{\prime}(N) : \mathbb{N} \rightarrow \mathbb{N}$
such that  $M^{\prime} (N) \rightarrow \infty$
 as $N \rightarrow \infty$ for which
$$d(\lambda , p) \geq N \Rightarrow d_{wel}(\beta_{ea},p) \geq M^{\prime}(N).$$

This proves the Lemma. \end{proof}

\subsection{Recovering Electro-ambient Quasigeodesics II}\label{recovery2}

This subsection is devoted to extracting an electro-ambient
quasigeodesic $\beta_{ea2}$ in ${\til{M}}_2= (\til{M}, d_{CH})$
from an electro-ambient
quasigeodesic $\beta_{ea}$ in ${\til{M}}_1=(\til{M_{wel}}, d_{G})$.
 $\beta_{ea2}$ shall satisfy  the property indicated by Lemmas \ref{adm-qgeod-props} and 
\ref{adm-ea-props} above.

Recall that ${\til{M}}_2= (\til{M}, d_{CH})$ denotes $\til{M}$ with the 
electric metric obtained by electrocuting the convex hulls $CH({\til K})$ of extended split components
$\til K$. 
Also, recall that
 an   electro-ambient $k$-quasigeodesic $\gamma$
in $(\til{M}, d_{CH})$  is
 a $k-$ quasigeodesic in $(\til{M}, d_{CH})$ such that in an ordering (from the left)
of the convex hulls of split components
that $\gamma$ meets, each $\gamma \cap CH({\til K})$ is a 
geodesic in the {\it intrinsic} metric on $CH({\til K})$ (which in turn is uniformly bi-Lipschitz to the hyperbolic metric
on $CH({\til K})$ under the bi-Lipschitz homeomorphism between the model manifold $M$ and the hyperbolic manifold $N$).

The underlying sets $\til{M_{wel}}$ (for ${\til{M}}_1$)
and $\til M$
(for ${\til{M}}_2$)  are homeomorphic as topological spaces.
Also, ${\til{M}}_1$ is obtained by electrocuting the welded metric, i.e. $(\til{M_{wel}}, d_{wel})$, whereas 
${\til{M}}_2$ is obtained by electrocuting the model metric, i.e. $(\til{M}, d_M)$.
Note further that the metrics $(\til{M}, d_{wel})$ and  $(\til{M}, d_M)$  coincide off Margulis tubes.

We need to set up a correspondence now between paths
in $(\til{M_{wel}}, d_{wel})$ and  $(\til{M}, d_M)$, and hence between ${\til{M}}_1 = (\til{M}, d_{G})$ and ${\til{M}}_2=(\til{M}, d_{CH})$.

\begin{rmk} {\rm   Paths $\alpha_i
\subset {\til{M}}_i$ ($i=1,2$) are said to {\bf correspond} if\\
1) They coincide off Margulis tubes \\
2) Each piece of $\alpha_2$ inside a (closed) Margulis tube is a geodesic in the model metric $d_M$. \\
It follows that any path $\alpha_1
\subset {\til{M}}_1$ corresponds to a unique $\alpha_2
\subset {\til{M}}_2$.}
\label{corr}
\end{rmk}

\begin{lemma} 
There exists $\kappa \geq 1$ and  a function $M^{\prime}(N) : \mathbb{N} \rightarrow \mathbb{N}$
such that  $M^{\prime}(N) \rightarrow \infty$
 as $N \rightarrow \infty$ for which the following holds:\\
Let $B_0$ denote the first block (thick or split) in $M_{wel}$ and let $S \times \{0\}$ denote its lower  boundary.
For a fixed reference point 
 $p \in \widetilde{S} \times \{ 0 \} \subset
\widetilde{B_0}$, and any geodesic $\lambda \subset \widetilde{S} \times \{ 0 \} \subset
\widetilde{B_0}$, 
there exists an
  electro-ambient $\kappa-$quasigeodesic $\beta_{ea}$  without 
backtracking in $(\til{M_{wel}}, d_G)$ and a path $\beta_{ea1}$ {\bf corresponding} to $\beta_{ea}$ in $(\til{M}, d_{CH})$, such that \\
$(1)$ $\beta_{ea1}$ joins
the end-points of $\lambda$. \\
$(2)$ 
$d(\lambda , p) \geq N \Rightarrow d(\beta_{ea1},p) \geq M^{\prime}(N)$. \\
\label{far-eas}
\end{lemma}

\begin{proof}
By Lemma \ref{adm-ea-props} (Section \ref{recovery1} above)  there exists an electro-ambient
$\kappa_0-$quasigeodesic $\beta_{ea}$  in 
${\til{M}}_1=(\til{M_{wel}}, d_{G})$ joining the end-points of $\lambda$ (where $\kappa_0$ is independent of $\lambda$).
By  Remark \ref{corr},  $\beta_{ea}$ {\it corresponds} to a unique path, which we call $\beta_{ea1}$,
in ${\til{M}}_2$. $\beta_{ea1}$ is  obtained by replacing intersections of $\beta_{ea}$ with tube-electrocuted Margulis tubes by
hyperbolic geodesics lying in the corresponding Margulis tubes as per   Remark \ref{corr}.
From Lemma \ref{qi12}, $(\til{M}, d_{CH}) (={\til{M}}_2)$
 is quasi-isometric to $(\til{M}, d_{G})$. Hence there exists $\kappa\geq 1$ such that for
any $\lambda$, the path $\beta_{ea1}$ is a $\kappa$-quasigeodesic in
${\til{M}}_2$.

Also by Lemma \ref{adm-ea-props}, there exists a function $M(N) : \mathbb{N} \rightarrow \mathbb{N}$
such that  $M(N) \rightarrow \infty$  as $N \rightarrow \infty$ for which the following holds:\\ If $d(\lambda , p) \geq N$, then
 $\beta_{ea}$ lies outside a large $M(N)$-ball about $p$ in $({\til{M}}_{wel},d_{wel})$.

It follows that the intersection
of $\beta_{ea}$ with the boundary $\partial \til{T}$ of the lift
$\til T$ of any Margulis tube ${T}$ lies outside an $M(N)$-ball about $p$. Each point 
$x \in \beta_{ea} \cap \partial \til{T}$ lies on a unique totally geodesic hyperbolic disk $D_x \subset  \til{T}$. Also, 
$\beta_{ea1} \cap  \til{T}\subset \bigcup_{x \in \beta_{ea} \cap \partial \til{T}}D_x$ by the convexity of $\bigcup_{x \in \beta_{ea} \cap \partial \til{T}}D_x$.
Let the maximum diameter of  Margulis tubes intersecting the $i$th block in $\til M$ be $t_i$. 
Then $d_M(\beta_{ea1} \cap \til{B_i} , p) \geq d_{wel}(\beta_{ea} \cap \til{B_i} , p) - t_i
\geq M(N) - t_i$. Now,  a reprise of the argument in Lemma \ref{far-nopunct}
shows that $\beta_{ea1}$ lies outside a large $M^{\prime}(N)$ ball about $p$,
where   $M^{\prime}(N) \rightarrow \infty$
 as $N \rightarrow \infty$. \end{proof}

To obtain an electro-ambient quasigeodesic $\beta_{ea2}$ in $(\til{M}, d_{CH})$ from $\beta_{ea1}$, first observe that
there exists $D_0$ such that the diameter in the  $d_G$ metric
$dia_G(\beta_{ea1}\cap CH({\til K}) ) \leq D_0$  for any $CH({\til K})$. This follows from the fact that $\beta_{ea1}$ is a 
$\kappa$-quasigeodesic in $(\til{M}, d_{G})$
and from Lemma \ref{qi12}, which says that $(\til{M}, d_{CH})$ and $(\til{M_{wel}}, d_{G})$ are quasi-isometric.

\begin{lemma} For every $D_0 \geq 0$ and split component $\til{K} \subset \til{M_{wel}}$, there exists $D_1 \geq 0$
such that the following holds:\\
Let $\alpha \subset CH({\til K})  \subset \til{M}$ be a path such that the path $\eta$ in   $\til{M_{wel}}$
corresponding to it is of length at most $D_0$ in the $d_G$ metric. 
Further suppose that \\
a) $\alpha \cap \til{C} ( \subset \til{M})$
for any split component $\til{C}$ is a geodesic in the intrinsic metric on $\til C$\\
b) $\alpha \cap \mathbb{T}$ is a hyperbolic geodesic
for any lift $\mathbb{T}$ of a Margulis tube.\\
 Let $\gamma = [a,b]$ be the (model) hyperbolic geodesic in $(\til{M},d_M)$  joining the end-points $a,b$ of $\alpha$. Then 
 $\gamma$ lies in a ($d_M-$) $D_1$ neighborhood of $\alpha$.
\label{g2ch}
\end{lemma}

\begin{proof} Note first that the complement in $\til M$ 
of the union of  split components is the union of the universal covers of thick blocks and
 Margulis tubes.
Hence by the hypotheses $\alpha$ can be described as the union of at most $3D_0$ pieces $\alpha_1, \cdots , \alpha_j (j \leq 3D_0)$, such that
each $\alpha_i$ is either a geodesic in the intrinsic metric on $\til C$ for some split component $\til C$, or a  geodesic  in $(\til{M},d_M)$.

 Let $\beta_i$ be the  geodesic  in $(\til{M},d_M)$  joining the
 end-points of $\alpha_i$. Then
$d(\gamma , \cup_i \beta_i) \leq j \delta_0 \leq 3D_0 \delta_0$, where $\delta_0$ is the (Gromov) hyperbolicity constant of $\til{M}$.

Since $\alpha$ meets a bounded number of split components, there exists $C_1 \geq 0$ such that each split component $\til C$ 
that  $\alpha$ meets  is $C_1$-quasiconvex.
Note that $C_1$ depends only on the convex hull $CH({\til K})$ 
 and the fact that any $CH({\til K})$ meets the lifts of only
a uniformly bounded number of split components
by graph quasiconvexity (Theorem \ref{gqc}). Hence for any $\alpha_i \subset \til{C}$, $d_M(\alpha_i,\beta_i) \leq C_1$. Choosing
$D_1 = C_1 + 3D_0 \delta_0$, we are through. \end{proof}

We are now in a position to obtain the last `recovery' Lemma of this section.
The main part of the argument is again a reprise of a similar argument in Lemma \ref{far-nopunct}.
We shall recount it briefly for completeness.

\begin{lemma}
There exists $\kappa \geq 1$
and  a function $M_0(N) : \mathbb{N} \rightarrow \mathbb{N}$
such that  $M_0(N) \rightarrow \infty$
 as $N \rightarrow \infty$ for which the following holds:\\
Let $B_0$ denote the first block (thick or split) in $M_{wel}$ and let $S \times \{0\}$ denote its lower  boundary.
For a fixed reference point 
 $p \in \widetilde{S} \times \{ 0 \} \subset
\widetilde{B_0}$, and any geodesic $\lambda \subset \widetilde{S} \times \{ 0 \} \subset
\widetilde{B_0}$,
there exists an 
  electro-ambient $\kappa-$quasigeodesic $\beta_{ea2}$  without 
backtracking in $(\til{M}, d_{CH})$, such that \\
$(1)$ $\beta_{ea2}$ joins
the end-points of $\lambda$. \\
$(2)$ 
$d(\lambda , p) \geq N \Rightarrow d_M(\beta_{ea2},p) \geq M_0(N)$. \\
\label{adm-ea2-props}
\end{lemma}

\begin{proof} By Lemma \ref{far-eas}, there exists $\kappa_0$  
and  a function $M^\prime(N) : \mathbb{N} \rightarrow \mathbb{N}$
such that for any
geodesic $\lambda \subset \widetilde{S} \times \{ 0 \} \subset
\widetilde{B_0}$ with $d(\lambda , p) \geq N$ 
there exists a path $\alpha$ in $(\til{M}, d_{CH})$ {\it corresponding} (as per Remark \ref{corr}) to an electro-ambient quasigeodesic in 
$(\til{M}, d_{G})$ satisfying the following:\\ 
a)  $\alpha$ joins
the end-points of $\lambda$. \\
b)  $d_M(\alpha,p) \geq M^{\prime}(N)$.\\
c) $N \rightarrow \infty \Rightarrow M^{\prime}(N) \rightarrow \infty $.\\

Let $\beta_{ea2}$ be an electro-ambient quasigeodesic in $(\til{M}, d_{CH})$ joining the end-points of $\alpha$.
Let $\mathcal{CH}(\til{\KK})$ be the collection of (images under the biLipschitz homeomorphism $F$ of)
convex hulls of extended split components. 

Recall that $\beta_{ea2}$ is obtained 
by looking at the intervals of intersection
of $\alpha$ with $CH(\til{K}) \in \mathcal{CH}(\til{\KK})$, ordered from the left, and replacing maximal  intersections
with (model) hyperbolic geodesics in $CH(\til{K}) $.

Let $x \in \beta_{ea2} \cap CH(\til{K}) $ for an extended split component $\til{K}$. Then by construction of the electro-ambient quasigeodesic
$\beta_{ea2}$ from $\alpha$ and Lemma \ref{g2ch} there exists $y \in \alpha \cap CH(\til{K}) $ and $D_1 = D_1(K)$ such that $d(x,y) \leq D_1$.

 By uniform graph quasiconvexity (Theorem \ref{gqc}), for each $i$ there exist finitely many extended split components
$K$ such that 
$ \til{B_i}  \cap CH(\til{K}) \neq \emptyset$. Let $D_i$ be the maximum value of the $D_1(K)$'s for these split components.
Hence $x \in \beta_{ea2}\cap\til{B_i} $ implies that $ d(x,p) \geq M^{\prime}(N) - D_i$. Also, by uniform $k_0$-separatedness of split surfaces, 
$x \in \til{B_i} $ implies that $ d(x,p) \geq ik_0$. Therefore\\
\begin{center}
 $ d(\beta_{ea2},p) \geq$ ${\rm{min}}_i$ max $(ik_0,M^{\prime}(N) -\sum_{j\leq i} D_j) $
\end{center}

Defining $M_0(N)$ to be $M_0(N) = {\rm{min}}_i$ max $(ik_0,M^{\prime}(N) - \sum_{j\leq i} D_j) $, and observing that   $M_0(N) \rightarrow \infty$
 as $N \rightarrow \infty$ (by the same argument as in Lemma \ref{far-nopunct}) we are through.
 \end{proof}

\subsection{Application to Sequences of Surface Groups} \label{seqs}
The main Proposition of this subsection will be used in \cite{mahan-series2}. 

The {\it proof} of  Lemma \ref{adm-ea2-props}
gives the following.

\begin{cor} 
Let $D$ be a positive integer. Let $B_{-D}, \cdots , B_0, \cdots , B_n, \cdots , B_{n+D}$
be a collection of split blocks and let $\BB_n^1$ be the union of these blocks glued along the common
boundary split surfaces (i.e. $B_{i-1}$ is glued to $B_{i}$ along $S_i$). We assume that this gluing can be done consistently (i.e. the Margulis tubes are compatible). Let $\BB_n =\bigcup_1^n B_i \subset \BB_n^1$.
Let $M$ be a manifold  of split geometry (not necessarily simply or doubly degenerate, i.e. we allow
$M$ to have finitely many split blocks), such that each split component
is $D$-graph quasiconvex and $\BB_n^1 \subset M$.  Then for all  $L \geq 0$ there exists  $N \geq 0$  such that
the following holds.\\
For all geodesic segments
 $\lambda$  lying outside an $N$-ball
around ${o}\in{\tilde{S_0}}$  and any  electro-ambient quasigeodesic $\beta_{ea2}^n$  without 
backtracking in $\til{M}$ joining
the end-points of $\lambda$, $\beta_{ea2}^n \cap \widetilde{\BB_n}$ lies outside the $L$-ball around 
$o \in \tilde{M}$.
\label{blocks}
\end{cor}

Corollary \ref{blocks} above will be used to prove the convergence of Cannon-Thurston maps for quasi-Fuchsian groups 
converging strongly to a simply degenerate group.

\begin{rmk} {\rm In Corollary \ref{blocks}, we could replace $\BB_n^1$ by $\BB_n^2 = B_{-n-D}, \cdots , B_0, \cdots , B_n, \cdots , B_{n+D}$
and the same conclusions follow. This will be used 
to prove the convergence of Cannon-Thurston maps for quasi-Fuchsian groups 
converging strongly to a doubly degenerate group.} \label{blocks2} \end{rmk}

\section{Cannon-Thurston Maps for Surfaces Without Punctures} \label{ct-nopunct}

We note the following properties of   $(\til{M},d_G)$ and $\KK$
where  $(\til{M},d_G)$  is the graph model of $\widetilde{M}$ and $\mathcal{K}$
consists of the split components.
 There exist $C, D, \Delta$ such that\\
1) Each split component  is $C$-graph quasiconvex by Theorem \ref{gqc}.\\
2)  $(\til{M},d_G)$   is $\Delta$-hyperbolic.\\
3) Given $K, \epsilon$, there exists $D_0$ such that if $\gamma$ 
be  a $(K, \epsilon)$ quasigeodesic   in  $(\til{M},d_M)$  joining $a, b$ and if
$\beta$ be a $(K, \epsilon)$ electro-ambient quasigeodesic  in  $(\til{M},d_G)$   joining $a,
b$, then $\gamma$ lies in a $D_0$ neighborhood of $\beta$  in  $(\til{M},d_M)$. This follows from Lemma \ref{ea-strong}.

We shall now assemble the proof
of the main Theorem.

\begin{theorem}
Let $M$ be a simply or doubly degenerate hyperbolic 3 manifold without parabolics, homeomorphic to $S \times J$ (for $J = [0,
  \infty ) $ or $( - \infty , \infty )$ respectively). Fix a base surface $S_0 = S \times \{ 0 \}$. Then the inclusion
  $i : \widetilde{S_0} \rightarrow \widetilde{M}$ extends continuously
  to a map between the compactifications
  $\hat{i} : \widehat{S_0} \rightarrow \widehat{M}$. Hence the limit set
  of $\widetilde{S_0}$ is locally connected.
\label{crucial}
\end{theorem}

\begin{proof} By Theorem \ref{gqc}, $M$ has split geometry and we may assume that $S_0 \subset B_0$, the first block. Let
$(\til{M}, d_{CH})$ and $(\til{M}, d_{G})$ be as above and let $d_M$ be the model metric on $\til{M}$.
Suppose $\lambda \subset \widetilde{S_0}$ lies outside
a large $N$-ball about $p$ in the (biLipschitz) hyperbolic metric on $\til{S_0}$. By Lemma 
\ref{adm-ea2-props} we obtain an electro-ambient
 quasigeodesic without backtracking
$\beta_{ea2}$ joining the end-points of $\lambda$ and
 lying outside an $M_0(N)$-ball about $p$ in $(\til{M},d_M)$, where $M_0(N) \rightarrow
 \infty $ as $N \rightarrow \infty $. 

Suppose that $\beta_{ea2}$ is
 a $(\kappa, \epsilon)$ electro-ambient quasigeodesic. Note that $\kappa, \epsilon$ depend on
 `the coarse Lipschitz constant' of $\Pi_\lambda$ and hence only on
 $\widetilde{S_0}$ and $\widetilde{M}$.

From  Lemma \ref{ea-strong} we know that
if $\beta^{h}$ denotes the (model) hyperbolic geodesic in
$\widetilde{M}$ joining the end-points of $\lambda$, then $\beta^h$
lies in a (uniform) $C^{\prime}$ neighborhood of $\beta_{ea2}$. 

Let $M_1(N) = M_0(N) - C^{\prime}$.
Then $M_1(N) \rightarrow
\infty$ as $N \rightarrow \infty$. Further, the (model)  hyperbolic geodesic 
 $\beta^h$ 
 lies outside an $M_1(N)$-ball around $p$. Hence, by Lemma
 \ref{contlemma}, 
the inclusion
  $i : \widetilde{S_0} \rightarrow \widetilde{M}$ extends continuously
  to a map 
  $\hat{i} : \widehat{S_0} \rightarrow \widehat{M}$. 

Since the
  continuous image of a compact locally connected set is locally
  connected  \cite{hock-young}
and the (intrinsic) boundary of $\widetilde{S_0}$ is a circle, we
  conclude that the limit set
  of $\widetilde{S_0}$ is locally connected.

This proves the theorem. \end{proof}

\section{Modifications for Surfaces with Punctures}

In this section, we shall describe the modifications necessary to prove
Theorem \ref{crucial}  for surfaces with punctures.

\subsection{ Partial Electrocution}
Two general references for this subsection are \cite{mahan-reeves},  \cite{mahan-pal}, where 
much
of what follows is done in a considerably more general
setting.

Let $M$ be a  convex hyperbolic 3-manifold
with a neighborhood of the cusps excised. Then each component of
the boundary of $M$ is of the form $\sigma \times P$, where $P$ is either
an interval or a circle, and $\sigma$ is a horocycle of some fixed 
length $e_0$. Each component of the boundary of the universal cover $\til M$, 
is a flat horosphere of the
form $\widetilde{\sigma} \times \tilde{P}$. Note that $\tilde{P} = P$ if
$P$ is an interval, and $\mathbb{R}$ if $P$ is a circle (the case for
a $(Z + Z)$-cusp).

The construction of {\it partially electrocuted}
horospheres below is half way between the spirit of
Farb's construction (in Lemmas \ref{farb1A}, \ref{farb2A},
where
the entire horosphere is coned off), and McMullen's Theorem \ref{ctm}
(where nothing is coned off, and properties of
 {\it ambient quasigeodesics} are investigated).

\smallskip

\noindent {\bf Partial Electrocution of Horospheres} \\
Let $Y$ be a convex simply connected hyperbolic 3-manifold.
Let $\mathcal{B}$ denote a collection of horoballs. Let $X$ denote
$Y$ minus the interior of the horoballs in $\mathcal{B}$. Let 
$\mathcal{H}$ denote the collection of boundary horospheres.Then each
$H_\alpha \in \mathcal{H}$ with the induced metric is isometric to a Euclidean
product $E^{1} \times L_\alpha$ for an interval $L_\alpha\subset \mathbb{R}$. Here $E^1$ denotes Euclidean $1$-space.

"Partially electrocute"  each 
$H_\alpha$ by giving it the product of the zero metric with the Euclidean metric,
i.e. on $E^{1}$ put the zero metric and on $L_\alpha$ put the Euclidean
metric.  Thus we are in the following situation:
\begin{enumerate}
\item $X$ is (strongly) hyperbolic relative to a collection $\mathcal{H}$ of horospheres.
\item Each horosphere $H_\alpha$ is equipped with a pseudometric making it isometric to a Euclidean
product $E^{1} \times L_\alpha$ for an interval $L_\alpha\subset \mathbb{R}$. We shall  denote the collection of 
$L_\alpha$'s by  $\LL$.
\end{enumerate}

The resulting pseudometric space is denoted $(X,d_{pel})$ and is called the {\bf partially electrocuted space} associated to the pair
$(X, \HH)$.

Its worth pointing out here that $(X,d_{pel})$ is essentially what one would get (in the spirit of \cite{farb-relhyp}) by gluing
to each $H_\alpha$ the mapping cylinder of the projection of $H_\alpha$ onto the $L_\alpha$-factor. Let $\mathcal G$ denote the collection of
these projections $g_\alpha : H_\alpha \rightarrow L_\alpha$. 
Thus, instead of coning all of a
horosphere down to a point 
we cone only horocyclic leaves of a foliation of the horosphere.
Effectively, therefore, we have a cone-line rather than a cone-point. We shall denote the union of $X$ and all the mapping cylinders of $g_\alpha$ by $\EE(X,\HH,\LL,\GG)$ in the spirit of the notation we have used
for electric spaces. As pointed out above, $\EE(X,\HH,\LL,\GG)$ and $(X,d_{pel})$ are quasi-isometric and both contain naturally embedded
copies of $X$ as a  subset (though not as a metric subspace). We shall therefore conflate $\EE(X,\HH,\LL,\GG)$ and $(X,d_{pel})$ in this subsection.
Geodesics and quasigeodesics
in the partially electrocuted space will be referred to as partially
electrocuted geodesics and quasigeodesics respectively.

In this
situation, we conclude as in Lemma \ref{farb1A}:

\begin{lemma} (Lemma 1.20 of \cite{mahan-pal})
For a 4-tuple $(X, \HH ,  \LL,\GG )$  as above,
$\EE(X,\HH,\LL,\GG)$ (resp. $(X,d_{pel})$) is a  hyperbolic metric spaces and  $L_\alpha \subset \EE(X,\HH,\LL,\GG)$
 (resp. $H_\alpha \subset (X,d_{pel})$) 
are uniformly quasiconvex.
\label{pel}
\end{lemma}

Recall that $X$ is obtained from a simply connected convex hyperbolic manifold $Y$ by excising a family of uniformly separated (open)
horoballs.

\begin{lemma}   (Lemma 1.21 of \cite{mahan-pal})
Let $(X, \HH ,  \LL,\GG )$  be a 4-tuple as above.
Given $K, \epsilon \geq 0$, there exists $C > 0$ such that the following
holds: \\
Let $\gamma_{pel}$ and $\gamma$ denote respectively a $(K, \epsilon)$
partially electrocuted quasigeodesic in $\EE(X,\HH,\LL,\GG)$ and a
$(K, \epsilon )$ hyperbolic quasigeodesic in $Y$ joining $a, b$. Then $\gamma \setminus
\bigcup_{H_\alpha\in\HH} H_\alpha$
lies in a  $C$-neighborhood of (any representative of)
$\gamma_{pel}$ in $(X,d)$. Further, outside of  the horoballs
that $\gamma$ meets, $\gamma$ and $\gamma_{pel}$ track each other, i.e. they lie in a $C$-neighborhood
of each other.
\label{pel-track}
\end{lemma}

\noindent {\bf Note:} $\EE(X,\HH,\LL,\GG)$ is strongly hyperbolic relative to
the sets $\{ L_\alpha \}$. In fact the space obtained by electrocuting the
sets $L_\alpha$ in $\EE(X,\HH,\LL,\GG)$ is just the space $\EE(X,\HH)$ obtained by 
electrocuting the sets $\{ H_\alpha \}$ in $X$.

Next, we show that partial electrocution preserves quasiconvexity. 

\begin{lemma} 
Given $C$ there exists $C_1$ such that if
 $A$ and $A \cap B$ (for any horoball
$B \in \mathcal{B}$) are $C$-quasiconvex in $Y$, then $(A\cap X,d_{pel})$ is $C_1$-quasiconvex 
in $(X,  d_{pel})$.

\label{pel-qc}
\end{lemma}

\begin{proof} It is given
that $A (\subset Y)$ as also $A \cap B$ for all $B \in \mathcal{B}$
are $C$-quasiconvex. Then given $a, b \in A \cap X$,
the  hyperbolic geodesic $\lambda$ in $Y$ joining $a, b$ lies 
in a $C$-neighborhood
of $A$. Since horoballs are convex, $\lambda$ cannot backtrack. We let $H = \partial B$ be the boundary horosphere of the horoball $B$,
and let $L$ be the element of $\LL$ corresponding to $H$.

Let $\lambda_{pel}$
be the partially electrocuted geodesic joining $a, b \in (X, d_{pel})$. Clearly, $\lambda_{pel}$ does not backtrack.  
Then by Lemma \ref{pel-track} above, we conclude that 
for all $H \in \mathcal{H}$
that $\lambda$ intersects, there exist points $a_H, b_H$ of $\lambda_{pel}$ close (in $Y$)
to the entry and exit points of $\lambda$ with respect to $H$. The points $a_H, b_H$
therefore lie close to
 $A \cap H$. Further,  the corresponding $L$ (resp. $H$) is quasiconvex in $\EE (X, \HH, \LL, \GG)$ (resp.
$(X, d_{pel})$)  by Lemma \ref{pel}.
It follows that $\lambda_{pel} \cap L $ (resp. $\lambda_{pel} \cap H$) lies within a uniformly bounded
distance of $A \cap H$ in $\EE (X, \HH, \LL, \GG)$ (resp.
$(X, d_{pel})$). The conclusion now
follows from Lemma \ref{pel-track}. \end{proof}

\subsection{Split geometry for Surfaces with Punctures} 
Recall that $N^h$ denotes (the convex core of) a simply or doubly degenerate hyperbolic 3-manifold {\it  with cusps}.  $N$ will 
 denote $N^h$ minus an open neighborhood of the cusps. $M$ will denote the model manifold
(Section \ref{min}) biLipschitz homeomorphic to $N$. {\it Since the proof in the case of surfaces with punctures is only
a small modification of the case of surfaces without punctures modulo known results (cf.
\cite{mahan-pal, mahan-pared}), we shall only sketch the proof, indicating
the necessary changes.}

It is worth noting here that the purpose of the partial electrocution operation
in the previous subsection is to ensure that successive split surfaces {\it with boundary} are uniformly separated so as
to ensure a model of weak split geometry
as defined in Remark \ref{wsplitrmk}. We shall proceed to construct a split geometry structure on $M$  outlined in the steps below. 
{\it In Steps (1)-(4) below we set up the model manifold of split geometry for $S$ with boundary.}

\noindent {\bf Step 1: Preliminary}\\
For a hyperbolic surface $S^h$ (possibly) with punctures, we fix a (small) 
$e_0$, and excise the cusps leaving horocyclic boundary components of
(ordinary or Euclidean) length 
$e_0$. We then take the induced {\em path metric} on $S^h$ minus cusps and call
the resulting surface $S$. This induced path metric will still be referred to
as the hyperbolic metric on $S$ (with the understanding that now $S$ possibly
has boundary). Note that the horocycle boundary components are now totally geodesic in $S$.

\noindent {\bf Step 2: Definition of Thick and Split Blocks and Hyperbolic Quasiconvexity of Split Components}\\
A thick block in $M$
is uniformly biLipschitz to $S \times I$ as before. 

The definitions and constructions of {\bf split building blocks} and 
{\bf split components} now go through with very little change. The
only difference is that $S$ now might have boundary curves of length $e_0$.

 There is one subtle point about hyperbolic quasiconvexity (in
$\til{M}$) of
split components. Hyperbolic quasiconvexity (cf. Lemma \ref{hypqc}) does not hold in the metric obtained by
merely excising the cusps and equipping the resulting horospheres with
the Euclidean metric. What we demand is that each split
component along with the parts of the horoballs that  abut it 
 be 
quasiconvex in $\tilde{N^h}$. Note that the intersection of split components in $\til M$ with horoballs
that abut it are (metric) products of horocycles with closed intervals. Lemma \ref{hypqc} furnishes the required quasiconvexity
in this case.

When we excise horoballs from $N^h$ to obtain $N$ and then partially electrocute horospheres in $N$
(or its biLipschitz model $M$)
in Step 3 below, and consider quasiconvexity in the resulting partially
electrocuted space, split components will remain
quasiconvex by Lemma \ref{pel-qc}. 

\noindent {\bf Step 3: Partially Electrocuting Horospherical Boundaries in $M$}\\
Next, we modify the metric on $M$ 
by partially electrocuting its boundary  horospherical  components so that the metric
on the horospherical boundary components of any (thick or split) block $S \times I$ is the product
of the zero metric on the horocycles of fixed (Euclidean) length $e_0$
and the Euclidean metric on the $I$-factor.
The resulting blocks will be called {\bf partially electrocuted
blocks}.  
Note that
$M_{pel}$ may also be constructed directly from $M$ by excising a
neighborhood of the cusps and partially electrocuting the resulting
horospheres. By Lemma \ref{pel} ${\til{M}}_{pel}$ is a hyperbolic
metric space and by  Lemma \ref{pel-qc}, partially electrocuted split components are quasiconvex in ${\til{M}}_{pel}$. 

\noindent {\bf Step 4: Split Blocks in ${\til{M}}_{pel}$ and Graph Quasiconvexity}\\
 Again,
the definitions and constructions of {\bf split blocks} and 
{\bf split components}  go through {\em mutatis mutandis} for the
partially electrocuted manifold ${\til{M}}_{pel}$.  By Lemma \ref{pel-qc},  quasiconvexity of split components as
well as   quasiconvexity of lifts of Margulis tubes is preserved by partial
electrocution.
 Hence in the model $M_{pel}$ obtained by gluing
together partially electrocuted blocks, the split components
are uniformly graph-quasiconvex. 

\smallskip

{\it In Steps (5)-(7) we indicate the modifications in the construction and use of the ladder $\LL_\lambda$ and the retract $\Pi_\lambda$.}\\

\noindent {\bf Step 5: Horo-ambient quasigeodesics}\\ Let $\lambda^h$ be a hyperbolic
geodesic in $\tilde{S^h}$. We replace pieces of
$\lambda^h$ that lie within horodisks by shortest
horocyclic segments joining its entry
and exit points (into the corresponding horodisk). Such a path
is called a {\bf horo-ambient quasigeodesic} cf. \cite{mahan-pared}. See
Figure below:

\begin{center}

\includegraphics[height=4cm]{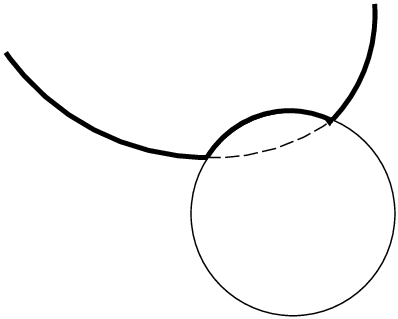}

\underline{Figure:{\it Horo-ambient quasigeodesic} }

\end{center}

A small modification might be introduced if we electrocute horocycles.
Geodesics and
quasigeodesics without backtracking
then travel for free along the zero metric horocycles.
This does not change matters much as the geodesics and
quasigeodesics in the two  constructions track each other by
Lemma \ref{farb2A}.

Thus, our starting point for the construction of the hyperbolic ladder
${\LL}_\lambda$ is not a hyperbolic geodesic $\lambda^h$ but a horoambient
quasigeodesic $\lambda$.

\noindent {\bf Step 6: Construction of the ladder ${\LL}_\lambda$}\\
The construction of ${\LL}_\lambda , \Pi_\lambda$ and their properties
go through {\it mutatis mutandis}
 and we conclude that ${\LL}_\lambda$ is quasiconvex
in the {\em graph metric} $({\til{M}}_{pel},d_G)$ on the
partially electrocuted space ${\til{M}}_{pel}$. As before,
${\widetilde{M_H}}_{pel}$ will denote the collection of horizontal sheets.

The modification of Theorem \ref{retract} in this context is given below:

\begin{theorem}
There exists $C > 0$ such that 
for any horo-ambient geodesic
 $\lambda = \lambda_0 \subset \widetilde{S} \times \{
0 \} \subset \widetilde{B_0}$, the retraction $\Pi_\lambda :
{\widetilde{M_H}}_{pel} \rightarrow {\LL}_\lambda $ satisfies: 
 $$d_{G}( \Pi_{\lambda} (x), \Pi_{\lambda} (y)) \leq C d_G(x,y) +
C.$$
\label{retract-punct}
\end{theorem}

\noindent {\bf Step 7: Decomposing the ladder ${\LL}_\lambda$ into $\LL^c_\lambda$ and ${\LL}^b_\lambda$}\\
From this step on, the modifications for punctured
surfaces follow \cite{mahan-pared}.
As in \cite{mahan-pared}, we decompose $\lambda$ into 
parts $\lambda^c$ and $\lambda^b$ consisting of (closures of) maximal segments  that lie
along horocycles and complementary pieces that do not  intersect 
 horocycles.  Accordingly, we decompose
${\LL}_\lambda$ into two
parts $\LL^c_\lambda$ and ${\LL}^b_\lambda$ consisting of parts that lie
along horocycles and those that do not. 
As in Lemma \ref{far-nopunct}, we get

\begin{lemma}
There exists a function $M(N) : \mathbb{N} \rightarrow \mathbb{N}$
such that  $M (N) \rightarrow \infty$
 as $N \rightarrow \infty$ for which the following holds:\\

For any horo-ambient quasigeodesic $\lambda \subset \widetilde{S}
\times \{ 0 \} \subset 
\widetilde{B_0}$, a fixed reference point 
 $p \in \widetilde{S} \times \{ 0 \} \subset
\widetilde{B_0}$ and any $x$ on  ${\LL}^b_\lambda$, 

\begin{center}

$d(\lambda^b , p) \geq N \Rightarrow d_{wel}(x,p) \geq M(N)$.

\end{center}

\label{far-punct}
\end{lemma}

 In Steps (8)-(10) we indicate the process of recovering a hyperbolic geodesic.\\
\noindent {\bf Step 8: Projecting and joining the dots}\\ 
Admissible
paths are constructed as in Section \ref{adm-sxn}. Now if $\lambda \subset \widetilde{S}
\times \{ 0 \} \subset 
\widetilde{B_0}$ is a horo-ambient geodesic joining $a, b$, let $\beta$ be an admissible path representing
 a $d_G$
geodesic in $\til{M_{pel}}$. Project $\beta\cap {\widetilde{M_H}}_{pel}$ onto ${\LL}_\lambda$ by $\Pi_\lambda$
 and "join the dots" as in Section \ref{join-sxn} to get
 a
connected ambient electric quasigeodesic $\beta_{amb}$.

\noindent {\bf Step 9: Recovery}\\ As in Sections \ref{recovery1} and \ref{recovery2}, construct from $\beta_{amb} \subset
 \widetilde{M}$ a partially electrocuted
quasigeodesic $\gamma$ in $({\widetilde{M}}_{pel}, d_{pel})$. 
Observe that  the parts of $\gamma$ that do not lie along partially 
electrocuted horospheres lie close
 to ${\LL}^b_\lambda$.
Hence, by Lemma \ref{far-punct} if $\lambda^h$ lies outside large
balls in $S^h$ then each point of  $\gamma\setminus \bigcup_{H_\alpha \in \HH} H_\alpha $   also lies outside
large balls in $\widetilde{M}$. 

At this stage we transfer the information to $\til N$ (=$\til{N^h}$ minus horoballs). Let $F: M \rightarrow N$ be the biLipschitz homeomorphism between $M$ and $N$ and
let $\til F$ denote its lift between universal covers. We thus conclude that  if $\lambda^h$ lies outside large
balls in $S^h$ then each point of  $\til{F}(\gamma\setminus \bigcup_{H_\alpha \in \HH} H_\alpha) $   also lies outside
large balls in $\widetilde{N}$. 

Note that in the case of surfaces without punctures, $\gamma$ itself
was a (biLipschitz) hyperbolic geodesic in $\til M$. However in the present situation of surfaces with punctures, one more step
of recovery is necessary.

\noindent {\bf Step 10: Conclusion}\\ Let $\gamma^h$ denote 
the hyperbolic geodesic in ${\widetilde{N}}^h$ joining the end-points
of $\til{F}(\gamma)$. By Lemma \ref{pel-track}  $\til{F}(\gamma)$ and $\gamma^h$
track each other away from   horoballs.  Then, every point of
$\gamma^h \cap \til{N}$ 
 must lie close to
 some point of $\til{F}(\gamma)$ lying outside partially 
electrocuted horospheres. Hence from Step (9), 
if $\lambda^h$ lies outside a large
ball about $p$ in $S^h$ then   $\gamma^h \cap  \til{N}$   also lies outside a
large ball about $p$ in $ \til{N}$. In particular,  $\gamma^h$ enters and
leaves horoballs at large distances from $p$. From this it follows (See Theorem 5.9 of \cite{mahan-pared} for instance)
that
$\gamma^h$ itself lies  outside a
large ball about $p$. Hence
by Lemma \ref{contlemma} there exists a Cannon-Thurston map and the
limit set is locally connected.
\\

We summarize the conclusion below:

\begin{theorem}
Let $N^h$ be a simply or doubly degenerate 3 manifold homeomorphic to $S^h \times J$ (for $J = [0,
  \infty ) $ or $( - \infty , \infty )$ respectively) for $S^h$ a finite volume hyperbolic surface such that  $i: S^h \rightarrow M^h$ is a proper map
inducing a homotopy equivalence. Then the inclusion
  $i : {\widetilde{S}}^h \rightarrow {\widetilde{N}}^h$ extends continuously
  to a map 
  $\hat{i} : {\widehat{S}}^h \rightarrow {\widehat{N}}^h$. Hence the limit set
  of ${\widetilde{S}}^h$ is locally connected.
\label{crucial-punct}
\end{theorem}

A part of the argument in  Lemmas \ref{adm-ea-props} and \ref{far-eas} and Step 9 above
does not use the full strength of the hypothesis that $M$ is a model for a surface group.
If we only assume that each end $E$ of a manifold $M$ is equipped with a split geometry structure where each split component
is incompressible, then the same arguments furnish the following.

\begin{lemma}   Let $N$ be the convex core of a complete hyperbolic $3-$manifold $N^h$ minus a neighborhood
of the cusps. Equip each degenerate end with a split geometry structure such that each split component
is incompressible. Let $M$ be the resulting model of split geometry and $F: N\rightarrow M$ be the bi-Lipschitz
homeomorphism between the two. Let $\til{F}$ be a lift of $F$ to the universal covers.
Then for all $C_0 > 0$, and $o \in \til{N}$
there exists a function $\Theta : \mathbb{N} \rightarrow \mathbb{N}$ satisfying
$\Theta (n) \rightarrow \infty$ as $n \rightarrow \infty$ such that the following holds.\\
For any $a, b \in \til{N}\subset \til{N^h}$, let $\lambda^h $ be the hyperbolic geodesic in $\til{N^h}$ joining
them and let $\lambda^h_{thick} = \lambda^h \cap  \til{N}$. Similarly let $\beta_{ea}^h$ be an electro-ambient
$C_0-$quasigeodesic without backtracking in $\til M \subset \EE(\til{M}, \KK^\prime )$
 joining $\til{F} (a)$, $\til{F} (b)$. Let
$\beta_{ea} = \beta_{ea}^h \setminus \partial \til{M}$ be the part of $\beta_{ea}^h$ lying away from
the (bi-Lipschitz) horospherical boundary of $ \til{M}$.

Then $d_M(\beta_{ea}, \til{F} (o)) \geq n $ implies that $d_{\Hyp^3}(\lambda^h_{thick}, o) \geq \Theta (n)$.
\label{contlemma2}
\end{lemma}

 This will be useful in \cite{mahan-kl}

\subsection{Local Connectivity of Connected Limit Sets}
Here we shall use a Theorem of Anderson and Maskit \cite{and-mask} along with Theorems \ref{crucial} and \ref{crucial-punct}
 above to prove that connected limit sets are locally connected. The connection between Theorems \ref{crucial} and \ref{crucial-punct} and Theorem \ref{lcfinal} below via Theorem \ref{and-mask} is similar to one discussed by Bowditch in \cite{bowditch-ct}. 

\begin{theorem} {\bf (Anderson-Maskit \cite{and-mask})} Let $\Gamma$ be an analytically finite Kleinian group with connected limit set. Then the limit set $\Lambda (\Gamma )$
is locally connected if and only if every simply  degenerate surface subgroup of $\Gamma$ without accidental parabolics has locally connected
limit set.
\label{and-mask}
\end{theorem}

 Combining Theorems \ref{crucial} and \ref{crucial-punct} with Theorem \ref{and-mask}, we have the following affirmative answer to
Question \ref{lcqn}.

\begin{theorem}
 Let $\Gamma$ be a finitely generated Kleinian group with connected limit set $\Lambda$. Then $\Lambda$ is locally connected.
\label{lcfinal}
\end{theorem}

Note that $\Lambda$ is connected if and only if the convex core of ${\mathbb{H}}^3/\Gamma$ is incompressible away from cusps.
In \cite{mahan-elct}, we prove that
for  surface groups without accidental parabolics, the point pre-images of the Cannon-Thurston map for points having multiple
pre-images
 are precisely the end-points of leaves of the ending lamination. In \cite{mahan-kl} we shall use the techniques developed in this paper to answer Question 
\ref{thurston-bams-qn} affirmatively.

\bibliography{ct}
\bibliographystyle{alpha}

\end{document}